\documentclass{article}
\usepackage[utf8]{inputenc}
\usepackage{mathtools} 
\usepackage{amssymb, amsthm} 
\usepackage{a4wide}
\usepackage{graphicx}
\usepackage[dvipsnames]{xcolor}

\usepackage{upgreek}

\usepackage{caption}
\usepackage{subcaption}

\usepackage{array,booktabs,multirow}

\usepackage{appendix}

\usepackage{todonotes}
\usepackage{xargs}

\usepackage{bbm}

\usepackage{algorithm}
\usepackage{algpseudocode}
\usepackage{enumitem}

\usepackage{hyperref}
\usepackage[nameinlink,capitalize]{cleveref}

\Crefname{assumption}{Assumption}{Assumption}
\crefname{assumption}{assumption}{assumption}

\newtheorem{lemma}{Lemma}
\Crefname{lemma}{Lemma}{Lemma}
\crefname{lemma}{lemma}{lemma}

\Crefname{corollary}{Corollary}{Corollary}
\crefname{corollary}{corollary}{corollary}
\newtheorem{remark}{Remark}
\Crefname{remark}{Remark}{Remark}
\crefname{remark}{remark}{remark}

\newtheorem{proposition}{Proposition}
\Crefname{proposition}{Proposition}{Proposition}
\crefname{proposition}{proposition}{proposition}

\Crefname{definition}{Definition}{Definition}
\crefname{definition}{definition}{definition}

\Crefname{example}{Example}{Example}
\crefname{example}{example}{example}

\DeclareMathAlphabet{\mathpzc}{OT1}{pzc}{m}{it}

\newcommand{\rme}{\mathrm{e}}

\newcommand{\rmd}{\mathrm{d}}
\newcommand{\rmO}{\mathrm{O}}
\newcommand{\rmF}{\mathrm{F}}
\newcommand{\rmG}{\mathrm{G}}
\newcommand{\rmI}{\mathrm{I}}

\newcommand{\bbR}{\mathbb{R}}

\newcommand{\bbE}{\mathbb{E}}
\newcommand{\bbT}{\mathbb{T}}
\newcommand{\bbP}{\mathbb{P}}

\newcommand{\calB}{\mathcal{B}}

\newcommand{\calN}{\mathcal{N}}

\newcommand{\calL}{\mathcal{L}}

\newcommand{\calC}{\mathcal{C}}
\newcommand{\calQ}{\mathcal{Q}}

\newcommand{\sfT}{\mathsf{T}}

\let\div\relax
\DeclareMathOperator{\div}{div}

\DeclareMathOperator{\Tr}{Tr}
\DeclareMathOperator{\Cov}{Cov}

\newcommand{\SV}{{\rm SV}}

\begin{document}

\title{Improving sampling by modifying the effective diffusion}
\author{T. Lelièvre\textsuperscript{1,2}, R. Santet\textsuperscript{1,2}, G. Stoltz\textsuperscript{1,2} \\
\small 1: CERMICS, ENPC, Institut Polytechnique de Paris, Marne-la-Vallée, France\\
\small 2: MATHERIALS project-team, Inria, Paris, France}

\maketitle

\begin{abstract}
    Markov chain Monte Carlo samplers based on discretizations of (overdamped) Langevin dynamics are commonly used in the Bayesian inference and computational statistical physics literature to estimate high-dimensional integrals. One can introduce a non-constant diffusion matrix to precondition these dynamics, and recent works have optimized it in order to improve the rate of convergence to stationarity by overcoming entropic and energy barriers. However, the introduced methodologies  to compute these optimal diffusions are generally not suited to high-dimensional settings, as they rely on costly optimization procedures. In this work, we propose to optimize over a class of diffusion matrices, based on one-dimensional collective variables (CVs), to help the dynamics explore the latent space defined by the CV. The form of the diffusion matrix is chosen in order to obtain an efficient effective diffusion in the latent space. We describe how this class of diffusion matrices can be constructed and learned during the simulation. We provide implementations of the Metropolis--Adjusted Langevin Algorithm and Riemann Manifold (Generalized) Hamiltonian Monte Carlo algorithms, and discuss numerical optimizations in the case when the CV depends only on a few degrees of freedom of the system. We illustrate the efficiency gains by computing mean transition durations between two metastable states of a dimer in a solvent.
\end{abstract}



\section{Introduction}

Computational statistical physics aims at estimating macroscopic quantities, such as the specific heat capacity of a material or the mean pressure of a system, by modelling matter at the microscopic scale and performing numerical simulations. These macroscopic quantities can be recast as thermodynamic averages, which are high-dimensional integrals. Molecular dynamics is then used to estimate these integrals, relying on trajectorial averages of well-chosen stochastic processes. These dynamics are actually used to sample high-dimensional probability distributions in many areas of applied mathematics, \emph{e.g.} for~Bayesian inference.

In this work, we focus on (overdamped) Langevin dynamics. Standard overdamped Langevin dynamics read
\begin{equation}
  \label{eq:standard_overdamped_langevin}
  \rmd q_t=-\nabla V(q_t)\rmd t+\sqrt{\frac{2}{\beta}}\rmd W_t,
\end{equation}
where~$V$ is the potential energy function defined on the position space~$\calQ$ of dimension~$d$,~$\beta=(k_BT)^{-1}$ is proportional to the inverse temperature ($k_B$ is the Boltzmann constant and~$T$ is the temperature of the heat bath) and~$(W_t)_{t\geqslant 0}$ is a standard~$d$-dimensional Brownian motion. Typically,~$\calQ=\bbR^{DN}$ or~$\calQ=\left( \ell\bbT\right)^{DN}$ (with~$\bbT$ the one-dimensional torus and~$\ell>0$) with~$DN=d$, where~$N$ is the number of particles and~$D\in\left\lbrace 1,2,3\right\rbrace$ is the dimension of the underlying space. Under standard assumptions on the potential energy function, the process~\eqref{eq:standard_overdamped_langevin} admits the Boltzmann-Gibbs measure~$\pi$ as invariant measure:
\begin{equation}
  \label{eq:gibbs_measure}
  \pi(\rmd q)=Z^{-1}\rme^{-\beta V(q)}\rmd q,\qquad Z=\int_{\calQ}\rme^{-\beta V}<+\infty.
\end{equation}
Discretizations of~\eqref{eq:standard_overdamped_langevin} therefore provide natural and efficient ways to sample Boltzmann-Gibbs measures~\eqref{eq:gibbs_measure}. One celebrated example is the Metropolis--Adjusted Langevin Algorithm (MALA)~\cite{rossky_1978,roberts_1996}, where the sampler is based on a Euler--Maruyama discretization of~\eqref{eq:standard_overdamped_langevin} along with a Metropolis--Hastings accept/reject procedure~\cite{metropolis_1953,hastings_1970} to remove the bias introduced by the time-discretization. This sampler achieves high performance by incorporating in the proposal the information of the derivatives of the log-density, which helps trajectories go towards low energy regions, and therefore reach regions of high-acceptance probability. 

\paragraph{Introducing a position-dependent diffusion in~\eqref{eq:standard_overdamped_langevin}.}
Overdamped Langevin dynamics can be generalized by introducing a diffusion operator~$D$ which outputs a positive definite symmetric matrix~$D(q)$ for any position~$q$:
\begin{equation}
  \label{eq:overdamped_langevin_diffusion}
  \rmd q_t = \left(-D(q_t)\nabla V(q_t) + \beta^{-1}\div D(q_t)\right)\rmd t + \sqrt{2\beta^{-1}}D(q_t)^{1/2}\rmd W_t.
\end{equation}
Here, the matrix~$D(q_t)^{1/2}$ is the square root of the diffusion matrix~$D(q_t)$ (defined by functional calculus) and~$\div D$ is defined as the vector whose~$i$-th component is the divergence of the~$i$-th column (or row) of the diffusion matrix~$D$. The dynamics~\eqref{eq:overdamped_langevin_diffusion} actually represent all the non-degenerate reversible diffusion processes that admit~$\pi$ as an invariant measure~\cite[Section~1.11.3, page~46]{bakry_2014}. Note that when the diffusion is set to the identity matrix for any position, the standard dynamics~\eqref{eq:standard_overdamped_langevin} are retrieved. 

In many practical cases of interest, and especially in high dimension, two characteristics of the target distribution impact the convergence of estimators: its multimodality (\emph{i.e.}~high probability regions are separated by low probability zones) and its anisotropy ({\em i.e.} high probability regions are elongated along some directions). This implies that isotropic local exploration methods, \emph{e.g.}~based on discretizations of~\eqref{eq:standard_overdamped_langevin}, generate metastable trajectories: the physical system stays in one conformation for an extensive period of time. One way to overcome these difficulties is to consider the dynamics~\eqref{eq:overdamped_langevin_diffusion} with a well chosen diffusion matrix~$D$.

To tackle anisotropy, note that the diffusion matrix~$D$ can be interpreted as the inverse of a position-dependent mass matrix, or as the inverse of a Riemannian metric. This diffusion operator is therefore a natural preconditioner for the overdamped Langevin dynamics. For strongly convex potential, one typical choice for~$D$ is the inverse of the Hessian of the potential~$V$, see {\em e.g.}~\cite{girolami_2011,graham_2022}. As for multimodality, we recently analyzed in~\cite{lelievre_2024} how to choose the diffusion~$D$ in simple one-dimensional settings. The aim of this paper is to generalize these results to high-dimensional settings.

\paragraph{Related works.}

Recent works have explored how to optimize the diffusion $D$ in order to accelerate convergence. In~\cite{cui_2024,lelievre_2024}, the optimizer is obtained by maximizing the spectral gap of the generator of~\eqref{eq:overdamped_langevin_diffusion}, which writes:
\begin{equation*}
  \calL = \left(-D\nabla V+\beta^{-1}\div D\right)\cdot\nabla + \beta^{-1}D:\nabla^{2}
\end{equation*}
with~$:$ the Frobenius inner product and~$\nabla^{2}$ the Hessian operator. The operator~$\calL$ is symmetric in the weighted space~$L^{2}(\pi)$ (this is the manifestation of the reversibility of the dynamics~\eqref{eq:overdamped_langevin_diffusion}):
\begin{equation*}
  \calL=-\beta^{-1}\nabla^{*}D\nabla,
\end{equation*}
where~$A^{*}$ denotes the adjoint of an operator~$A$ in~$L^{2}(\pi)$. Under appropriate assumptions on $V$, the operator $\calL$ has compact resolvent, and therefore a discrete, real, non-negative spectrum, the smallest eigenvalue being 0 with multiplicity 1. The spectral gap of~$\calL$, denoted by~$\Lambda$, is then given by the first nonzero eigenvalue. This quantity controls the convergence rate of the law of the process~\eqref{eq:overdamped_langevin_diffusion} at time~$t$, denoted by~$\pi_t$, towards the Boltzmann-Gibbs measure~\eqref{eq:gibbs_measure}: it holds~\cite{lelievre_2013,bakry_2014,lelievre_2016}
\begin{equation*}
\forall t \ge 0, \, \qquad
  \left\lVert \frac{\pi_t}{\pi}-1\right\rVert_{L^{2}(\pi)}\leqslant \rme^{-\beta^{-1}\Lambda t}\left\lVert\frac{\pi_0}{\pi}-1\right\rVert_{L^{2}(\pi)},
\end{equation*}
where we identified probability distributions with their associated densities. Larger spectral gaps therefore lead to faster convergence towards equilibrium. A natural idea followed in~\cite{cui_2024,lelievre_2024} is then to maximize~$\Lambda$ with respect to~$D$. In order for this maximization problem to be well-posed, normalization constraints need to be imposed on the diffusion. Indeed, the spectral gap associated with the diffusion~$\alpha D$ is~$\alpha\Lambda$ for any~$\alpha>0$: multiplying the diffusion by a  positive constant multiplies the rate of convergence to equilibrium. In order to maintain the level of numerical error, one however needs to divide the time step by a factor~$\alpha$, so that, at the numerical level, no improvement is actually observed. In~\cite{cui_2024}, $\calQ=\bbR^{d}$ and the constraint is chosen to be
\begin{equation*}
  \int_{\calQ}\Tr D(q)\pi(\rmd q) = \Tr\left( \Cov\pi\right),\qquad \Cov\pi=\bbE_{\pi}\left[\left(X-\bbE_{\pi}\left[X\right]\right)\left(X-\bbE_{\pi}\left[X\right]\right)^{\sfT}\right]\in\bbR^{d\times d}.
\end{equation*}
In that case, the optimizer can be related to the so-called Stein kernels. In~\cite{lelievre_2024}, $\calQ=\bbT^{d}$ and the normalization constraint is chosen to be a~$L^{p}$ norm:
\begin{equation}
  \label{eq:constraint_lp_norm}
  \int_{\calQ}\left\lvert \rme^{-\beta V(q)}D(q)\right\rvert_{\rmF}^{p}\,dq = 1,
\end{equation}
where~$\left\lvert\cdot\right\rvert_{\rmF}$ denotes the Frobenius norm. Moreover, an analytical expression of an optimizer can be obtained in an homogenized regime (in one dimension for the~$L^{p}$ constraint~\eqref{eq:constraint_lp_norm} and any dimension for linear constraints, see~\cite[Section~5]{lelievre_2024}), which reads
\begin{equation}
\label{eq:optimal_diffusion_homogenization}
  D_{\rm Hom}(q)=\rme^{\beta V(q)}\rmI_d.
\end{equation}
This diffusion has also been proposed in~\cite{roberts_2002} to efficiently sample Boltzmann-Gibbs measures. This optimizer seems to be a good approximation of the optimal diffusion at least in one-dimensional cases~\cite{lelievre_2024}.

In high-dimensional settings, solving numerically the maximization of the spectral gap of $\calL$ over the diffusion is not an easy task. For example, Finite Element Methods are used in~\cite{lelievre_2024}, but these do not scale well with dimension. Using the optimal diffusion in the homogenized regime~\eqref{eq:optimal_diffusion_homogenization} is not a solution either, since the variations of the potential energy~$V$ are prohibitively large in high dimensions (since~$V$ is an extensive quantity), and this requires to drastically reduce the time step in order to maintain the accuracy. 

\paragraph{Constructing the diffusion.}

In this work, we propose to use a diffusion whose analytical expression is based on the optimal homogenized diffusion~\eqref{eq:optimal_diffusion_homogenization}, but replacing the potential energy~$V$ by an \textit{effective} potential energy in a low-dimensional space, namely the free energy~$F$ associated with a collective variable~$\xi:\calQ\to\bbR^{m}$ (with $m \ll d$). Collective variables are maps commonly used in computational statistical physics to build a low-dimensional representation of the molecular system. Usually,~$\xi(q)$ is defined to be a slow variable of the system, \emph{i.e.}~the characteristic evolution time of~$\xi(q_t)$, where~$q_t$ solves for example~\eqref{eq:standard_overdamped_langevin} or~\eqref{eq:overdamped_langevin_diffusion}, is much larger than the characteristic time needed for the diffusion process to sample level sets of~$\xi$. In this work, we focus on scalar-valued collective variables ($m=1$), the extension to the multivalued case being left for future works, see Remark~\ref{rem:multidimendisional_case}.

To properly introduce the class of diffusions considered in this work, let us first introduce various objects. For~$z\in\xi(\calQ)$, denote by~$\Sigma(z)=\xi^{-1}\left(\left\lbrace z\right\rbrace\right)$ the level set of the collective variable for the level~$z$. The probability measure~$\pi$ conditioned at a fixed value~$z$ of the collective variable is given by (see \emph{e.g}~\cite{lelievre_2010})
\begin{equation}
  \label{eq:conditional_measure}
  \pi^{\xi}(\rmd q|z)=\frac{\rme^{-\beta V(q)}\left\lVert\nabla\xi(q)\right\rVert^{-1}\sigma_{\Sigma(z)}(\rmd q)}{\displaystyle\int_{\Sigma(z)}\rme^{-\beta V}\left\lVert\nabla\xi\right\rVert^{-1}\rmd\sigma_{\Sigma(z)}},
\end{equation}
where the measure~$\sigma_{\Sigma(z)}(\rmd q)$ is the surface measure on~$\Sigma(z)$ (namely the measure on~$\Sigma(z)$ induced by the Lebesgue measure in the ambient space~$\calQ$ and the Euclidean scalar product), and~$\left\lVert\cdot\right\rVert$ denotes the Euclidean norm on~$\bbR^{d}$. Note that the measure~$\left\lVert\nabla\xi\right\rVert^{-1} \sigma_{\Sigma(z)}(\rmd q)$ is sometimes called the \emph{delta measure} and denoted by~$\delta_{\xi(q)-z}(\rmd q)$. The free energy~$F$ is given by
\begin{equation}
  \label{eq:free_energy}
  F(z)=-\beta^{-1}\ln\left(\int_{\Sigma(z)}Z^{-1}\rme^{-\beta V}\left\lVert\nabla\xi\right\rVert^{-1}\rmd\sigma_{\Sigma(z)}\right),
\end{equation}
where~$Z$ is the normalizing constant defined in~\eqref{eq:gibbs_measure}. For any~$q\in\calQ$, let us define the matrix
\begin{equation}
  \label{eq:P}
  P(q)=\frac{\nabla\xi(q)\otimes\nabla\xi(q)}{\left\lVert\nabla\xi(q)\right\rVert^{2}}\in\bbR^{d\times d}.
\end{equation}
For~$q\in\Sigma(z)$, the orthogonal projection operator~$P^{\perp}$ onto the tangent space~$T_q\Sigma(z)$ to~$\Sigma(z)$ at~$q$ is given by (see for instance~\cite[Section~3.2.3.1]{lelievre_2010})
\begin{equation}
  \label{eq:P_perp}
  P^{\perp}(q)=\rmI_d-P(q)=\rmI_d-\frac{\nabla\xi(q)\otimes\nabla\xi(q)}{\left\lVert\nabla\xi(q)\right\rVert^{2}}.
\end{equation}
Note that
\begin{equation}
  \label{eq:P_nabla_xi_relations}
  P^{\perp}(q)\nabla\xi(q)=0,\qquad P(q)\nabla\xi(q)=\nabla\xi(q).
\end{equation}
Lastly, let us introduce
\begin{equation}
    \label{eq:intro_effective_diffusion}
  \sigma^{2}(z)=\int_{\Sigma(z)}\left\lVert\nabla\xi\right\rVert^{2}\rmd\pi^{\xi}(\cdot|z).
\end{equation}
The map~$\sigma$ can be seen as the (multiplicative) noise of the effective dynamics on the latent space~$\xi(\calQ)$, see Section~\ref{sec:effective_dynamics} for further details. The class of diffusions we are interested in is then defined as
\begin{equation}
  \label{eq:optimal_diffusion_class}
  D_{\alpha}(q)=\kappa_{\alpha}\left[P^{\perp}(q)+a_{\alpha}(\xi(q))P(q)\right]=\kappa_{\alpha}\left[\rmI_d+\left(a_{\alpha}(\xi(q))-1\right)P(q)\right],\qquad a_{\alpha}(z)=\frac{\rme^{\alpha\beta F(z)}}{\sigma^{2}(z)},
\end{equation}
where~$\alpha\in\bbR$ is a numerical parameter and~$\kappa_{\alpha}>0$ is a normalization constant. Note that~$D_{\alpha}$ indeed defines a diffusion, in the sense that it has values in the set of symmetric positive definite matrices. Indeed, the eigenvalues of~$D_{\alpha}(q)$ are~$\kappa_{\alpha}>0$ and~$\kappa_{\alpha}(1+a_{\alpha}(\xi(q))-1)=\kappa_{\alpha}a_{\alpha}(\xi(q))>0$ for any~$\alpha\in\bbR$.

\begin{figure}
    \centering
    \begin{tikzpicture}[xscale=0.8, yscale=0.8]
        \node at (0,0) {\includegraphics[draft=false,width=0.78\columnwidth]{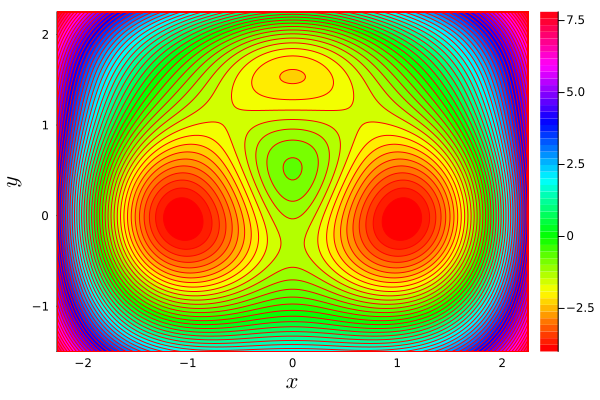}};
        \draw[->, thick, blue, line width=1mm] (-3,-0.4) --++ (3.5,0);
        \draw[->, line width=0.7mm] (-3,-0.4) node {\Large$\cdot$} node [left] {\Large$\mathbf{q}$} --++ (1.5,0) node [below] {\Large$\mathbf{P(q)}$};
        \draw[->, line width=0.7mm] (-3,-0.4) --++ (0, 1.5) node [above right, yshift=-1mm] {\Large $\mathbf{P}^\perp(q)$};
        \draw[dashed, line width=0.7mm] (-3,4.45) -- (-3,-3.65) node [below, yshift=-4mm] {\Large$\mathbf{\Sigma(\xi(q))}$};
        \draw[->, thick, blue, line width=1mm] (2.35,-4.75) --++ (1.5,0) node [right] {\Large$\left(a_\alpha(\xi(q))-1\right)\mathbf{P(q)}$};
        \node[rotate=-90] at (7.7, 0.25) {\Large Potential Energy};
    \end{tikzpicture}
    \caption{Graphical illustration of the diffusion~\eqref{eq:optimal_diffusion_class}. The potential energy function is symmetric with respect to the $x$-axis and exhibits two deep wells separated by an energy barrier. An example of a collective variable is~$\xi(q)\equiv \xi(x,y)=x$. The function~$a_\alpha$ modifies the amplitude of the diffusion along the direction of~$\nabla\xi(q)$, which is orthogonal to the level set~$\Sigma(\xi(q))$. This helps the particle at point~$q$ to visit new values of the $x$-coordinate.}
    \label{fig:diffusion_graphical_motivation}
\end{figure}

Let us give a brief description and motivation for the diffusion~\eqref{eq:optimal_diffusion_class}. The crucial point in our construction is to modulate the diffusion in the directions which are ``difficult to explore'' (those parametrized by~$\xi$), with some factor~$a_\alpha-1$; while leaving the diffusion unchanged in the associated orthogonal directions. The modulation factor itself is a function of~$\xi$ only. Concretely, this is done by modulating only the $P$ part of the diffusion. We refer to Figure~\ref{fig:diffusion_graphical_motivation} for a graphical illustration of this idea. The analytical expression of the map~$a_{\alpha}$ is based on the expression of the optimal homogenized diffusion~\eqref{eq:optimal_diffusion_homogenization}. It is built such that the effective dynamics in the latent space~$\xi(\calQ)$ is governed by the optimal homogenized effective diffusion~$\rme^{\beta F}$ (at least when~$\alpha=1$); see Section~\ref{sec:effective_dynamics} for further precisions. Using this diffusion should therefore favor exploration in the latent space which, if the collective variable is well-chosen, should accelerate the convergence towards equilibrium. The factor~$\alpha$ present in the definition of~$a_{\alpha}$ is introduced to ensure that~$\alpha\beta F$ is of order~1 in order to scale properly the argument in the exponential and avoid too large diffusions which would require to decrease the time step used in the numerical integration. Note that when~$\alpha=0$ and~$\sigma^{2}$ is constant equal to~1, the standard overdamped Langevin dynamics~\eqref{eq:standard_overdamped_langevin} are retrieved. Lastly, the constant~$\kappa_{\alpha}$ acts as (the inverse of) a normalization constant for the diffusion. When discretizing the dynamics with time step~$\Delta t$, this constant simply scales the time step. Introducing the constant~$\kappa_\alpha$ is beneficial when learning the diffusion on the fly, as described in Section~\ref{sec:adaptive_learning} below. When the diffusion does not need to be learned, one can simply ignore the constant~$\kappa_\alpha$ and tune the time step~$\Delta t$, \emph{e.g.}~to meet a specific rejection probability in Metropolis algorithms.

To construct the diffusion~\eqref{eq:optimal_diffusion_class}, we need two quantities that we do not know \textit{a priori}: the free energy~$F$ defined in~\eqref{eq:free_energy} and the effective diffusion~$\sigma^{2}$ defined in~\eqref{eq:intro_effective_diffusion}. These quantities are both averages with respect to conditional measures. Therefore, they can be computed using standard techniques to estimate conditional expectations. One can for example use constrained sampling methods, as in the estimation of mean force and free energy for thermodynamic integration (see Appendix~\ref{app:TI} and~\cite{kirkwood_1935,lelievre_2010}). These conditional expectations can also be learned on the fly as in free energy adaptive biasing techniques: this will be presented in Section~\ref{sec:adaptive_learning} where we use Adapted Biasing Force (ABF) methods~\cite{darve_2001,henin_2004}.

We will show in Section~\ref{sec:MALA_numerics} that the quantities required to run the sampling dynamics~\eqref{eq:overdamped_langevin_diffusion} with the diffusion~\eqref{eq:optimal_diffusion_class}
are easily obtained (in particular its square root, inverse and determinant). Moreover, when the collective variable~$\xi$ is a function of only~$k\leqslant d$ components of~$q$ (such as for bond lengths or dihedral angles), then only those~$k$ components are modified by the diffusion. The implementation of samplers based on discretizations of overdamped Langevin dynamics can therefore be tailored to limit computational overheads, so that the associated computational costs are typically negligible compared to force computations, see Section~\ref{sec:MALA_numerics} for details.

Finally, note that this class of diffusion can also be used in samplers related to the Langevin dynamics, such as (Generalized) Hamiltonian Monte Carlo ((G)HMC) algorithms~\cite{duane_1987,horowitz_1991}. These algorithms add a momentum variable, generating a Markov chain in the phase space~$\calQ\times\bbR^{d}$. The Markov chain is constructed such that the probability measure~$\mu(\rmd q\,\rmd p)\propto\rme^{-\beta H(q,p)}\rmd q\,\rmd p$ is an invariant measure for the Markov chain. The Hamiltonian function~$H:\calQ\times\bbR^{d}\to\bbR$ is defined so that the marginal in position of~$\mu$ is exactly~$\pi$. To introduce nonconstant diffusions in (G)HMC algorithms, the relevant framework is based on Riemann Manifold (Generalized) Hamiltonian Monte Carlo algorithms (RM(G)HMC)~\cite{girolami_2011}. In that case, the diffusion acts as the inverse of a position-dependent mass tensor, which preconditions the Hamiltonian dynamics~\cite{beskos_2011,bou-rabee_2018,bou-rabee_2021}. For good choices of the numerical parameters, these algorithms actually provide weakly consistent discretizations of the overdamped Langevin dynamics~\eqref{eq:overdamped_langevin_diffusion} (see~\cite[Section~3.3]{lelievre_2023_i}). Since the Hamiltonian function in the RM(G)HMC algorithms is not separable, implicit problems have to be solved (\emph{e.g.}~using Newton's method) and reversibility checks have to be implemented in order to perform an unbiased sampling~\cite{graham_2022,noble_2023,lelievre_2023_i}. In the case when the collective variable~$\xi$ is a function of a small number of components of~$q$, one can optimize the implementation of Newton's method in order to limit the computational costs (more details are given in Section~\ref{sec:RMHMC}). For our numerical illustration, the best results are actually obtained with RMGHMC algorithms.



\paragraph{Outline of the work.} In Section~\ref{sec:overdamped}, we motivate the choice of the analytical expression of the diffusion~\eqref{eq:optimal_diffusion_class} and detail how samplers based on discretizations of the overdamped Langevin dynamics~\eqref{eq:overdamped_langevin_diffusion} can be implemented. We also present an adaptive scheme to update the diffusion when the free energy and effective diffusion are not available when the simulation starts. In Section~\ref{sec:HMC}, we describe how the diffusion can be introduced in samplers based on discretizations of the Langevin dynamics, utilizing the RMHMC algorithm and its generalized variant. All numerical results are presented for the same physical system composed of a dimer in a solvent, as described in Section~\ref{sec:MALA_numerics}. All the methods and experiments are provided in an open source Julia code available at~\url{https://github.com/rsantet/Improving_Sampling_By_Modifying_The_Effective_Diffusion}.

\paragraph{Assumptions on the collective variable.} We assume that the collective variable is smooth (at least~$\calC^{2}$) and that the gradient of~$\xi$ is nonzero everywhere. It should be noted that the gradient and Hessian of the collective variable are needed in order to run the algorithms presented in this work. These quantities can be either derived by hand or obtained numerically using automatic differentiation tools. 


\section{Optimizing the diffusion: the overdamped Langevin case}
\label{sec:overdamped}

In this section, we make precise how the diffusion~\eqref{eq:optimal_diffusion_class} can be used in combination with samplers based on discretizations of the overdamped Langevin dynamics~\eqref{eq:overdamped_langevin_diffusion}. In Section~\ref{sec:effective_dynamics}, we provide a theoretical motivation for the diffusion~\eqref{eq:optimal_diffusion_class} by computing the effective dynamics associated with the overdamped Langevin dynamics~\eqref{eq:overdamped_langevin_diffusion}. We then provide in Section~\ref{sec:MALA} one possible implementation of a sampler based on discretizations of the overdamped Langevin dynamics using the MALA algorithm. Associated numerical results are presented in Section~\ref{sec:MALA_numerics}. We next describe in Section~\ref{sec:adaptive_learning} a methodology to learn the diffusion along the simulation, utilizing standard methods used in ABF algorithms. Associated numerical results are presented in Section~\ref{sec:adaptive_MALA_numerics}.

\subsection{Effective dynamics}
\label{sec:effective_dynamics}

The central motivation for the choice of the analytical expression of the diffusion~\eqref{eq:optimal_diffusion_class} is that the associated effective dynamics in the latent space is governed by the optimal homogenized diffusion~$\rme^{\beta F}$ (when~$\alpha=1$). Let us recall that for well-chosen collective variables, effective dynamics are good approximations of the dynamics~$t\mapsto \xi(q_t)$ in the latent space~$\xi(\calQ)$ (see~\cite{legoll_2010}). In particular, whatever the map~$\xi$, their stationary probability measure is the image of the measure~$\pi$ by~$\xi$, denoted by~$\xi\star\pi(\rmd z)$ and defined by
\begin{equation}
  \label{eq:image_measure}
  \xi\star\pi(\rmd z)=\rme^{-\beta F(z)}\rmd z=Z^{-1}\left(\int_{\Sigma(z)}\rme^{-\beta V}\left\lVert\nabla\xi\right\rVert^{-1}\rmd\sigma_{\Sigma(z)}\right)\rmd z.
\end{equation}
The effective dynamics obtained for the standard overdamped Langevin dynamics~\eqref{eq:standard_overdamped_langevin} is recalled in the next proposition (see~\cite[Section~2.3]{legoll_2010} for the derivation; recall that, in our work, the collective variable~$\xi$ is a scalar-valued function).
\begin{proposition}
    \label{prop:standard_effective_dynamics}
  Let~$q_t$ solve~\eqref{eq:standard_overdamped_langevin}. Then the process~$t\mapsto\xi(q_t)$ satisfies
  \begin{equation*}
    \rmd\xi(q_t)=\left(-\nabla V(q_t)\cdot\nabla\xi(q_t)+\beta^{-1}\Delta\xi(q_t)\right)\rmd t+\sqrt{2\beta^{-1}}\nabla\xi(q_t)\cdot\rmd W_t.
  \end{equation*}
  The effective dynamics is defined by
  \begin{equation}
    \label{eq:standard_overdamped_langevin_dynamics_effective_dynamics}
    \rmd z_t=b(z_t)\,\rmd t+\sqrt{2\beta^{-1}}\sigma(z_t)\,\rmd B_t,
  \end{equation}
  where
  \begin{equation}
    \label{eq:standard_effective_drift_noise}
    b(z)=\int_{\Sigma(z)}\left(-\nabla V\cdot\nabla\xi+\beta^{-1}\Delta \xi\right)\rmd\pi^{\xi}(\cdot|z),\qquad \sigma^{2}(z)=\int_{\Sigma(z)}\left\lVert\nabla\xi\right\rVert^{2}\rmd\pi^{\xi}(\cdot|z),
  \end{equation}
  and~$(B_t)_{t\geqslant0}$ is a standard one-dimensional Brownian motion. Furthermore, it holds
  \begin{equation}
    \label{eq:standard_drift_noise_relation_effective_dynamics}
    b(z)=-\sigma^{2}(z)F'(z)+\beta^{-1}\left(\sigma^{2}\right)'(z),
  \end{equation}
  so that the effective dynamics can be rewritten as a one-dimensional overdamped Langevin dynamics of the form~\eqref{eq:overdamped_langevin_diffusion} with effective potential~$F$ and effective diffusion~$\sigma^{2}$.
\end{proposition}
In~\eqref{eq:standard_drift_noise_relation_effective_dynamics}, the derivative of the free energy, called the \emph{mean force}, is given by
\begin{equation}
  \label{eq:mean_force}
  F'(z)=\int_{\Sigma(z)}f\rmd\pi^{\xi}(\cdot|z),\qquad f=\frac{\nabla V\cdot\nabla\xi}{\left\lVert\nabla\xi\right\rVert^{2}}-\beta^{-1}\div\left(\frac{\nabla\xi}{\left\lVert\nabla\xi\right\rVert^{2}}\right).
\end{equation}
The map~$f$ is called the \emph{local mean force}. It follows from~\eqref{eq:standard_drift_noise_relation_effective_dynamics} that the effective dynamics~\eqref{eq:standard_overdamped_langevin_dynamics_effective_dynamics} admits~$\xi\star\pi$ as a stationary probability measure. The proof of the identity~\eqref{eq:standard_drift_noise_relation_effective_dynamics} is given in Appendix~\ref{app:proposition_1} (see also~\cite[Lemma~2.4]{legoll_2010} for similar computations).

The effective dynamics~\eqref{eq:standard_overdamped_langevin_dynamics_effective_dynamics} is governed by the effective diffusion~$\sigma^{2}$, which may not be optimal in order to favor exploration in the latent space~$\xi(\calQ)$. It therefore makes sense to modify the diffusion of the original dynamics and thus this effective diffusion in order to obtain better convergence towards equilibrium for the effective dynamics. As the following result shows, one possible solution to modify the effective diffusion is to change the (original) diffusion only in the direction of~$\nabla\xi$. The proof is given in Appendix~\ref{app:proof_effective_dynamics}.
\begin{proposition}
  \label{prop:effective_dynamics}
  Let~$q_t$ solve~\eqref{eq:overdamped_langevin_diffusion} with diffusion~$D(q)=P^{\perp}(q)+a(\xi(q))P(q)$ where~$a:\xi\left(\calQ\right)\to(0,+\infty)$. Then the process~$t\mapsto\xi(q_t)$ satisfies
  \begin{align}
    \label{eq:overdamped_langevin_dynamics_effective_dynamics}
    \rmd\xi(q_t)&=\left(
      -a(\xi(q_t))\left(
        \nabla V(q_t)\cdot\nabla\xi(q_t)-\beta^{-1}\Delta\xi(q_t)
      \right)+\beta^{-1}a'(\xi(q_t))\left\lVert\nabla\xi(q_t)\right\rVert^{2}
    \right)\rmd t\\
    &\nonumber\quad+\sqrt{2\beta^{-1}a(\xi(q_t))}\nabla\xi(q_t)\cdot\rmd W_t.
  \end{align}
  The effective dynamics is then given by
  \begin{equation}
    \label{eq:effective_dynamics}
    \rmd z_t=b_a(z_t)\,\rmd t+\sqrt{2\beta^{-1}}\sigma_{a}(z_t)\,\rmd B_t,
  \end{equation}
  where
  \begin{equation}
    \label{eq:drift_noise_effective_dynamics}
    \left\lbrace
    \begin{aligned}
      b_a(z)&=a(z)b(z)+\beta^{-1}a'(z)\sigma^{2}(z),\\
      \sigma_a^{2}(z)&=a(z)\sigma^{2}(z),
    \end{aligned}
    \right.
  \end{equation}
  and~$(B_t)_{t\geqslant0}$ is a standard one-dimensional Brownian motion. As in~\eqref{eq:standard_drift_noise_relation_effective_dynamics}, the maps~$b_a$ and~$\sigma_a$ are related by the following identity:
  \begin{equation}
    \label{eq:drift_noise_relation_effective_dynamics}
    b_a(z)=-\sigma_a^{2}(z)F'(z)+\beta^{-1}\left(\sigma_a^{2}\right)'(z).
  \end{equation}
\end{proposition}
The identity~\eqref{eq:drift_noise_relation_effective_dynamics} shows that the probability measure~$\xi\star\pi$ is still an invariant probability measure for the effective dynamics~\eqref{eq:overdamped_langevin_dynamics_effective_dynamics}, whatever the choice of the map~$a$. A natural candidate for the effective diffusion~$\sigma_a$ in~\eqref{eq:effective_dynamics}-\eqref{eq:drift_noise_relation_effective_dynamics} is the optimal homogenized diffusion~$\rme^{\beta F}$, which was shown to be optimal in some sense in~\cite[Section~5]{lelievre_2024}. In practice, one may want to scale the free energy by a factor~$\alpha\in\bbR$ so that the argument of the exponential is not too large, thus preventing unstable dynamics. This leads us to introduce the diffusion~\eqref{eq:optimal_diffusion_class}. We demonstrate on a numerical example below that using this diffusion indeed leads to a faster exploration of the latent space than with the standard overdamped Langevin dynamics~\eqref{eq:standard_overdamped_langevin}, hence a better exploration of the potential energy surface.


\begin{remark}
  \label{rem:free_energy_constant}
    Note that any diffusion of the form $D(q)=g(q)P^\perp(q)+a(\xi(q))P(q)$ with $g$ a general function of the position leads to the same evolution process~\eqref{eq:overdamped_langevin_dynamics_effective_dynamics}, thus to the same effective dynamics~\eqref{eq:effective_dynamics}. There is therefore a degree of freedom corresponding to the choice of this function $g$, which amounts to balancing the relative (pointwise) contributions of~$P$ and~$P^\perp$ in~$D$. In particular, multiplying the function $a \circ \xi$ by a constant in~\eqref{eq:optimal_diffusion_class}, which amounts to adding a constant to the free energy, has an effect on this contribution. Optimally balancing the weights between $P$ and $P^\perp$ is left for future work.
\end{remark}

\begin{remark}[The multidimensional case]
\label{rem:multidimendisional_case}
    Given a multivariate collective variable~$\xi:\bbR^d\to\bbR^m$, the effective diffusion writes~$\bbE_{\pi^\xi}\left[\nabla\xi^\top D \nabla\xi\right]$, with~$\nabla\xi\in\bbR^{d\times m}$ and~$D$ the diffusion of the original overdamped Langevin dynamics~\cite{zhang_2016}. Given a target symmetric positive definite matrix~$\Sigma\in\bbR^{m\times m}$, the diffusion~$D(q)=P^\perp(q)+P(q)A(q)P(q)=\rmI_d-P(q)+A(q)$ with
    \begin{equation*}
        \left\lbrace
        \begin{aligned}
          P(q)&=\sum_{k,\ell=1}^{m}G_{k,\ell}^{-1}(q)\nabla\xi_{k}^\sfT(q)\nabla\xi_{\ell}(q),\quad G(q)=\nabla\xi^\sfT(q)\nabla\xi(q),\quad G_{k,\ell}^{-1}(q):=[G(q)^{-1}]_{k,\ell},\\
          A(q)&=\sum_{k,k',\ell,\ell'=1}^{m}G_{k,k'}^{-1}(q)\Sigma_{k,\ell}(\xi(q))G_{\ell,\ell'}^{-1}(q)\nabla\xi_{k'}(q)\nabla\xi_{\ell'}^{\sfT}(q),
        \end{aligned}
        \right.
    \end{equation*}
    is symmetric positive definite and such that the diffusion governing the effective dynamics is~$\Sigma$. Here~$P$ is the generalization of the projection operator~\eqref{eq:P} in the multidimensional case. Choosing a `good` target effective diffusion~$\Sigma$ remains an open question when~$m\geqslant 2$, and is left for future endeavours. 
    
    In the case~$m=1$, this diffusion writes~$P^\perp(q)+\frac{\Sigma(\xi(q))}{\left\lVert\nabla\xi(q)\right\rVert^2}P(q)$, which is similar to the choice we make in this work except that we enforce the denominator~$\left\lVert\nabla\xi(q)\right\rVert^2$ to already be averaged under~$\pi^\xi$ for the original overdamped Langevin dynamics.
\end{remark}

\subsection{Implementation using MALA}
\label{sec:MALA}

In this section, we show how to build samplers based on discretizations of the overdamped Langevin dynamics, using MALA as an example. Note that a discretization for Langevin processes with position-dependent diffusions has recently been proposed in~\cite{bronasco_2025}, which improves on the weak error scaling with respect to the time step compared to MALA. We also provide at the end of the section the analytical expressions of all the quantities needed to run the algorithm, such as the inverse, determinant, and divergence of the diffusion~\eqref{eq:optimal_diffusion_class}.

MALA is built on two blocks: (i) a proposal computed using an Euler--Maruyama discretization of~\eqref{eq:overdamped_langevin_diffusion} with time step~$\Delta t$ and (ii) a Metropolis--Hastings accept/reject procedure~\cite{metropolis_1953,hastings_1970}.  Even though the name MALA is usually referring to the case~$D \equiv \rmI_d$, we also use this name when using a nonconstant diffusion. For a fixed configuration~$q^{n}\in\calQ$, the proposal is defined by
\begin{equation}
  \label{eq:MALA_proposal}
  \widetilde{q}^{n+1}=q^{n}+\left[-D(q^{n})\nabla V(q^{n})+\beta^{-1}\div D(q^{n})\right]\Delta t + \sqrt{2\beta^{-1}\Delta t}\,D(q^{n})^{1/2}\rmG^{n+1},
\end{equation}
where~$\rmG^{n+1}\sim\calN(0,\rmI_{d})$. The associated transition kernel~$T$ has the following density:
\begin{equation*}
  T(q,q')=\left(\frac{\beta}{4\pi\Delta t}\right)^{d/2}\det\left(D(q)\right)^{-1/2}\exp\left(
  -\frac{\beta}{4\Delta t}\left(q'-\mu_{\Delta t}(q)\right)^{\sfT}D(q)^{-1}\left(q'-\mu_{\Delta t}(q)\right)
  \right),
\end{equation*}
where
\begin{equation*}
  \mu_{\Delta t}(q)=q+\left[-D(q)\nabla V(q)+\beta^{-1}\div D(q)\right]\Delta t.
\end{equation*}
MALA is summarized in Algorithm~\ref{alg:MALA}. It generates a Markov chain~$(q^{n})_{n\geqslant0}$ which is reversible with respect to~$\pi$.

\begin{algorithm}
  \caption{MALA.}
  \label{alg:MALA}
  Consider an initial condition~$q^{0}\in\calQ$, and set~$n=0$.
  \begin{enumerate}[label={[\thealgorithm.\roman*]}, align=left]
      \item \label{step:MALA_1} Compute the proposal~$\widetilde{q}^{n+1}$ as in~\eqref{eq:MALA_proposal};
      \item \label{step:MALA_2} Draw a random variable~$U^{n}$  with uniform law on~$[0,1]$:
      \begin{itemize}[label=$\bullet$]
          \item if~$U^{n}\leqslant r(q^{n},\widetilde{q}^{n+1})$ where
          \begin{equation*}
            r(q,q')=\min\left(1,\frac{\pi(q')T(q',q)}{\pi(q)T(q,q')}\right)
          \end{equation*}
          accept the proposal and set~$q^{n+1}=\widetilde{q}^{n+1}$;
          \item else reject the proposal and set~$q^{n+1}=q^{n}$;
      \end{itemize}
      \item Increment~$n$ and go back to~\ref{step:MALA_1}.
  \end{enumerate}
\end{algorithm}

\paragraph{Square root, divergence, inverse and determinant of the diffusion.} To implement MALA with the diffusions~\eqref{eq:optimal_diffusion_class}, one needs to compute the square root and the divergence of the diffusion in Step~\ref{step:MALA_1}, and the determinant and the inverse of the diffusion in Step~\ref{step:MALA_2}. The particular form of the diffusion operator~\eqref{eq:optimal_diffusion_class} is such that these quantities are available analytically without extra computational cost. Indeed, it is easily verified that
\begin{equation}
  \label{eq:sqrt_inv_det_diffusion}
  \left\lbrace
  \begin{aligned}
    D_\alpha(q)^{1/2}&=\sqrt{\kappa_\alpha}\left[\rmI_d+\left(\sqrt{a_\alpha(\xi(q))}-1\right)P(q)\right],\\
    D_\alpha(q)^{-1}&=\kappa_\alpha^{-1}\left[\rmI_d+\left(\frac{1}{a_\alpha(\xi(q))}-1\right)P(q)\right],\\
    \det D_\alpha(q)&=\kappa_\alpha^{d}a_\alpha(\xi(q)).
  \end{aligned}
  \right.
\end{equation}
As for the divergence, one readily checks that
\begin{equation}
  \label{eq:div_P}
  \div P(q)=-\div P^{\perp}(q)=\frac{\nabla^{2}\xi(q)\nabla\xi(q)}{\left\lVert\nabla\xi(q)\right\rVert^{2}}+\frac{\Delta\xi(q)}{\left\lVert\nabla\xi(q)\right\rVert^{2}}\nabla\xi(q)-2\frac{\nabla\xi(q)^{\sfT}\nabla^{2}\xi(q)\nabla\xi(q)}{\left\lVert\nabla\xi(q)\right\rVert^{4}}\nabla\xi(q),
\end{equation}
where~$\Delta$ is the Laplacian operator. Besides, it holds
\begin{equation}
  \label{eq:chain_rule_div}
  \div(a_{\alpha}(\xi(q))P(q))=a_{\alpha}(\xi(q))\div P(q)+a_{\alpha}'(\xi(q))P(q)\nabla\xi(q)=a_{\alpha}(\xi(q))\div P(q)+a_{\alpha}'(\xi(q))\nabla\xi(q),
\end{equation}
where we used~\eqref{eq:P_nabla_xi_relations} for the second equality. Therefore, the divergence of the diffusion~\eqref{eq:optimal_diffusion_class} is given by
\begin{align}
  \label{eq:divergence_diffusion}
  \div D_{\alpha}(q)&=\kappa_{\alpha}
    \left(a_{\alpha}(\xi(q))-1\right)\left(\frac{\nabla^{2}\xi(q)\nabla\xi(q)}{\left\lVert\nabla\xi(q)\right\rVert^{2}}+\frac{\Delta\xi(q)}{\left\lVert\nabla\xi(q)\right\rVert^{2}}\nabla\xi(q)-2\frac{\nabla\xi(q)^{\sfT}\nabla^{2}\xi(q)\nabla\xi(q)}{\left\lVert\nabla\xi(q)\right\rVert^{4}}\nabla\xi(q)\right)\\
    &\nonumber\quad+\kappa_{\alpha}a_{\alpha}'(\xi(q))\nabla\xi(q).
\end{align}
Moreover, the derivative of the map~$a_{\alpha}$ introduced in~\eqref{eq:optimal_diffusion_class} is given by
\begin{equation}
  \label{eq:derivative_a}
  a_{\alpha}'(z)
  =\frac{\rme^{\alpha\beta F(z)}}{\sigma^{4}(z)}\left[
    \alpha\beta F'(z)\sigma^{2}(z)-\left(\sigma^{2}\right)'(z)
  \right]
  =\frac{\beta\rme^{\alpha\beta F(z)}}{\sigma^{4}(z)}\left[
    \left(\alpha-1\right) \sigma^{2}(z)F'(z)-b(z)
  \right],
\end{equation}
where we used~\eqref{eq:standard_drift_noise_relation_effective_dynamics} for the second equality. The last expression is interesting since it involves only averages with respect to conditional measures.

\begin{remark}
  \label{rem:MALA_gradient_xi_constant}
  When~$\left\lVert\nabla\xi\right\rVert$ is constant, which is the case for the numerical example of Section~\ref{sec:MALA_numerics}, the quantity~$\nabla^{2}\xi\nabla\xi$ vanishes, so that~\eqref{eq:divergence_diffusion} simplifies as
  \begin{equation*}
    \div D_{\alpha}(q)=\kappa_{\alpha}\left(
    \left[a_{\alpha}(\xi(q))-1\right]\frac{\Delta\xi(q)}{\left\lVert\nabla\xi(q)\right\rVert^{2}}+a_{\alpha}'(\xi(q))
  \right)\nabla\xi(q).
  \end{equation*}
  Likewise, the derivative of the map~$a_{\alpha}$ reduces to~$a_{\alpha}'(z)=\alpha\beta F'(z)a_{\alpha}(z)$. In that case, only the free energy and the mean force are needed to construct the diffusion~\eqref{eq:optimal_diffusion_class} and compute its derivatives.
\end{remark}

To run MALA, one therefore needs to have access to the free energy~\eqref{eq:free_energy}, the mean force~\eqref{eq:mean_force} as well as the effective drift and diffusion defined in~\eqref{eq:standard_effective_drift_noise}. Two approaches can be undertaken to estimate these quantities:
\begin{enumerate}[align=left, label=(\roman*)]
  \item They can be precomputed, by sampling the conditional expectations~$\rmd\pi^{\xi}(\cdot|z)$ using constrained sampling methods (see Appendix~\ref{app:TI} where we make precise the thermodynamic integration method for our numerical experiment);
  \item They can be learned on the fly, as presented below in Section~\ref{sec:adaptive_learning}.
\end{enumerate}

Finally, the scalar~$\kappa_{\alpha}$ in~\eqref{eq:optimal_diffusion_class} is given by
\begin{equation}
  \label{eq:kappa_alpha}
  \kappa_{\alpha}=\left(\int_{\xi(\calQ)}\sqrt{d-1+a_{\alpha}(z)^{2}}\,\rme^{-\beta F(z)}\rmd z\right)^{-1}.
\end{equation}
The computations leading to the formula~\eqref{eq:kappa_alpha} are detailed in Appendix~\ref{app:kappa_alpha}. Note that the quantity~$\kappa_{\alpha}$ simply rescales the time step~$\Delta t$ when discretizing the overdamped Langevin dynamics. It is therefore not strictly necessary to compute it in order to run the MALA algorithm. In particular, in the numerical experiments of Section~\ref{sec:MALA_numerics}, we optimize over the time step~$\Delta t$ for each value of~$\alpha$, so that computing~$\kappa_{\alpha}$ is only performed in order to renormalize the time steps in the presentation of the numerical results. However, the computation of~$\kappa_{\alpha}$ becomes relevant when learning the free energy~$F$ on the fly (as described in Section~\ref{sec:adaptive_learning}), in order to properly scale the updated values of the diffusion with respect to the time step~$\Delta t$ and avoid that it takes too large values.

\subsection{MALA: numerical results}
\label{sec:MALA_numerics}

We illustrate efficiency gains in sampling when discretizing overdamped Langevin dynamics~\eqref{eq:overdamped_langevin_diffusion} with diffusion~\eqref{eq:optimal_diffusion_class}, utilizing the MALA algorithm described in Algorithm~\ref{alg:MALA}. The system we consider has been used as a toy model for solvation, see for instance~\cite{straub_1988,dellago_1999, lelievre_2007,lelievre_2010, nilmeier_2011,lelievre_2012,das_2022}.

\paragraph{Physical system.}

We consider a two-dimensional system composed of~$N=16$ particles (the dimension is then~$d=32$). The particles are placed in a periodic square box with side length~$\ell$ such that the particle density is~$N/\ell^{2}=0.7$ (thus~$\ell\approx 4.78$). The particles interact through the following repulsive WCA pair potential:
\begin{equation*}
  V_{\rm WCA}(r)=\left\lbrace\begin{aligned}
    4\varepsilon\left[\left(\frac{R}{r}\right)^{12}-\left(\frac{R}{r}\right)^{6}\right]+\varepsilon &  & \text{if }r\leqslant r_0, \\
    0                                                                                                         &  & \text{if }r>r_0.
  \end{aligned}
  \right.
\end{equation*}
Here,~$r$ denotes the distance between two particles,~$\varepsilon=1$ and~$R=1$ are two positive parameters and~$r_0=2^{1/6}R$. The first two particles are designated to form a solute dimer while the others are solvent particles. For these two particles, the WCA potential is replaced by the following double-well potential
\begin{equation}
  \label{eq:dimer_potential}
  V_{\rm DW}(r)=h\left(1-\frac{(r-r_1-w)^{2}}{w^{2}}\right)^{2},
\end{equation}
where $w=0.35$ and~$r_1=\ell/4-w$ is chosen such that $V_{\rm DW}(0)=V_{\rm DW}(\ell/2)$ which yields $r_1\approx 0.85$. This potential admits two minima, one corresponding to the compact state~$r=r_1$, and one corresponding to the stretched state~$r=r_1+2w\approx1.55$. We set $h=2$ ($h$ represents the energy barrier separating the two states in the double-well potential~\eqref{eq:dimer_potential}). Note that~$r_0<\ell/2$, so that a particle cannot interact with any of its periodic copies. We set~$\beta=1$. The values we chose (rather dense system in terms of solvent particles, relatively high energy barrier) ensure that we have a sufficiently difficult sampling problem allowing to demonstrate the benefit of using a modified diffusion.

To make the potential energy precise, we denote by~$q=(q_1,\dots,q_N)\in\left(\ell\bbT\right)^{2N}$ the positions of the system, with~$q_1$ and~$q_2$ forming the dimer. The potential energy therefore writes
\begin{equation}
    \label{eq:numerical_example_potential_energy_function}
    V(q)=V_{\rm DW}\left(\left\lVert q_2-q_1\right\rVert\right)+\sum_{\substack{i\in\left\lbrace1,2\right\rbrace\\3\leqslant j\leqslant N}}V_{\rm WCA}\left(\left\lVert q_i-q_j\right\rVert\right)+\sum_{3\leqslant i<j\leqslant N}V_{\rm WCA}\left(\left\lVert q_i-q_j\right\rVert\right).
\end{equation}
Note that the distances appearing in~\eqref{eq:numerical_example_potential_energy_function} are computed taking the periodic boundary conditions into account.

\paragraph{Collective variable.}
The collective variable~$\xi$ is defined by
\begin{equation*}
  \xi(q)=\frac{\left\lVert q_2-q_1\right\rVert-r_0}{2w}.
\end{equation*}
The value of the collective variable is 0 when the dimer is in a compact state, and 1 when it is in a stretched state. We therefore define two subsets of the configuration space that correspond to these two metastable states:
\begin{equation}
  \label{eq:metastable_states}
  C_{0}=\xi^{-1}(-\infty,\eta),\qquad C_{1}=\xi^{-1}(1-\eta,+\infty),
\end{equation}
where~$\eta=0.1$ is a numerical parameter that sets a tolerance to identify the two states. The latent space~$\xi(\calQ)$ is approximated by the interval~$[z_{\rm min},z_{\rm max}]=[-0.2,1.225]$ (which is the range of values of the collective variable observed during test runs). Note that the collective variable is not smooth everywhere, but only on~$\left\lbrace q\in (\ell\bbT)^{2N}\,\middle|\,q_1\neq q_2\right\rbrace$. When performing simulations, the repulsive interactions between the atoms composing the dimer prevent the dynamics from coming close to the region of singularity.

\begin{remark}
  Since~$\xi$ is related to a distance up to periodic boundary conditions, it is important to check that the dimer can actually expand enough in the simulation box so that~$\xi$ can be larger than 1. In fact, one requires that~$r_1+2wz_{\rm max}<\ell/2$. In practice, we define a tolerance~$z_{\rm tol}$ such that~$r_1+2wz_{\rm tol}<\ell/2$ with~$z_{\rm max}\leqslant z_{\rm tol}$, and check whether the collective variable exceeds this tolerance during the simulation. For our numerical experiments, we choose~$z_{\rm tol}=2.0$.
\end{remark}

Note that the collective variable is a function of only the first two particles of this two-dimensional system, meaning that its derivatives vanish for components other than the first four. In fact, it holds
\begin{equation}
  \label{eq:xi_first_derivatives}
  \nabla\xi(q) = \frac{1}{2w\left\lVert q_2-q_1\right\rVert} \begin{pmatrix}
    x_1-x_2 \\ y_1-y_2 \\ -(x_1-x_2) \\ -(y_1-y_2) \\ 0\\\vdots\\0
  \end{pmatrix}\in \bbR^{2N}.
\end{equation}

\paragraph{Diffusion.}
Since the collective variable~$\xi(q)$ is related to a bond length, the map~$\sigma$ defined in~\eqref{eq:standard_effective_drift_noise} is constant, and is simply equal to (see~\eqref{eq:xi_first_derivatives})
\begin{equation*}
  \left\lVert\nabla\xi\right\rVert^{2}=\frac{1}{2w^{2}}\approx 4.08.
\end{equation*}
We can also use the simplifications described in Remark~\ref{rem:MALA_gradient_xi_constant}. The numerical parameter~$\alpha\geqslant0$ appearing in~\eqref{eq:optimal_diffusion_class} should typically be chosen such that~$\alpha\beta F$ is of order 1 over~$[z_{\rm min},z_{\rm max}]$ to reduce numerical instabilities. In our numerical setting, the variations of~$F$ are of order~$h$ over this interval, so that~$\alpha$ should to be of order~$1/(\beta h)$. We use a left-Riemann rule to approximate the value of~$\kappa_{\alpha}$ defined by~\eqref{eq:kappa_alpha}.

\paragraph{Efficient implementation of MALA.} Since~$\xi$ is a function of only~4 components, the additional computations and storage (compared to the standard MALA algorithm with identity diffusion) is limited as the diffusion~\eqref{eq:optimal_diffusion_class} is block diagonal, with a~$4\times4$ block followed by a scalar~$\left(2N-4\right)\times\left(2N-4\right)$ matrix which does not depend on the position. More precisely,
\begin{equation}
  \label{eq:diffusion_block_diagonal}
  D_{\alpha}(q)=\kappa_{\alpha}\begin{pmatrix}
    \rmI_{4}+(a_{\alpha}(\xi(q))-1)\widetilde{D}(q) & 0_{4,2N-4}  \\
    0_{2N-4,4}                & \rmI_{2N-4}
  \end{pmatrix}\in\bbR^{2N\times 2N},
\end{equation}
with
\begin{equation*}
  \widetilde{D}(q)=\frac{1}{2}\begin{pmatrix}
    A(q)  & -A(q) \\
    -A(q) & A(q)
  \end{pmatrix}\in\bbR^{4\times 4},\quad
  A(q) = \frac{q_1 - q_2}{\left\lVert q_1-q_2\right\rVert}\otimes\frac{q_1 - q_2}{\left\lVert q_1-q_2\right\rVert}\in\bbR^{2\times 2}.
\end{equation*}
Additionally, this means that
\begin{itemize}
  \item the gradient of~$\xi$ has at most 4 nonzero components (see~\eqref{eq:xi_first_derivatives}), and the only nonzero part of the Hessian of~$\xi$ is a~$4\times 4$ submatrix. In fact, for our numerical example, the gradient of~$\xi$ only has two degrees of freedom while its Hessian only has three, thanks to the relations between the various derivatives of~$\xi$ in~\eqref{eq:xi_first_derivatives} (and in~\eqref{eq:xi_second_derivatives}, see Appendix~\ref{app:computations});
  \item the divergence of the diffusion~\eqref{eq:divergence_diffusion} only acts through derivatives of~$\xi$, so that only its first four components can be nonzero. In fact, for our numerical experiment, the divergence only has two degrees of freedom (see~\eqref{eq:divergence_diffusion_numerical_experiment} in Appendix~\ref{app:computations}).
\end{itemize}
Therefore, only the update related to the first two particles is notably modified compared to the standard case~$D=\rmI_d$. These observations can be generalized to any collective variable which is a function of only~$k\leqslant d$ components of the positions.

\begin{figure}
  \centering
  \begin{subfigure}{0.45\textwidth}
    \centering
    \includegraphics[width=\textwidth]{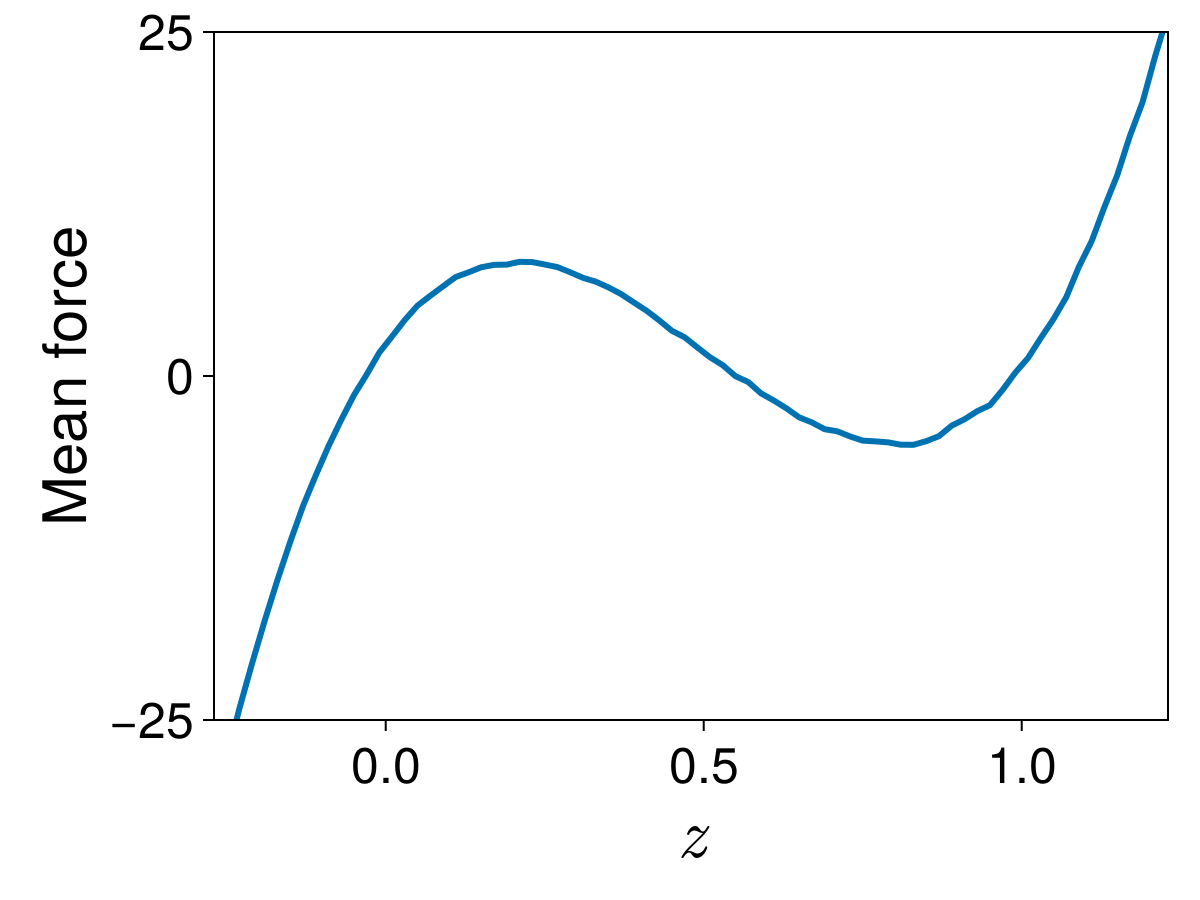}
    \caption{Mean force~$F'$ vs. collective variable~$\xi$.}
    \label{fig:mean_force_abf}
  \end{subfigure}\hfill
  \begin{subfigure}{0.45\textwidth}
    \centering
    \includegraphics[width=\textwidth]{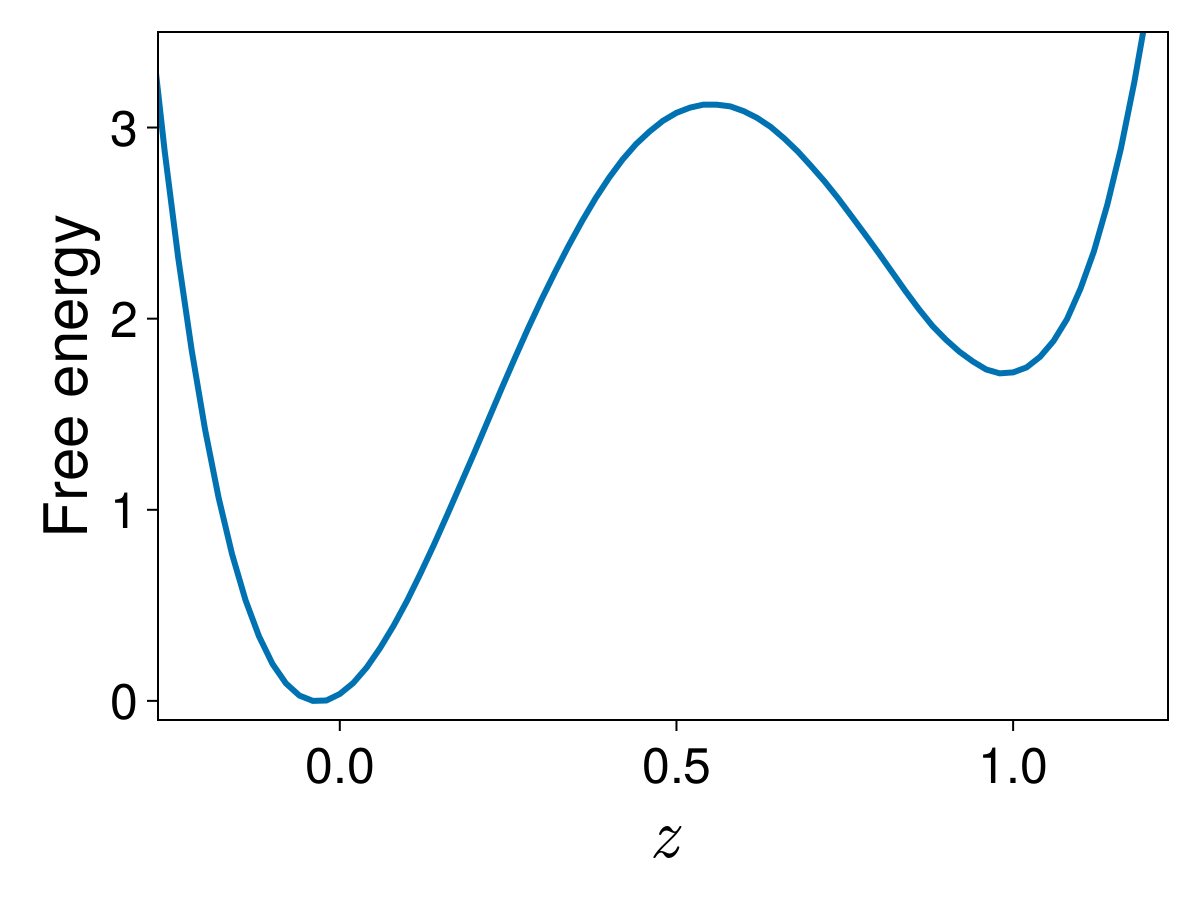}
    \caption{Free energy~$F$ vs. collective variable~$\xi$.}
    \label{fig:free_energy_abf}
  \end{subfigure}
  \caption{Mean force and free energy computed with thermodynamic integration.}
\end{figure}

\paragraph{Setting of the numerical experiment.} We first assume that the free energy~$F$ and mean force~$F'$ have been precomputed (we will explain in Section~\ref{sec:adaptive_learning} how to learn these quantities on the fly). To estimate~$F'$ (and thus~$F$), we used thermodynamic integration as described in Appendix~\ref{app:TI}. The mean force and the free energy obtained with this procedure are plotted in Figures~\ref{fig:mean_force_abf} and~\ref{fig:free_energy_abf}, respectively. Observe that the free energy well associated with the compact state is lower than the one corresponding to the stretched state, which means that the compact state is more likely than the stretched one. The free energy is defined up to an additive constant: we fix this constant so that the minimum of the free energy over the interval~$[z_{\rm min},z_{\rm max}]$ is~0. We therefore deviate from the definition given in~\eqref{eq:free_energy}, and the quantity~$\rme^{-\beta F(z)}$ does not normalize to~1 on~$\xi(\calQ)$. Note that this does not impact the form of the optimal homogenized effective diffusion, but only alters the definition of~$a_{\alpha}$ by introducing a multiplicative constant, see Remark~\ref{rem:free_energy_constant}.


\paragraph{Efficiency metric.}

To observe the efficiency of using nonconstant diffusion matrices, we compute the mean number of iterations to go back and forth from one metastable state to the other. More precisely, we consider the initial condition~$q^{0}=(q^{0}_1,\dots,q^{0}_N)\in(\ell\bbT)^{2N}$ defined by
\begin{equation}
  \label{eq:q0}
  \forall\, 1\leqslant i \leqslant N,\qquad q^{0}_i=(x^{0}_i,y^{0}_i),\qquad
  \left\lbrace
  \begin{aligned}
    x^{0}_{i}&=a(0.5+\left\lfloor (i-1)/4\right\rfloor)\\
    y^{0}_{i}&=a(0.5+((i-1)\bmod{4})),
  \end{aligned}
  \right.
\end{equation}
where~$a=\ell/\sqrt{N}\approx1.20$, except for~$y^{0}_{2}$ which is set to~$y^{0}_{1}+r_1$ so that~$\xi(q^{0})=0$. This particular configuration is such that the dimer is in a compact state, and the other particles are initially set on a lattice. We introduce the additional variable~$\Theta^{0}=0$, which keeps track of whether the dimer last visited the compact (0) or the stretched (1) state. We define the following stopping times and values of the additional variable:~$\tau^{0}=0$, and
\begin{equation*}
  \tau^{k+1}=\min\limits_{n\geqslant 1}\left\lbrace q^{n+\sum_{i=1}^{k}\tau^{i}}\in C_{\Theta^{k+1}}\right\rbrace,\qquad \Theta^{k+1}=1-\Theta^{k},
\end{equation*}
where the sets~$C_{0},C_{1}$ are defined in~\eqref{eq:metastable_states}. Each value~$\tau^{i}$ corresponds to the number of iterations of the dynamics used to transition between~$C_0$ to~$C_1$ or~$C_1$ to~$C_0$, depending on the parity of~$i$. The simulation is run until~$K=10^{5}$ transitions are observed. The mean number of iterations to observe a transition, for a given value of~$\alpha$ and time step~$\Delta t$, is then estimated as
\begin{equation*}
  \widehat{\tau}(\alpha,\Delta t)=\frac{1}{K}\sum_{k=1}^{K}\tau^{k}.
\end{equation*}
We denote by~$\widehat{\tau}(\alpha)^{\star}$ the minimum value of~$\widehat{\tau}(\alpha,\cdot)$ over the range of time steps~$\Delta t$ used in this numerical experiment. 

We perform the same numerical experiment with the constant diffusion~$D_{\rm cst}(q)=\kappa\rmI_d$, where~$\kappa>0$ acts as the same normalizing constant as~$\kappa_\alpha$ in~\eqref{eq:kappa_alpha}: it is a left-Riemman rule approximation of
\begin{equation*}
  \left(\int_{\xi(\calQ)}\sqrt{d}\,\rme^{-\beta F(z)}\rmd z\right)^{-1}.
\end{equation*}
Since we used the same normalization, we can therefore compare the numerical results when using the various diffusions. When using the constant diffusion,~the estimate~$\widehat{\tau}(\alpha,\Delta t)$ (respectively~$\widehat{\tau}(\alpha)^{\star}$) is denoted by~$\widehat{\tau}_{\rm cst}(\Delta t)$ (respectively~$\widehat{\tau}_{\rm cst}^{\star}$).

To perform the numerical experiment, we choose 16 values of the time step evenly spaced log-wise in the interval~$[5\times10^{-4},10^{-2}]$. We present the results for values of~$\alpha$ such that
\begin{equation*}
  \alpha\in\left\lbrace 0, 0.1, 0.2,\dots,2.4\right\rbrace.
\end{equation*}
For larger values of~$\alpha$, the mean number of iterations dramatically increases, and eventually the numerical scheme becomes unstable (\emph{i.e.~}we observe numerical overflows for the time steps we consider). 

\begin{figure}
  \centering
  \begin{subfigure}[t]{0.475\columnwidth}
    \includegraphics[width=\columnwidth]{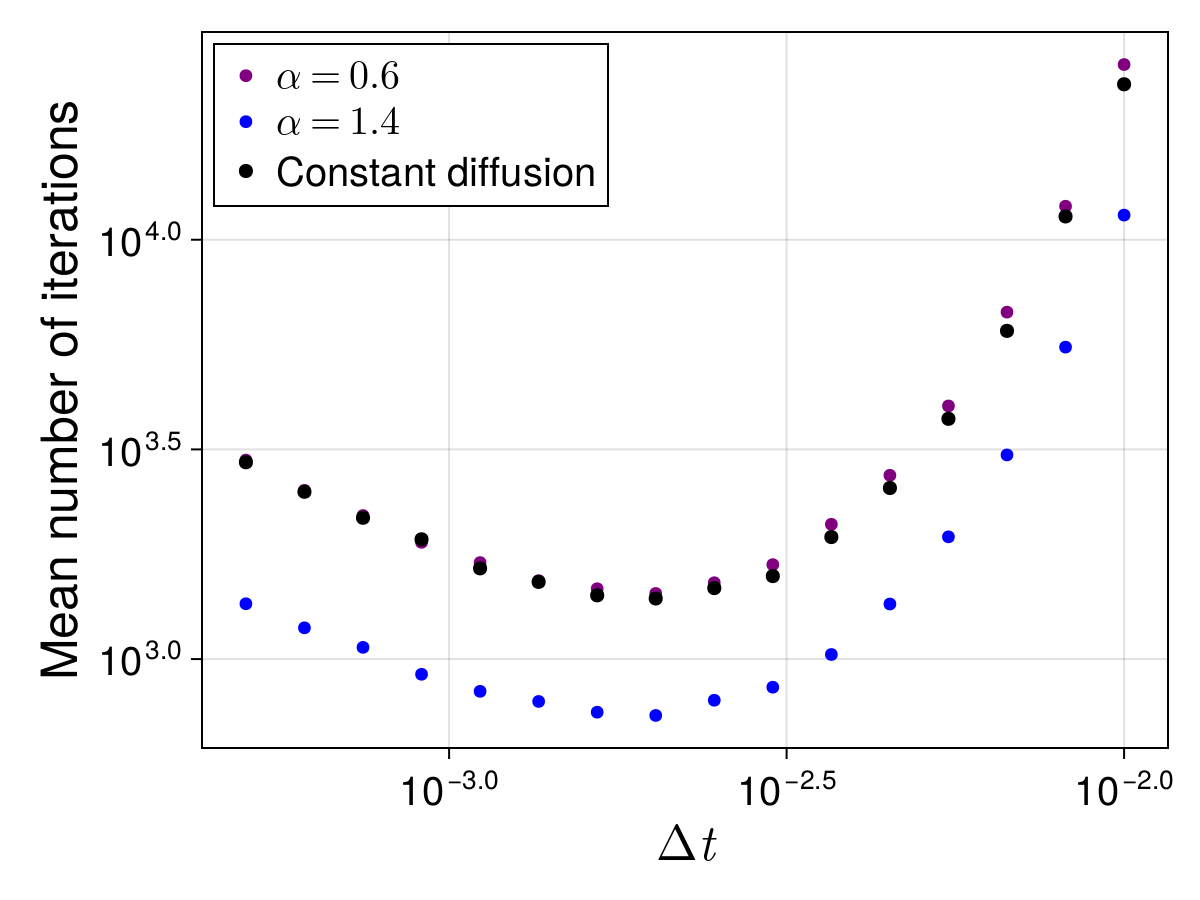}
    \caption{Mean number of iterations to observe a transition as a function of~$\Delta t$ for two values of~$\alpha$, and for the constant diffusion.}
    \label{fig:some_values_n_iterations}
  \end{subfigure}\hfill
  \begin{subfigure}[t]{0.475\columnwidth}
    \centering
    \includegraphics[width=\columnwidth]{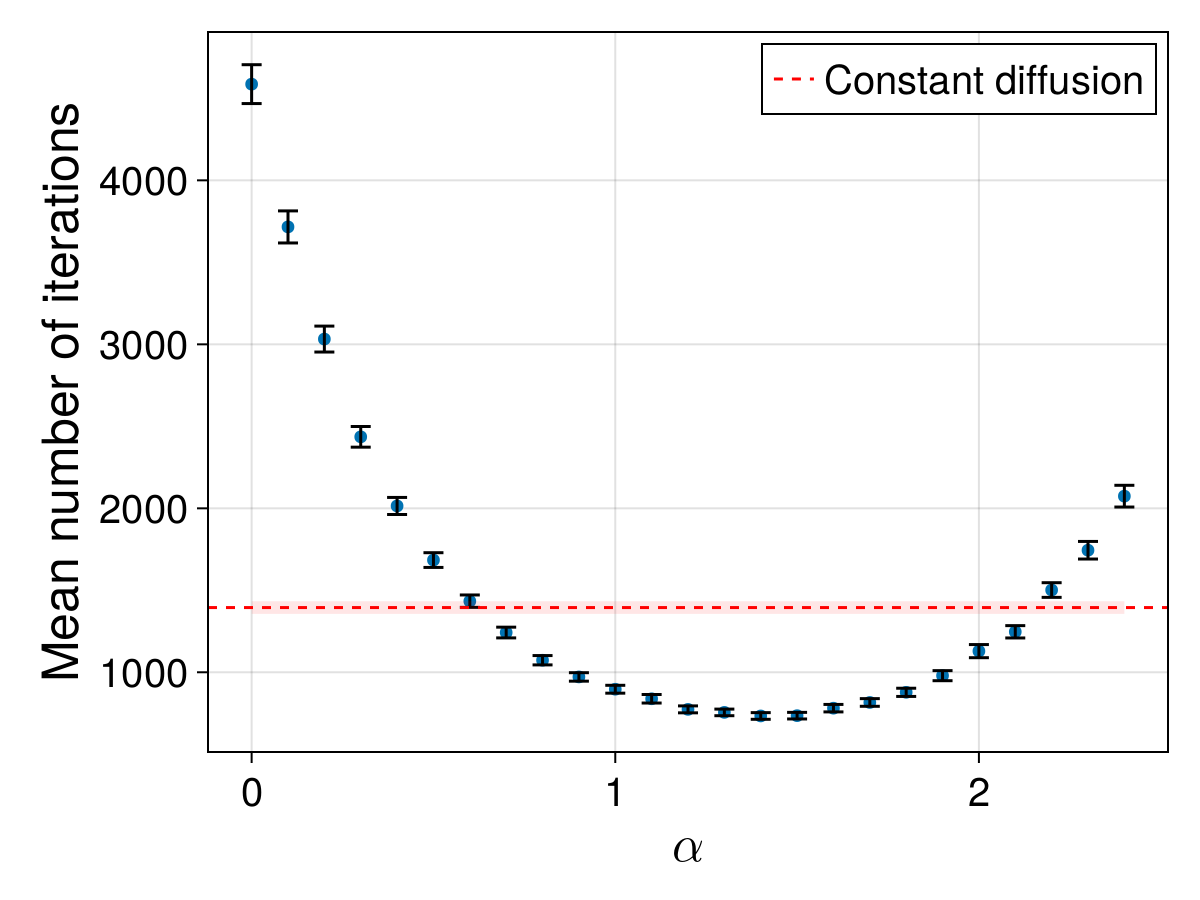}
    \caption{Minimum values of the transition times over the time steps as a function of~$\alpha$, and minimum transition time for the constant diffusion.}
    \label{fig:minimum_values}
  \end{subfigure}
  \caption{Sampling efficiency using the MALA algorithm~\ref{alg:MALA}.}
  \label{fig:MALA}
\end{figure}

\paragraph{Results.} We present in Figure~\ref{fig:MALA} the mean number of iterations to observe a transition as a function of the time step~$\Delta t$ for various values of~$\alpha$. The error bars represent 95\% confidence intervals.  The minimal number of iterations across all simulations is obtained for~$\alpha_{\rm opt}=1.4$. For this value, it holds~$\widehat{\tau}(\alpha_{\rm opt})^{\star}=733.71$, while~$\widehat{\tau}_{\rm cst}^{\star}=1394.74$. This shows that introducing multiplicative noise with the diffusion~\eqref{eq:optimal_diffusion_class} is useful to cross (free) energy barriers, thus enhancing the exploration efficiency of the configuration space. Figure~\ref{fig:minimum_values} shows the behavior of~$\alpha\mapsto\widehat{\tau}(\alpha)^{\star}$, and compares it to the value~$\widehat{\tau}_{\rm cst}^{\star}$. The minimum values were all obtained for time steps which were similar, around~$\Delta t\approx 2\times10^{-3}$ (thanks to the normalization of the diffusion). We observe that the numerical results obtained using the nonconstant diffusions are better for a large range of values of~$\alpha$. Figure~\ref{fig:some_values_n_iterations} shows that, when~$\widehat{\tau}(\alpha_1)^{\star}\leqslant\widehat{\tau}(\alpha_2)^{\star} $ (respectively~$\widehat{\tau}(\alpha_1)^{\star}\leqslant \widehat{\tau}_{\rm cst}^{\star}$) for two values of~$\alpha\in\left\lbrace\alpha_1,\alpha_2\right\rbrace$, then~$\widehat{\tau}(\alpha_1,\Delta t)\leqslant\widehat{\tau}(\alpha_2,\Delta t)$ (respectively~$\widehat{\tau}(\alpha_1,\Delta t)\leqslant\widehat{\tau}_{\rm cst}(\Delta t)$) for any time step~$\Delta t$ used in this numerical experiment. This means that the choice of the time step is in fact irrelevant, and only the value of~$\alpha$ impacts the numerical results. Finally, note that for small time steps, the mean number of iterations scales as~$\Delta t^{-1}$, since the numerical scheme converges to the continuous limit. 

\subsection{Adaptive learning of the diffusion}
\label{sec:adaptive_learning}

When the free energy~\eqref{eq:free_energy}, the mean force~\eqref{eq:mean_force}, or the effective drift and diffusion~\eqref{eq:standard_effective_drift_noise} are not precomputed, they can be learned during the simulation, in the spirit of free energy adaptive biasing techniques  such as the Adapted Biasing Force (ABF) methods~\cite{darve_2001,henin_2004}. These methods are used to efficiently compute free energy differences between metastable states, \emph{i.e.}~relative likelihood of states. More precisely, standard dynamics such as~\eqref{eq:standard_overdamped_langevin} are modified by adding a force which biases the dynamics towards unexplored regions. This additional force is an estimate of the mean force~$F'$, which is learned on the fly, and the bias is updated in order for the dynamics to escape metastable states.

In practice, ABF methods uses empirical averages which are updated on the fly to estimate conditional expectations with respect to~$\pi^{\xi}$ (see~\eqref{eq:conditional_measure}). The diffusion~\eqref{eq:optimal_diffusion_class} depends on these conditional expectations, through the effective diffusion function~$\sigma$ and free energy~$F$ (recovered by integrating the mean force~$F'$). The divergence of the diffusion, whose analytical expression is given in~\eqref{eq:divergence_diffusion}, requires to additionally compute the conditional expectations of the effective drift~$b$, see~\eqref{eq:derivative_a}. We can therefore learn these quantities when performing the simulation, and update the diffusion and its divergence on the fly. Let us make precise one way to approximate expectations with respect to the conditional measure~$\pi^{\xi}$. We refer to~\cite[Section~5.1.3]{lelievre_2010} for general discretization methods to approximate these quantities.

Let~$q_t$ solves the overdamped Langevin dynamics~\eqref{eq:overdamped_langevin_diffusion}. To compute a conditional expectation of the form
\begin{equation}
  \label{eq:conditional_expectation}
  \int_{\Sigma(z)}h\,\rmd\pi^{\xi}(\cdot|z),
\end{equation}
where~$h:\calQ\to\bbR$, the idea is to update values of a map~$H_t(z)$ along the simulation (where~$t$ represents the simulation time) in order for~$H_t(z)$ to converge to~\eqref{eq:conditional_expectation} in the long time limit. To naturally construct these maps, we rely on the formal identity
\begin{equation*}
  \int_{\Sigma(z)}h\,\rmd\pi^{\xi}(\cdot|z)=\lim\limits_{\Delta z\to 0}\frac{\bbE\left[h(q_t)\mathbbm{1}_{\left\lbrace \xi(q_t)\in(z,z+\Delta z)\right\rbrace}\right]}{\bbP\left(\xi(q_t)\in (z,z+\Delta z)\right)}.
\end{equation*}
Both the numerator and denominator can be approximated using empirical averages over a trajectory and/or multiple replicas, introducing a binning of the collective variable values to get a piecewise constant (in~$z$) approximation. This is made precise in the next paragraph.

\paragraph{The discretization in practice.} We now specify the methodology used in the numerical experiments. Let us consider an interval~$[z_{\rm min},z_{\rm max}]$ of the latent space~$\xi(\calQ)$ over which the map~$H_t$ is approximated. The interval~$[z_{\rm min},z_{\rm max}]$ should be chosen so that most of the values of the collective variable lie in it. Outside this interval, the value of~$H_t$ is simply set to a constant so that $H_t$ is continuous. We next discretize this interval using~$N_z$ bins by defining
\begin{equation*}
  z^{i}=z_{\rm min}+i\Delta z,\qquad \Delta z=\frac{z_{\rm max}-z_{\rm min}}{N_z}.
\end{equation*}
We thus construct an approximation of the map~$H_t$ using a piecewise constant map on~$[z_{\rm min},z_{\rm max}]$. Denote by~$(q^{k})_{0\leqslant n\leqslant N}$ a trajectory obtained from a discretization of~\eqref{eq:overdamped_langevin_diffusion} (\emph{e.g.}~using MALA). The value of~$H_{t}(z)$ at time~$t_n=n\Delta t$ is then approximated as
\begin{equation}
  \label{eq:updating_formulas}
  H_{n\Delta t}^{\Delta z}\left(\left\lfloor\frac{z-z_{\rm min}}{\Delta z}\right\rfloor\right),
  \qquad H_{n\Delta t}^{\Delta z}(i)=\frac{\displaystyle\sum_{j=0}^{n}h(q^{j})\mathbbm{1}_{\left\lbrace \xi(q^{j})\in(z^i,z^{i+1})\right\rbrace}}{\displaystyle\sum_{j=0}^{n}\mathbbm{1}_{\left\lbrace \xi(q^{j})\in(z^i,z^{i+1})\right\rbrace}}.
\end{equation}

\paragraph{Hyperparameters.} ABF methods usually use additional hyperparameters, which prevent numerical instabilities when discretizing the dynamics. We consider two of them here:
\begin{itemize}
  \item The learned diffusion is used in a given bin only if this bin has been visited a sufficiently large number of times, the threshold being given by~$N_{\rm min}$.
  \item The update of the conditional expectation is performed every~$N_{\rm update}$ steps of the dynamics.
\end{itemize}

\paragraph{Conditional expectations computed for the diffusion~\eqref{eq:optimal_diffusion_class}.} 
Two conditional expectations are needed for the diffusion~\eqref{eq:optimal_diffusion_class}: one related to~$\sigma^{2}$ defined by~\eqref{eq:standard_effective_drift_noise}, and one related to the mean force~$F'$ defined by~\eqref{eq:mean_force}. Note that when~$\left\lVert\nabla\xi\right\rVert$ is constant, one only needs to approximate the mean force. Moreover, the scalar~$\kappa_{\alpha}$ appearing in~\eqref{eq:optimal_diffusion_class} also depends on the free energy and the effective diffusion (see~\eqref{eq:kappa_alpha}): it therefore should be updated along the simulation as well.

When approximating the mean force~$F'$ (respectively the effective diffusion function~$\sigma$), we choose 0 (respectively~1) as the value used for the mean force (respectively the effective diffusion) when there are not enough observations in a bin or when the value of the collective variable is outside the interval~$[z_{\rm min},z_{\rm max}]$.

To approximate the values of the free energy~$F$, we simply integrate the approximation of the mean force over the interval~$[z_{\rm min},z_{\rm max}]$. The values of the free energy outside the mesh are set in order to obtain a continuous function. Lastly, we fix the additive constant of the free energy so that the minimum value of the free energy over the interval~$[z_{\rm min},z_{\rm max}]$ is~0.

\paragraph{Conditional expectations computed for the divergence of the diffusion~\eqref{eq:optimal_diffusion_class}.} In view of~\eqref{eq:divergence_diffusion}-\eqref{eq:derivative_a}, the effective drift~$b$ defined in~\eqref{eq:standard_effective_drift_noise} also needs to be approximated, using similar strategies as for the mean force~$F'$ (actually, the mean force and effective drift~$b$ are opposite one another when~$\left\lVert\nabla\xi\right\rVert$ is constant and equal to~1, see~\eqref{eq:standard_drift_noise_relation_effective_dynamics}). Note that the identity~\eqref{eq:drift_noise_relation_effective_dynamics} therefore does not hold exactly in this discretized setting.

\paragraph{Adaptive MALA.} We are now in position to describe the adaptive MALA algorithm, see Algorithm~\ref{alg:adaptive_MALA}. The symbols~$\left(F_{n\Delta t}^{\Delta z}\right)',F_{n\Delta t}^{\Delta z},\sigma_{n\Delta t}^{\Delta z},b_{n\Delta t}^{\Delta z}$, which are maps/vectors defined on~$\left\lbrace 0,N_z-1\right\rbrace$ with values in~$\bbR$, should be understood as the approximations at time~$n\Delta t$ (\emph{i.e.}~iteration~$n$) of the mean force~\eqref{eq:mean_force}, free energy~\eqref{eq:free_energy} and effective diffusion and drift functions~\eqref{eq:standard_effective_drift_noise} respectively. The symbols~$\kappa_{\alpha,n\Delta t}^{\Delta z},D_{\alpha,n\Delta t}^{\Delta z}$ should be understood as the current approximation of the normalization constant of the diffusion (see~\eqref{eq:kappa_alpha}) and the diffusion (see~\eqref{eq:optimal_diffusion_class}) respectively. We use a left-Riemann rule in Step~\ref{step:adaptive_MALA_1} to approximate the integral defining the scalar~$\kappa_{\alpha}$, other choices are of course possible.

Note that Algorithm~\ref{alg:adaptive_MALA} belongs to the class of adaptive MCMC algorithms (see for instance~\cite{haario_1999,haario_2001,roberts_2009}). It is therefore not clear that the Boltzmann-Gibbs measure~$\pi$ is correctly sampled, even though each iteration yields a kernel which is reversible with respect to~$\pi$. A sufficient condition to ensure unbiased sampling in the long-time limit is that all the quantities learned over the simulation eventually converge, see~\cite[Theorem~1]{roberts_2007}. Proving that these methods actually converge is usually not straightforward, as was noted for the long-time convergence of ABF methods, see for instance~\cite{lelievre_2008, lelievre_2011,comer_2015}.

Of course, one can stop the learning process at any iteration, and plug the learned quantities in Algorithm~\ref{alg:MALA}, which will ensure an unbiased sampling. The convergence of these quantities is usually observed to be rather fast in most  applications (compared to the total number of iterations the practitioner uses), which limits the bias. This will be illustrated in Section~\ref{sec:adaptive_MALA_numerics} using our numerical experiments.

\begin{algorithm}
  \caption{Adaptive MALA.}
  \label{alg:adaptive_MALA}
  Consider an initial condition~$q^{0}\in\calQ$, and set~$n=0$. Additionally, set
  \begin{equation*}
    \left\lbrace
    \begin{aligned}
      \left(\left(F_0^{\Delta z}\right)'(i)\right)_{0\leqslant i\leqslant N_z-1}&=0,\\
      \left(F_0^{\Delta z}(i)\right)_{0\leqslant i\leqslant N_z-1}&=0,\\
      \left(b_0^{\Delta z}(i)\right)_{0\leqslant i\leqslant N_z-1}&=0,\\
      \left(\sigma_0^{\Delta z}(i)\right)_{0\leqslant i\leqslant N_z-1}&=\left(1\right)_{0\leqslant i\leqslant N_z-1}.
    \end{aligned}
    \right.
  \end{equation*}
  \begin{enumerate}[label={[\thealgorithm.\roman*]}, align=left]
      \item \label{step:adaptive_MALA_1} Perform one step of the MALA algorithm~\ref{alg:MALA} using the diffusion
      \begin{equation*}
        D_{\alpha,n\Delta t}^{\Delta z}(q)=\kappa_{\alpha,n\Delta t}^{\Delta z}\left(
          P^{\perp}(q)+a_{\alpha,n\Delta t}^{\Delta z}\left(\left\lfloor\frac{\xi(q)-z_{\rm min}}{\Delta z}\right\rfloor\right)P(q)
        \right),\qquad a_{\alpha,n\Delta t}^{\Delta z}(i)=\frac{\rme^{\alpha\beta F_{n\Delta t}^{\Delta z}(i)}}{\left(\sigma_{n\Delta t}^{\Delta z}\right)^{2}(i)}.
      \end{equation*}
      $$
      \kappa_{\alpha,n\Delta t}^{\Delta z}=\left(\Delta z \sum_{i=0}^{N_z-1}\sqrt{d-1+\left(a^{\Delta z}_{\alpha,n\Delta t} (i)\right)^{2}}\,\rme^{-\beta F_{n\Delta t}^{\Delta z}(i)}\right)^{-1}$$
      The divergence of the diffusion is approximated using~\eqref{eq:divergence_diffusion} where~$\kappa_{\alpha},F,F',b$ and~$\sigma$ are replaced by their current approximations.
      \item Increment~$n$. Store the values of the local mean force~$f$ (see~\eqref{eq:mean_force}), of the integrand of~$b$ (see~\eqref{eq:standard_effective_drift_noise}), and of~$\left\lVert\nabla\xi\right\rVert^{2}$. Update the histogram of the collective variable values. If~$n = k N_{\rm update}$ for some integer $k$:
      \begin{enumerate}
        \item Update the effective drift~$b_{n\Delta t}^{\Delta z}$, effective diffusion~$\sigma_{n\Delta t}^{\Delta z}$ and the mean force~$\left(F_{n\Delta t}^{\Delta z}\right)'$ using the formulas~\eqref{eq:updating_formulas}, while considering the additional rules related to~$N_{\rm min}$ and the histogram values in each bin;
        \item Update the free energy~$F_{n\Delta t}^{\Delta z}$;
        \item Update the normalization constant~$\kappa_{\alpha,n\Delta t}^{\Delta z}$.
      \end{enumerate}
      \item Go back to~\ref{step:adaptive_MALA_1}.
  \end{enumerate}
\end{algorithm}

\subsection{Adaptive MALA: numerical results}
\label{sec:adaptive_MALA_numerics}

The numerical experiment is the same as in Section~\ref{sec:MALA_numerics}: the only difference is that the algorithm is Adaptive MALA instead of MALA. Since~$\left\lVert\nabla\xi\right\rVert$ is constant, we only need to learn the mean force~$F'$. We set~$N_{\rm min}=100$ and~$N_{\rm update}=20$. The current approximation of the scalar~$\kappa_{\alpha}$ is simply computed by using a left-Riemann rule every time the free energy and effective diffusion are updated.

\begin{figure}
  \centering
  \begin{subfigure}[t]{0.475\columnwidth}
    \centering
    \includegraphics[width=\columnwidth]{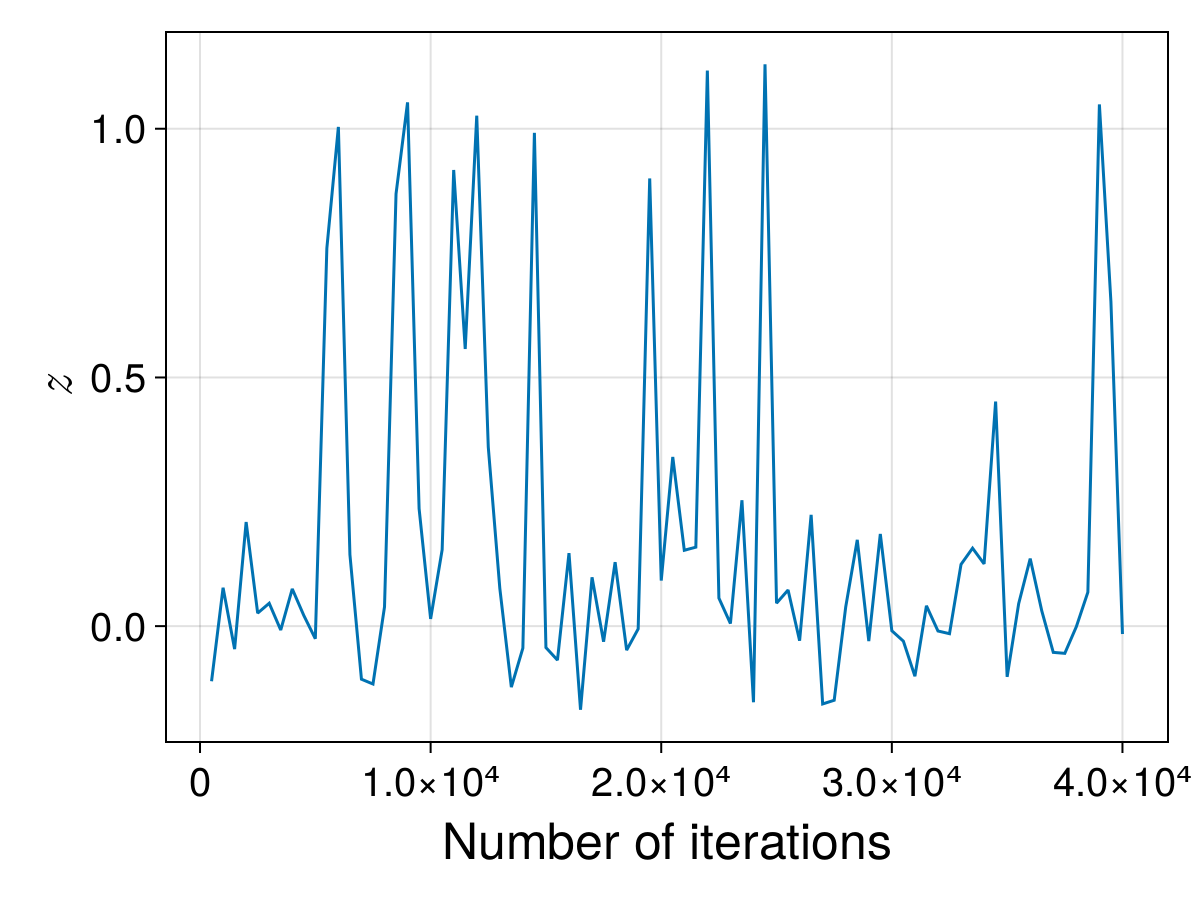}
    \caption{Values of the collective variable as a function of the number of iterations.}
    \label{fig:adaptive_CV}
  \end{subfigure}\hfill
  \begin{subfigure}[t]{0.475\columnwidth}
    \includegraphics[width=\columnwidth]{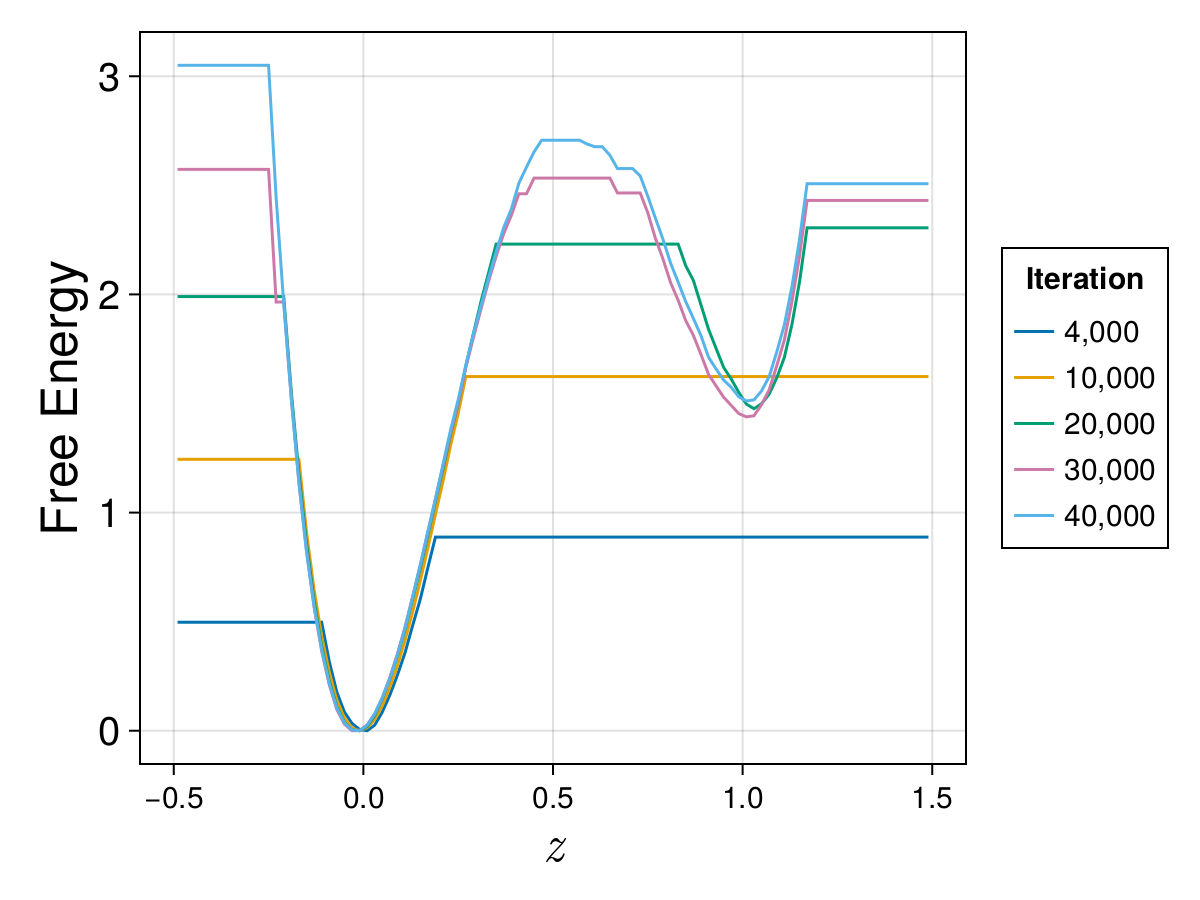}
    \caption{Free energy at different iterations of the simulation.}
    \label{fig:adaptive_FE}
  \end{subfigure}
  \caption{Behavior of the free energy and the values of the collective variable during a sample trajectory using Algorithm~\ref{alg:adaptive_MALA}.}
  \label{fig:adaptive_behavior}
\end{figure}

We present the results in Figure~\ref{fig:MALA}. The optimal value obtained is~$\alpha_{\rm opt}=1.5$, with~$\widehat{\tau}(\alpha_{\rm opt})^{\star}=672$, while~$\widehat{\tau}_{\rm cst}^{\star}=1402$. The comments on the results are the same as the ones in Section~\ref{sec:MALA_numerics}. This shows that the same efficiency is obtained with the adaptive scheme. This is expected, as the learned quantities usually rapidly converge. This is depicted in Figure~\ref{fig:adaptive_behavior}: we run a sample trajectory using~$4\times10^{4}$ iterations with time step~$\Delta t=2.5\times10^{-3}$ and~$\alpha=1.5$ (we checked that similar results were obtained using other values of~$\alpha$). We plot in Figure~\ref{fig:adaptive_CV} the values of the reaction coordinate as a function of the number of iterations. We present in Figure~\ref{fig:adaptive_FE} the values of the free energy over the range~$[z_{\rm min},z_{\rm max}]$ for 5 values of the number of iterations~$n$. We observe that after a small number of iterations (relatively to the original dynamics), the learned diffusion coefficient allows for rapid transitions between the two metastable states. The free energy rapidly converges to its stationary state in the two metastable states. The convergence of the free energy in the vicinity of the transition state takes more time, since the process dashes from one state to the other, without spending time over the transition state.

\begin{figure}
  \centering
  \begin{subfigure}[t]{0.475\columnwidth}
    \includegraphics[width=\columnwidth]{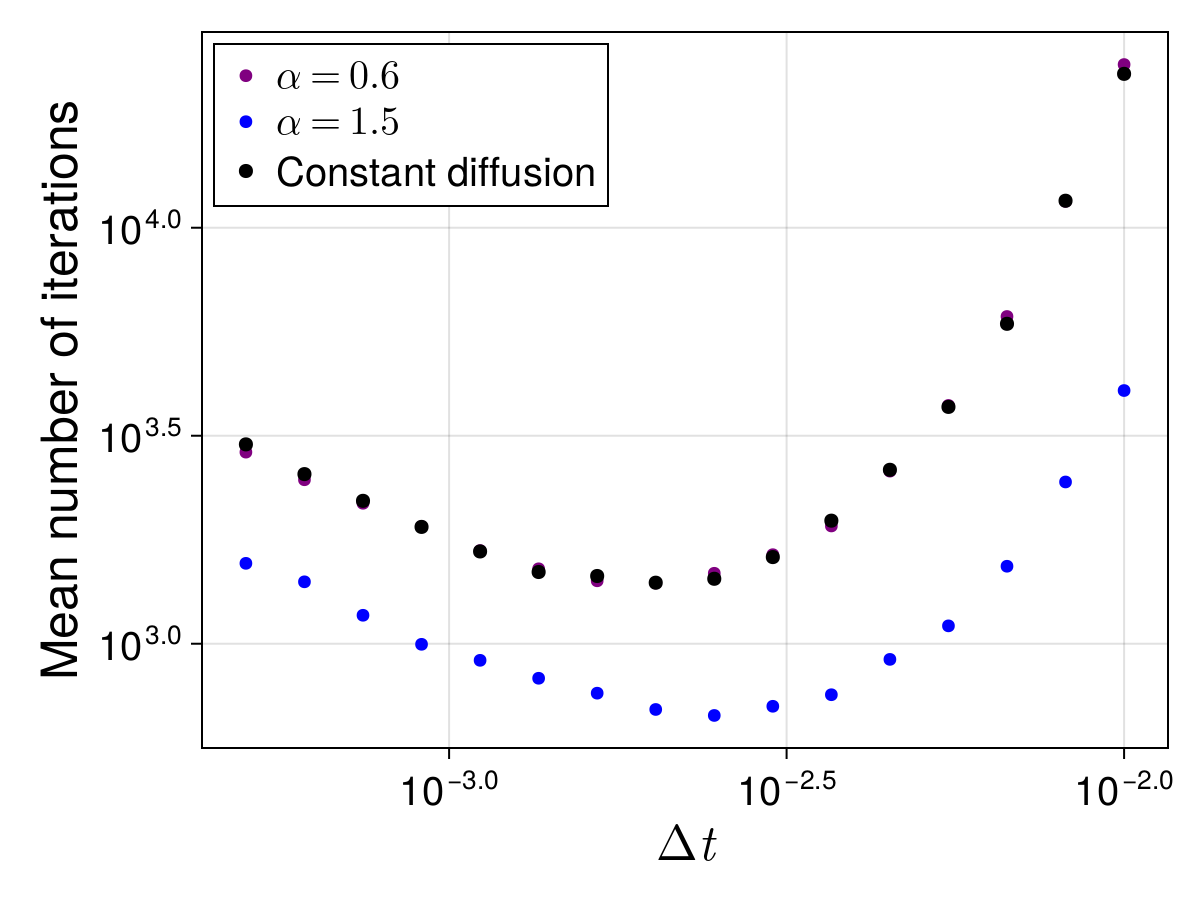}
    \caption{Mean number of iterations to observe a transition as a function of~$\Delta t$ for two values of~$\alpha$ and for the constant diffusion.}
    \label{fig:adaptive_some_values_n_iterations}
  \end{subfigure}\hfill
  \begin{subfigure}[t]{0.475\columnwidth}
    \centering
    \includegraphics[width=\columnwidth]{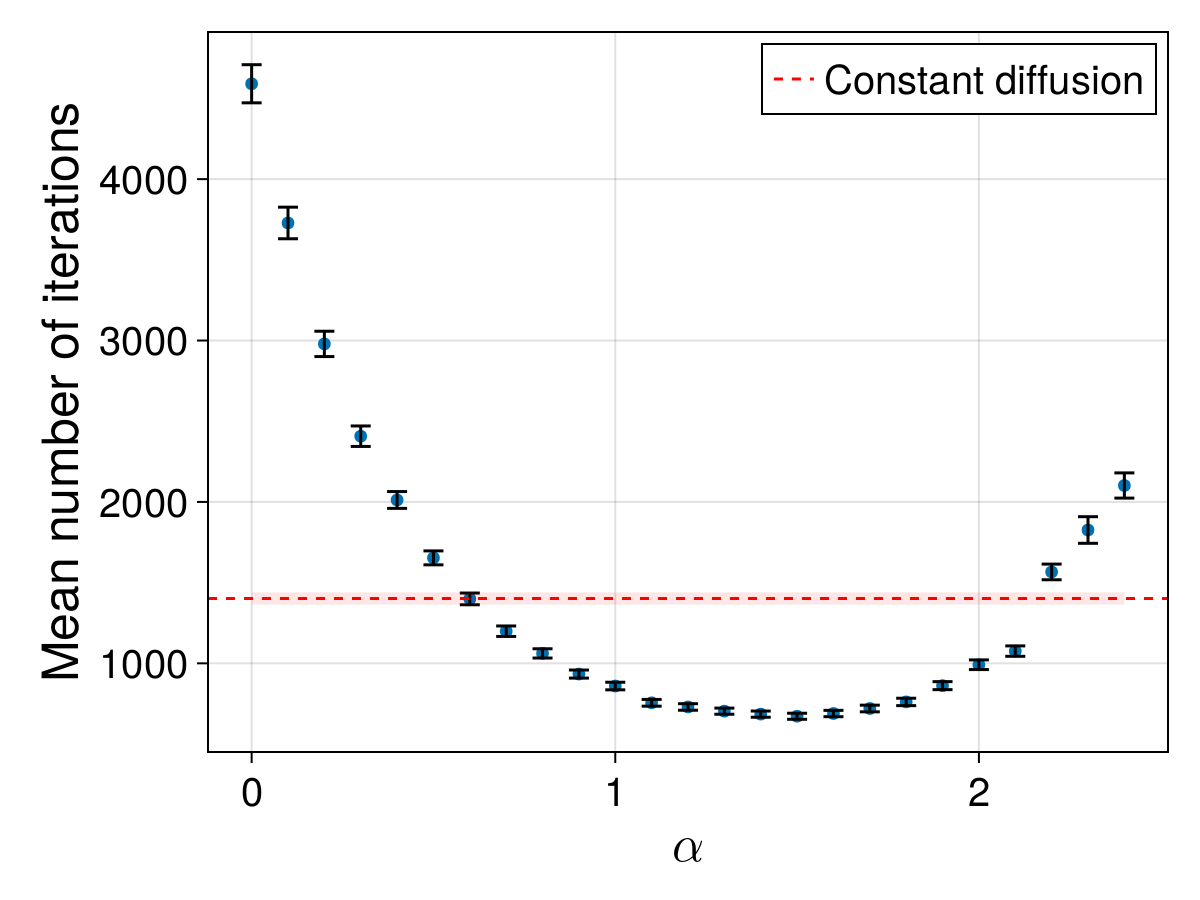}
    \caption{Minimum values of the transition times over the time steps as a function of~$\alpha$, and minimum transition time for the constant diffusion.}
    \label{fig:adaptive_minimum_values}
  \end{subfigure}
  \caption{Sampling efficiency using the Adaptive MALA algorithm~\ref{alg:adaptive_MALA}.}
  \label{fig:adaptive_MALA}
\end{figure}

\section{Optimizing the diffusion: the kinetic case}
\label{sec:HMC}

The objective of this section is to introduce position dependent diffusions for dynamics involving positions and momenta. In this context, the diffusion is actually the inverse of the mass.  In Section~\eqref{sec:GHMC}, we recall basics about  Generalized Hamiltonian Monte Carlo algorithms for position dependent diffusions. These algorithms are then introduced in details in Section~\ref{sec:algorithms}. We finally present numerical results in Section~\ref{sec:HMC_numerics}.

\subsection{Generalized Hamiltonian Monte Carlo algorithms for position dependent diffusions}\label{sec:GHMC}
First introduced in~\cite{duane_1987}, HMC algorithms generate a Markov chain~$(q^n,p^n)_{n\geqslant0}$ in the augmented space~$\calQ\times\bbR^{d}$ such that the marginal in position of the equilibrium probability distribution is the target Boltzmann-Gibbs probability distribution~$\pi$ defined in~\eqref{eq:gibbs_measure}. The additional variable~$p^n$ is a momentum associated with the position~$q^n$. HMC algorithms are built on three blocks: (i) a partial or full resampling of the momentum variable, (ii) the integration of Hamiltonian dynamics using the (possibly Generalized) Störmer--Verlet scheme~\cite{hairer_2006} and (iii) a Metropolis--Hastings accept/reject procedure. When the momentum variable is partially refreshed, the algorithm is referred to as a Generalized HMC algorithm~\cite{horowitz_1991}. HMC algorithms sample measures in the phase space of the form ~$\tilde Z^{-1} \rme^{-\beta H(q,p)} \rmd q \rmd p$ for some Hamiltonian function $H$ chosen such that the marginal in the position variable $q$ is the measure~$\pi$ defined by~\eqref{eq:gibbs_measure}. Here and in the following $\tilde Z=\int_{\calQ\times\bbR^{d}} \rme^{-\beta H(q,p)} \rmd q \rmd p$ denotes the normalizing constant in phase space, assumed to be finite.
Standard HMC algorithms rely on the Hamiltonian function
\begin{equation}
  \label{eq:standard_hamiltonian_function}
  H(q,p)=V(q)+\frac{1}{2}p^{\sfT}M^{-1}p,
\end{equation}
where~$M$ is a positive definite symmetric matrix which can be tuned to precondition the Hamiltonian dynamics.  For such a separable Hamiltonian, the Hamiltonian dynamics is integrated using the Störmer--Verlet (SV), or \textit{leapfrog}, numerical scheme.  This is an explicit scheme, and the associated numerical flow is defined as follows:~$(q^{n+1},p^{n+1})=\varphi_{\Delta t}^{\SV}(q^{n},p^{n})$ where
\begin{equation*}
  \left\lbrace
  \begin{aligned}
    p^{n+1/2}&=p^{n}-\frac{\Delta t}{2}\nabla V(q^{n}),\\
    q^{n+1}&=q^{n}+\Delta t M^{-1}p,\\
    p^{n+1}&=p^{n+1/2}-\frac{\Delta t}{2}\nabla V(q^{n+1}).
  \end{aligned}
  \right.
\end{equation*}
The flow~$\varphi_{\Delta t}^{\SV}$ is such that
\begin{equation}
  \label{eq:SV_property}
  S\circ\varphi_{\Delta t}^{\SV}\circ S\circ\varphi_{\Delta t}^{\SV} = {\rm id}_{\calQ},
\end{equation}
where~$S$ is the momentum reversal map:~$S(q,p)=(q,-p)$ and $\circ$ denotes the composition operator. Using the fact that $\varphi_{\Delta t}^{\SV}$ is an involution (see~\eqref{eq:SV_property}
 and that the flow $\varphi_{\Delta t}^{\SV}$ preserves the phase-space Lebesgue measure (it is actually symplectic), one can check that HMC samples $\tilde Z^{-1} \rme^{-\beta H(q,p)} \rmd q \rmd p$ (see for instance~\cite[Sections~2.1.3 and~2.1.4]{lelievre_2010}).

To introduce a nonconstant diffusion~$D$, one considers the Hamiltonian function
\begin{equation}
  \label{eq:RMHMC_hamiltonian_function}
  H(q,p)=V(q)-\frac{1}{2}\ln\det D(q)+\frac{1}{2}p^{\sfT}D(q)p.
\end{equation}
 One can check that the marginal in position of measure $\tilde Z^{-1} \rme^{-\beta H(q,p)} \rmd q \rmd p$ is the Boltzmann-Gibbs measure~$\pi$, whatever the diffusion matrix~$D$.
This Hamiltonian is for example  considered in the Riemann Manifold HMC (RMHMC) algorithm which was introduced in~\cite{girolami_2011}. Of course, \eqref{eq:RMHMC_hamiltonian_function} reduces to~\eqref{eq:standard_hamiltonian_function} when the diffusion is constant and equal to~$M^{-1}$. For well-chosen numerical parameters in the refreshment step of the momenta, the RMHMC 
algorithm can be shown to yield a (weakly) consistent discretization of the overdamped Langevin dynamics~\eqref{eq:overdamped_langevin_diffusion}, see~\cite[Section~3.3]{lelievre_2023_i}. This is why the inverse mass is denoted by $D$ in~\eqref{eq:RMHMC_hamiltonian_function}. The interest of using RMHMC instead of MALA is that the rejection probability is smaller in the limit of small timesteps (see~\cite[Section 3.3]{lelievre_2023_i}).

Since the Hamiltonian function~\eqref{eq:RMHMC_hamiltonian_function} is not separable, the integration of the Hamiltonian part is done using the Generalized St\"ormer Verlet (GSV) integrator, which is implicit. The associated numerical flow is formally defined as follows: for a given configuration~$(q^{n},p^{n})\in\calQ\times\bbR^{d}$,
\begin{equation}
  \label{eq:GSV}
  \left\lbrace
  \begin{aligned}
    p^{n+1/2}&=p^{n}-\frac{\Delta t}{2}\nabla_q H(q^{n},p^{n+1/2}),\\
    q^{n+1}&=q^{n}+\frac{\Delta t}{2}\left(\nabla_p H(q^{n},p^{n+1/2})+\nabla_p H(q^{n+1},p^{n+1/2})\right),\\
    p^{n+1}&=p^{n+1/2}-\frac{\Delta t}{2}\nabla_q H(q^{n+1},p^{n+1/2}).
  \end{aligned}
  \right.
\end{equation}
The first two steps are implicit and are typically solved in practice using Newton's method (see Appendix~\ref{app:RMHMC_implementation} as well as the discussion in Section~\ref{sec:HMC_numerics} for details of implementation and discussions of the computational cost). As discussed in~\cite{lelievre_2023_i}, this may introduce a bias, related to the fact that the implicit problem may have multiple solutions, or no solution, and the property~\eqref{eq:SV_property} may not hold for the numerical flow implemented in practice. To remedy this issue, reversibility checks have been introduced, first to sample probability measures on submanifolds~\cite{zappa_2018,lelievre_2019,lelievre_2022}, then in the context of RMHMC~\cite{graham_2022,noble_2023} and for general Hamiltonian functions~\cite{lelievre_2023_i}. Given a numerical flow~$\varphi_{\Delta t}^{\rm GSV}$ that attempts at integrating~\eqref{eq:GSV} (\emph{e.g.}~using Newton's method), one defines a \textit{numerical flow with reversibility checks} by
\begin{equation}
  \label{eq:GSV_rev}
  \forall (q,p)\in\calQ\times\bbR^{d},\qquad \varphi_{\Delta t}^{\rm rev}(q,p)=\varphi_{\Delta t}^{\rm GSV}(q,p)\,\mathbbm{1}_{(q,p)\in \calB_{\Delta t}}+S(q,p)\,\mathbbm{1}_{(q,p)\notin \calB_{\Delta t}},\\
\end{equation}
where the set~$\calB_{\Delta t}$ is defined as the set of configurations~$(q,p)\in\calQ\times\bbR^{d}$ such that the property~\eqref{eq:SV_property} holds, namely:
\begin{itemize}
  \item the Newton method for the \textit{forward} problem is converging so that~$\varphi_{\Delta t}^{\rm GSV}(q,p)$ is well-defined;
  \item the Newton method for the \textit{backward} problem, namely starting from~$(q,-p)$, is converging, so that~$\varphi_{\Delta t}^{\rm GSV}\circ S\circ\varphi_{\Delta t}^{\rm GSV}(q,p)$ is well-defined;
  \item the solutions to the forward and backward problems satisfy the reversibility property: $S \circ \varphi_{\Delta t}^{\rm GSV}\circ S\circ\varphi_{\Delta t}^{\rm GSV}(q,p)=(q,p)$.
\end{itemize}
It is then easy to check that the property~\eqref{eq:SV_property} holds for the flow with reversibility checks~$\varphi_{\Delta t}^{\rm rev}$. Moreover, one can prove that~$\varphi_{\Delta t}^{\rm rev}$ is a measurable map which preserves the Lebesgue measure over~$\calQ$. These two ingredients ensure that the RM(G)HMC algorithms are unbiased. We refer to~\cite{lelievre_2023_i} for precise definitions and results.  In the following, we refer to these algorithms with reversibility checks as RM(G)HMC even though in the literature, this denomination does not necessarily means that reversibility checks are used.

In our context, for the diffusion~\eqref{eq:optimal_diffusion_class}, the Hamiltonian function~\eqref{eq:RMHMC_hamiltonian_function} rewrites
\begin{equation}
  \label{eq:H_diffusion}
  H(q,p)=V(q)-\frac{d}{2}\ln\kappa_{\alpha}-\frac{1}{2}\ln a_{\alpha}(\xi(q))+\frac{\kappa_{\alpha}}{2}p^{\sfT}p+\frac{\kappa_{\alpha}}{2}(a_{\alpha}(\xi(q))-1)\frac{\left(\nabla\xi(q)\cdot p\right)^{2}}{\left\lVert\nabla\xi(q)\right\rVert^{2}}.
\end{equation}
The second term in~\eqref{eq:H_diffusion} is constant and can be omitted.

\subsection{Algorithms and reversibility checks}
\label{sec:algorithms}

We make precise in Section~\ref{sec:RMHMC} the RMHMC algorithm. We then describe in Section~\ref{sec:RMGHMC} the RMGHMC algorithm, which relies on the integration of an Ornstein--Uhlenbeck process to partially refresh the momenta. The specific implementation of these algorithms for the diffusion~\eqref{eq:optimal_diffusion_class} is made precise in Appendix~\ref{app:RMHMC_implementation}: it consists in integrating the Hamiltonian dynamics with Hamiltonian function~\eqref{eq:H_diffusion} using Newton's method. We also describe how these algorithms can be optimized in order to limit computational overheads, especially when the collective variable is a function of only a few components of the position variable.

\subsubsection{RMHMC}
\label{sec:RMHMC}

\begin{algorithm}
  \caption{RMHMC algorithm with reversibility checks.}
  \label{alg:RMHMC}
  Consider an initial condition~$(q^{0},p^{0})\in\calQ\times\bbR^{d}$, and set~$n=0$.
  \begin{enumerate}[label={[\thealgorithm.\roman*]}, align=left]
      \item \label{step:hmc_1} Sample~$\widetilde{p}^{n}\sim\calN(0,D(q^{n})^{-1})$;
      \item \label{step:hmc_2} Apply one step of the Hamiltonian dynamics with momentum reversal and~$S$-reversibility check:
      \begin{equation*}
          (\widetilde{q}^{n+1}, \widetilde{p}^{n+1})=S \circ \varphi_{\Delta t}^{\rm rev}(q^{n},\widetilde{p}^{n}),
      \end{equation*}
      where~$\varphi_{\Delta t}^{\rm rev}$ is defined by~\eqref{eq:GSV_rev};
      \item \label{step:hmc_3} Draw a random variable~$U^{n}$ with uniform law on~$[0,1]$:
      \begin{itemize}[label=$\bullet$]
          \item if~$U^{n}\leqslant \exp\left(\beta\left(H(q^{n},\widetilde{p}^{n})-H(\widetilde{q}^{n+1},\widetilde{p}^{n+1})\right)\right)$, accept the proposal and set $(q^{n+1},p^{n+1})=(\widetilde{q}^{n+1},\widetilde{p}^{n+1})$;
          \item else reject the proposal and set~$(q^{n+1},p^{n+1})=(q^{n},\widetilde{p}^{n})$;
      \end{itemize}
      \item Increment~$n$ and go back to~\ref{step:hmc_1}.
  \end{enumerate}
\end{algorithm}


The RMHMC algorithm with reversibility checks is detailed in Algorithm~\ref{alg:RMHMC}. Let us comment on the various steps of this algorithm. In Step~\ref{step:hmc_1}, the momenta is fully resampled according to the marginal in momenta of the measure $\tilde Z \rme^{-\beta H(q,p)} \rmd q \rmd p$. In view of~\eqref{eq:RMHMC_hamiltonian_function}, one easily sees that it is a Gaussian probability distribution with zero mean and covariance matrix~$D(q)^{-1}$. In Step~\ref{step:hmc_2}, the Hamiltonian dynamics is integrated using the modified GSV flow with reversibility checks~\eqref{eq:GSV_rev}. Lastly, Step~\ref{step:hmc_3} is a Metropolis--Hastings accept/reject procedure. Note that the sign of the momentum variable is not relevant when accepting or rejecting a proposal as the momenta are fully resampled at each iteration. 

\subsubsection{RMGHMC}
\label{sec:RMGHMC}

This section is based on~\cite[Section~3.2.1]{lelievre_2023_i}. The GHMC algorithm is related to a discretization of the Langevin dynamics, which writes, for the Hamiltonian~\eqref{eq:RMHMC_hamiltonian_function},
\begin{equation}
  \label{eq:langevin_dynamics}
  \left\lbrace
  \begin{aligned}
    \rmd q_t &= D(q_t)p_t\,\rmd t,\\
    \rmd p_t &= -\nabla_q H(q_t,p_t)\,\rmd t-\gamma D(q_t)p_t\,\rmd t+\sqrt{\frac{2\gamma}{\beta}}\,\rmd W_t,
  \end{aligned}
  \right.
\end{equation}
where~$\gamma>0$ is a friction parameter, and~$(W_t)_{t\geqslant 0}$ is a standard~$d$-dimensional Brownian motion. The GHMC algorithm is built using a splitting technique, considering separately the Hamiltonian and fluctuation-dissipation parts of the dynamics. Here, the fluctuation-dissipation part of the dynamics is simply an Ornstein--Uhlenbeck process, inducing a partial refreshment of the momenta instead of a full resampling as in the RMHMC algorithm~\ref{alg:RMHMC}. To integrate the Ornstein--Uhlenbeck process for a fixed position~$q^{n}$, namely
\begin{equation}
  \label{eq:OU}
  \rmd p_t = -\gamma D(q^{n})p_t\,\rmd t+\sqrt[]{\frac{2\gamma}{\beta}}\,\rmd W_t,
\end{equation}
we use a midpoint Euler scheme: for a time increment~$\Delta t/2$, this scheme reads
\begin{equation*}
  p^{n+1/4}=p^{n}-\frac{\Delta t}{4}\gamma D\left(q^{n}\right)\left(p^{n}+p^{n+1/4}\right)+\sqrt{\gamma\beta^{-1}\Delta t}\,\rmG^{n},
\end{equation*}
which yields
\begin{equation}
  \label{eq:OU_update_rmhmc}
    p^{n+1/4}= \left[
        \rmI_{d} + \frac{\Delta t}{4}\gamma D\left(q^{n}\right)
    \right]^{-1}\left[
        \left(\rmI_{d}-\frac{\Delta t}{4}\gamma D\left(q^{n}\right)\right)p^{n}+\sqrt{\gamma\beta^{-1}\Delta t}\,\rmG^{n}
    \right].
\end{equation}
This defines the map~$\varphi_{\Delta t/2}^{\rm OU}$:
\begin{equation*}
  \varphi_{\Delta t/2}^{\rm OU}(q^{n},p^{n},\rmG^{n})=(q^{n},p^{n+1/4}).
\end{equation*}
A Strang splitting based on these elements leads to Algorithm~\ref{alg:RMGHMC}. The main difference with Algorithm~\ref{alg:RMHMC} is that the momenta are not discarded, but partially refreshed at Steps~\ref{step:ghmc_1} and~\ref{step:ghmc_5}. Notice that the momenta are reversed at Step~\ref{step:ghmc_4}: this is crucial to get a consistent discretization of the Langevin dynamics~\eqref{eq:langevin_dynamics} in the small $\Delta t$ regime. More precisely, the momentum flip $S$ in Step~\ref{step:ghmc_2}, and the momentum flip at Step~\ref{step:ghmc_4} cancel if there is no rejection, neither in Step~\ref{step:ghmc_2} (\emph{i.e.}~$\varphi_{\Delta t}^{\rm rev}(q^{n},p^{n+1/4})=\varphi_{\Delta t}^{\rm GSV}(q^{n},p^{n+1/4})$) nor in Step~\ref{step:ghmc_3}.

\begin{algorithm}
  \caption{RMGHMC algorithm with reversibility checks.}
  \label{alg:RMGHMC}
  Consider an initial condition~$(q^{0},p^{0})\in\calQ\times\bbR^{d}$, and set~$n= 0$.
  \begin{enumerate}[label={[\thealgorithm.\roman*]}, align=left]
      \item\label{step:ghmc_1} Evolve the momenta by integrating the fluctuation-dissipation part with time increment~$\Delta t/2$:~$(q^{n},p^{n+1/4})=\varphi_{\Delta t/2}^{\rm OU}(q^{n},p^{n},\rmG^n)$.
      \item\label{step:ghmc_2} Apply one step of the Hamiltonian dynamics with momentum reversal and~$S$-reversibility check:
      \begin{equation*}
          (\widetilde{q}^{n+1}, \widetilde{p}^{n+3/4})=S\circ\varphi_{\Delta t}^{\rm rev}(q^{n},p^{n+1/4}),
      \end{equation*}
      where~$\varphi_{\Delta t}^{\rm rev}$ is defined in~\eqref{eq:GSV_rev};
      \item\label{step:ghmc_3} Draw a random variable~$U^{n}$ with uniform law on~$(0,1)$:
      \begin{itemize}[label=$\bullet$]
          \item if~$U^{n}\leqslant \exp(-H(\widetilde{q}^{n+1},\widetilde{p}^{n+3/4})+H(q^{n},p^{n+1/4}))$, accept the proposal and set $(q^{n+1},p^{n+3/4})=(\widetilde{q}^{n+1},\widetilde{p}^{n+3/4})$;
          \item else reject the proposal and set~$(q^{n+1},p^{n+3/4})=(q^{n},p^{n+1/4})$.
      \end{itemize}
      \item\label{step:ghmc_4} Reverse the momenta:~$\widetilde{p}^{n+1}=-p^{n+3/4}$.
      \item\label{step:ghmc_5} Evolve the momenta by integrating the fluctuation-dissipation part with time increment~$\Delta t/2$:~$(q^{n+1},p^{n+1})=\varphi_{\Delta t/2}^{\rm OU}(q^{n+1},\widetilde{p}^{n+1},\rmG^{n+1/2})$.
      \item Increment~$n$ and go back to~\ref{step:ghmc_1}.
  \end{enumerate}
\end{algorithm}


\subsection{Numerical results}
\label{sec:HMC_numerics}

The numerical experiment is the same as in Section~\ref{sec:MALA_numerics}. The values of~$\alpha$ are kept the same, but the time interval over which the optimal time step is searched for is scheme dependent: we choose 16 values of the time step evenly spaced log-wise on the interval~$[\Delta t_{\min},\Delta t_{\max}]$ with
\begin{itemize}
  \item for RMHMC:~$\Delta t_{\min}=3\times10^{-2}$ and~$\Delta t_{\max}=10^{-1}$;
  \item for RMGHMC:~$\Delta t_{\min}=5\times10^{-3}$ and~$\Delta t_{\max}=5\times10^{-2}$.
\end{itemize}
These intervals are chosen so that, for each scheme, an optimal timestep can be identified in terms of the mean transition duration, see below.

\paragraph{Hyperparameters.} Hyperparameters for Newton's algorithm are given in Appendix~\ref{app:RMHMC_implementation}. For the RMGHMC algorithm, we set~$\gamma=1$. This hyperparameter could be made position-dependent and be optimized as well, as was done for instance in~\cite{chak_2023}. This is left for future works. The uninformed choice~$\gamma=1$ can be seen as a balance between the limit~$\gamma\to0$, where the Langevin dynamics reduces to the Hamiltonian dynamics, and~$\gamma\to+\infty$, where the Langevin dynamics converge to the overdamped Langevin dynamics~\eqref{eq:overdamped_langevin_diffusion} (see~\cite{hottovy_2015,wang_2024}). Actually, hypocoercivity results on weighted~$L^{2}$ spaces for the continuous-in-time Langevin dynamics show that the rate of convergence towards equilibrium is lower bounded by~$\min\left(\gamma,\gamma^{-1}\right)$ (up to a multiplicative constant), see for instance~\cite{dolbeault_2012,grothaus_2016,roussel_2018}, so that an optimal non zero and finite value for $\gamma$ exists.


\begin{figure}
  \centering
  \begin{subfigure}{0.475\columnwidth}
    \includegraphics[width=\columnwidth]{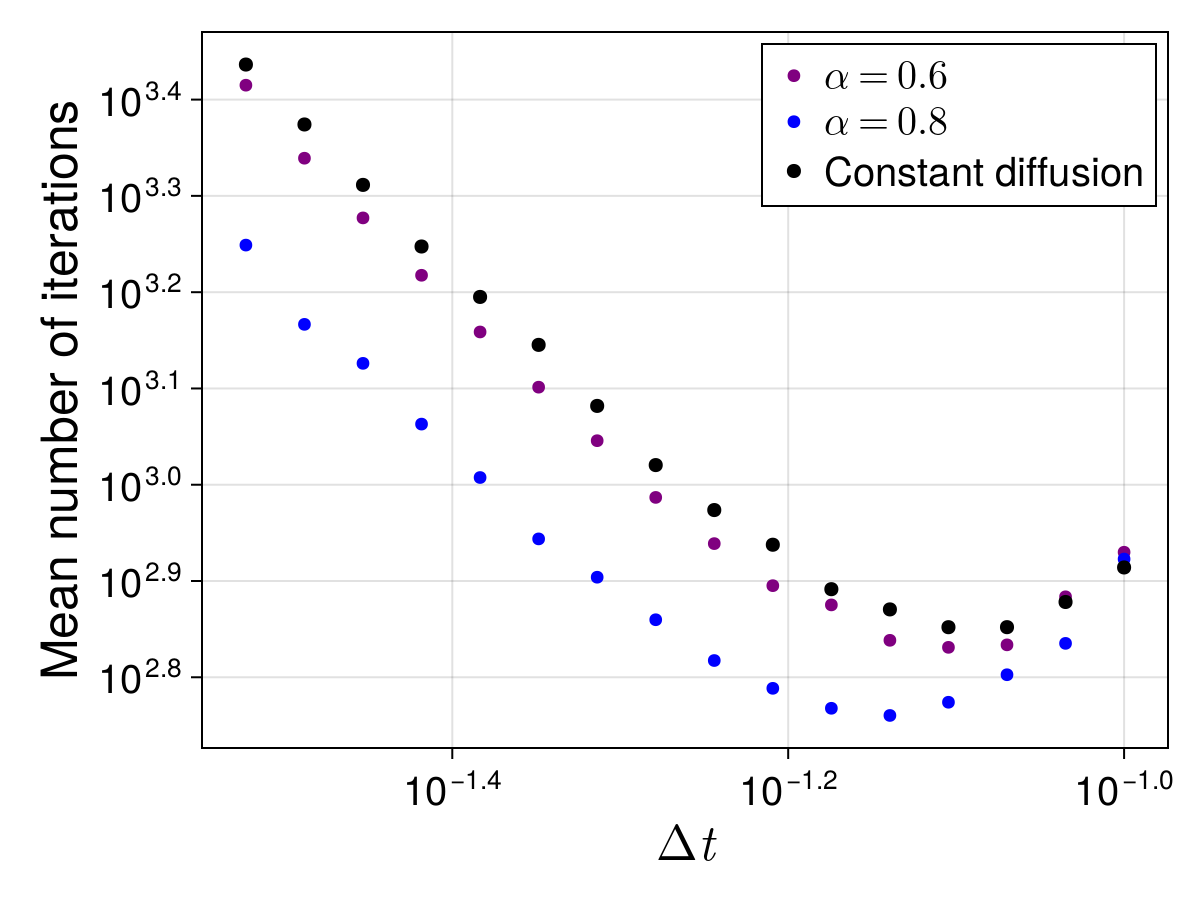}
    \caption{Mean number of iterations to observe a transition as a function of~$\Delta t$ for some values of~$\alpha$, and for the constant diffusion.}
    \label{fig:RMHMC_some_values_n_iterations}
  \end{subfigure}\hfill
  \begin{subfigure}{0.475\columnwidth}
    \centering
    \includegraphics[width=\columnwidth]{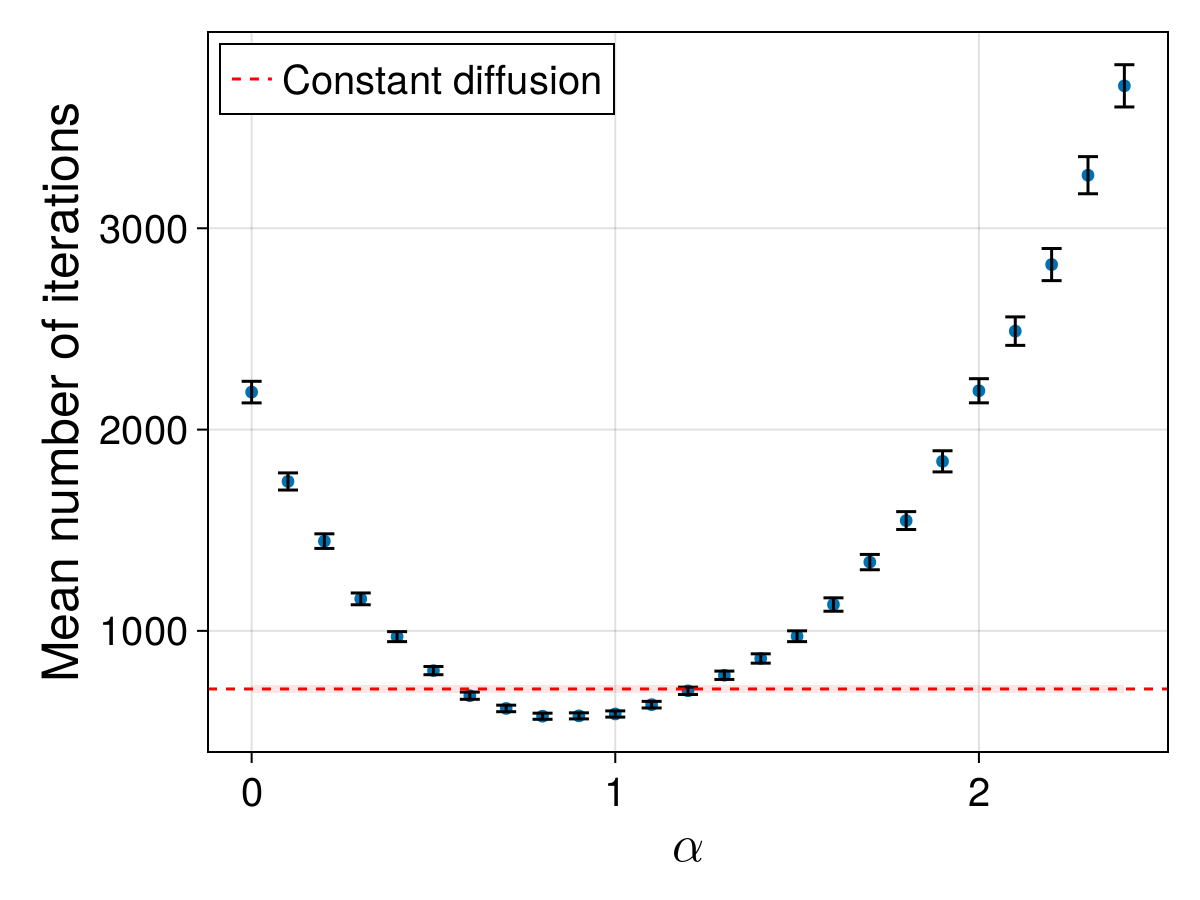}
    \caption{Minimum values of the transition times over the time steps as a function of~$\alpha$, and minimum transition time for the constant diffusion.}
    \label{fig:RMHMC_minimum_values}
  \end{subfigure}
  \caption{Sampling efficiency using the RMHMC algorithm~\ref{alg:RMHMC}.}
  \label{fig:RMHMC}
\end{figure}

\begin{figure}
  \centering
  \begin{subfigure}{0.475\columnwidth}
    \includegraphics[width=\columnwidth]{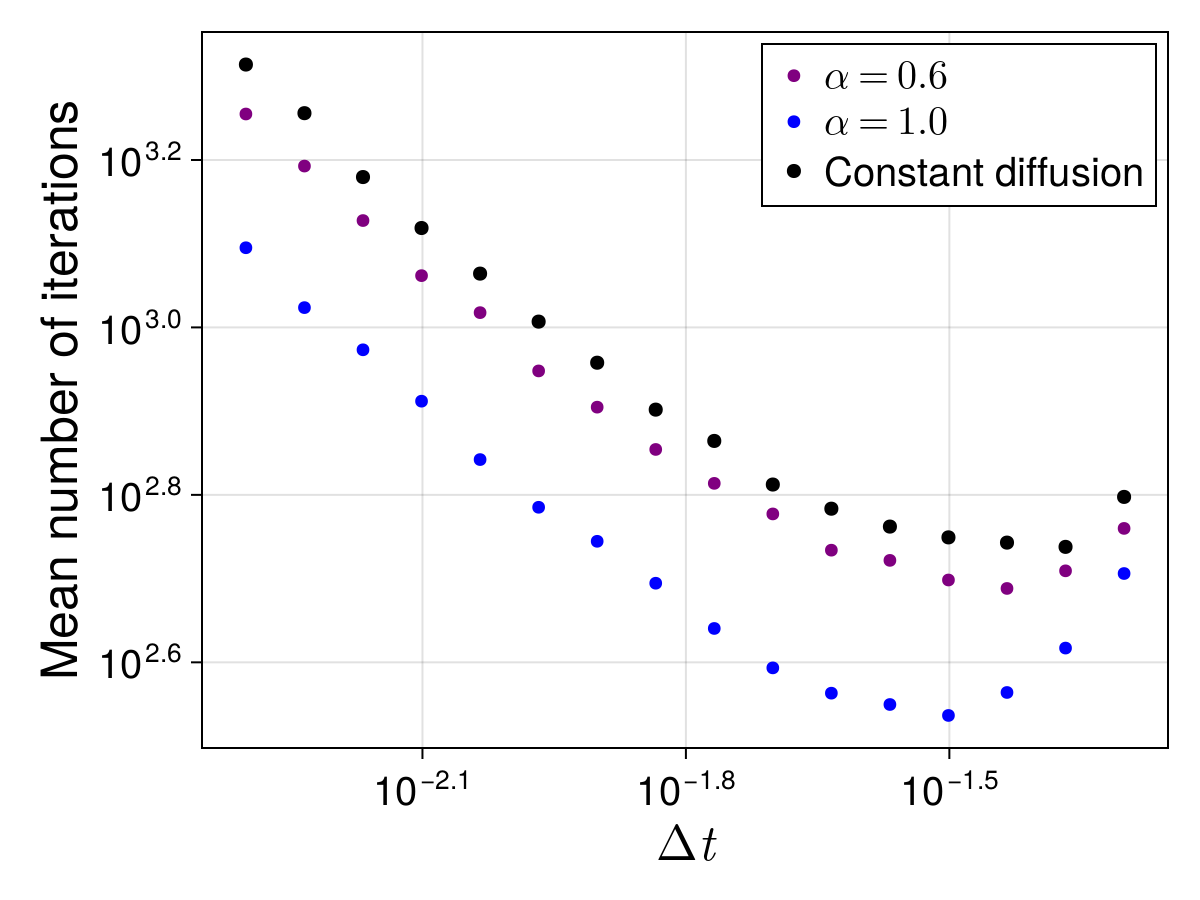}
    \caption{Mean number of iterations to observe a transition as a function of~$\Delta t$ for some values of~$\alpha$, and for the constant diffusion.}
    \label{fig:RMGHMC_some_values_n_iterations}
  \end{subfigure}\hfill
  \begin{subfigure}{0.475\columnwidth}
    \centering
    \includegraphics[width=\columnwidth]{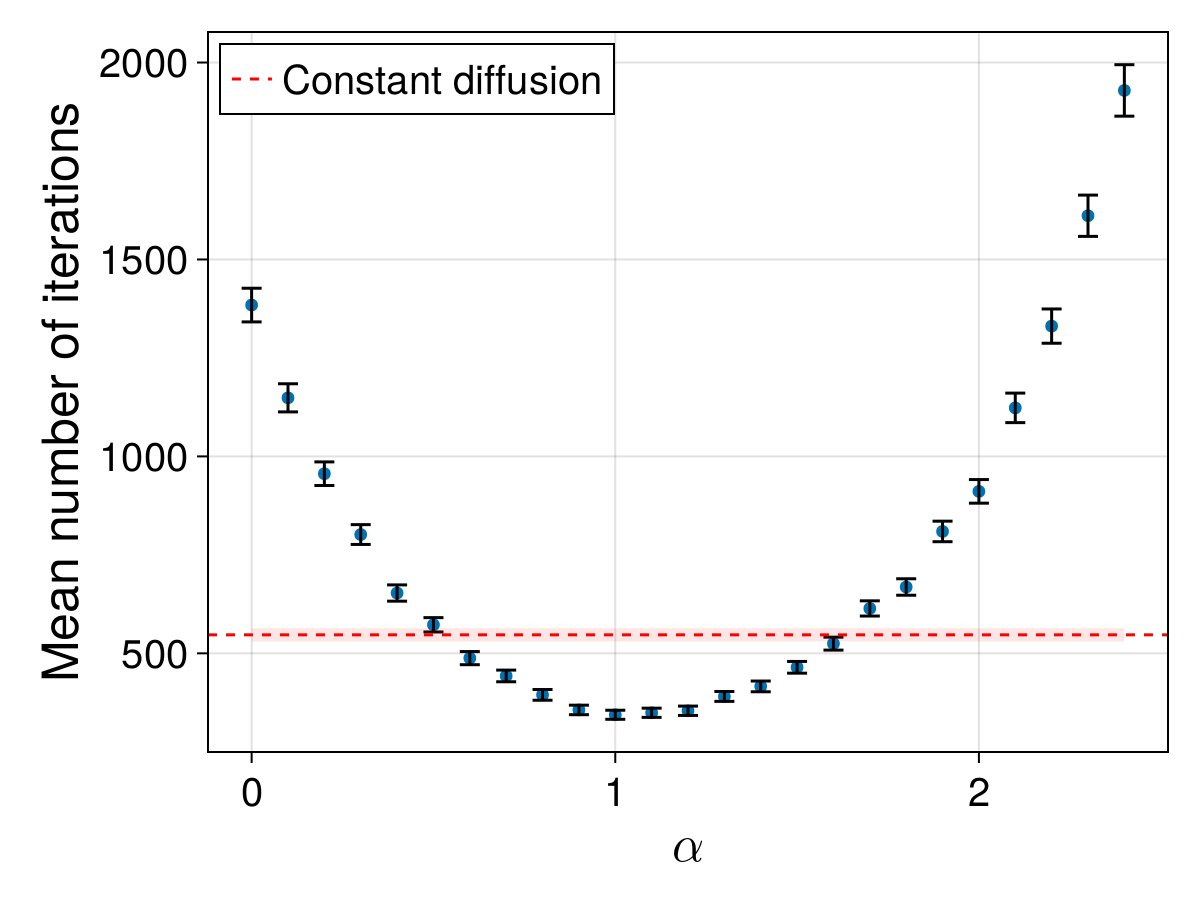}
    \caption{Minimum values of the transition times over the time steps as a function of~$\alpha$, and minimum transition time for the constant diffusion.}
    \label{fig:RMGHMC_minimum_values}
  \end{subfigure}
  \caption{Sampling efficiency using the RMGHMC algorithm~\ref{alg:RMGHMC}.}
  \label{fig:RMGHMC}
\end{figure}

\paragraph{Mean number of iterations for transitions.} The results using the RMHMC (respectively the RMGHMC) algorithm are presented in Figure~\ref{fig:RMHMC} (respectively in Figure~\ref{fig:RMGHMC}). For RMHMC, we obtain~$\widehat{\tau}(\alpha_{\rm opt})^{\star}=576$ for~$\alpha_{\rm opt}=0.8$, while~$\widehat{\tau}_{\rm cst}^{\star}=711$. For RMGHMC, we obtain~$\widehat{\tau}(\alpha_{\rm opt})^{\star}=344$ for~$\alpha_{\rm opt}=1.0$, while~$\widehat{\tau}_{\rm cst}^{\star}=547$. The results for RMHMC are comparable (though slightly better) to what has been obtained for MALA in Sections~\ref{sec:MALA_numerics} and~\ref{sec:adaptive_MALA_numerics}. The results for RMGHMC are significantly better, which shows the relevance of building inertial samplers in phase space and performing the very modest computational overheads related to the Ornstein--Uhlenbeck dynamics~\eqref{eq:OU}. For both algorithms, we observe a great reduction in the mean transition durations when choosing multiplicative noise instead of additive noise. In practice, is seems that a good rule of thumb to set the value of~$\alpha$ is to make the argument of the exponential in~\eqref{eq:optimal_diffusion_class} of the order of~1 in low free energy regions. 

Let us emphasize that the Newton steps involved in the RM(G)HMC do not significantly increase  the cost of one iteration compared with MALA for two reasons. First, the Newton steps do not require the computation of new forces. Second, for collective variables which involve only a few degrees of freedom, the matrices needed in the Newton updates are the identity up to low-dimensional block matrices.

Note that the mean transition duration scales, in the limit~$\Delta t\to0$, as~$\Delta t^{-1}$ for the RMGHMC algorithm while it scales as~$\Delta t^{-2}$ for the RMHMC algorithm. This is consistent with the fact that
\begin{itemize}
  \item RMGHMC yields a consistent discretization of the Langevin dynamics~\eqref{eq:langevin_dynamics} (see~\cite[Section~3.2.2]{lelievre_2023_i}), so that the physically simulated time is~$t_{\rm sim}=N_{\rm iter}\Delta t$;
  \item RMHMC yields a weakly consistent discretization of the overdamped Langevin dynamics~\eqref{eq:overdamped_langevin_diffusion} with effective time step~$h=\Delta t^{2}/2$ (see~\cite[Section~3.3]{lelievre_2023_i}), so that the physically simulated time is~$t_{\rm sim}=N_{\rm iter}h$.
\end{itemize}
We observe that the optimal time step for the RMHMC algorithm~$\Delta t_{\rm RMHMC}^{\star}$ is related to the optimal time step for the RMGHMC algorithm~$\Delta t_{\rm RMGHMC}^{\star}$ as~$\Delta t_{\rm RMHMC}^{\star}\approx\sqrt{\Delta t_{\rm RMGHMC}^{\star}}$ (which is somewhat in line with the consistency results mentioned above). This is why, in order to obtain an optimal time step in terms of the mean number of iterations for transitions, we chose different interval of values for RMHMC and RMGHMC. Let us mention that the extra computations incurred by the integration of the Ornstein--Uhlenbeck process~\eqref{eq:OU} in the RMGHMC algorithm are negligible compared to gain in terms of performance with respect to RMHMC.

\paragraph{Rejection probabilities.}
When the time step~$\Delta t$ gets larger, it is expected that rejections because of nonconvergence when solving the implicit problems increase. Let us quantify that by computing the various rejection probabilities when using the RM(G)HMC algorithms. We computed the rejection probabilities for each value of the time step used in the numerical experiment above (\emph{i.e.}~using 16 values evenly spaced log-wise in the interval~$[\Delta t_{\rm min},\Delta t_{\rm max}]$), and also for smaller values of the time step in order to analyze the behaviours in the small $\Delta t$ regime. For the RMHMC algorithm, these smaller values are
\begin{equation*}
  \left\lbrace 10^{-4}, 2\times10^{-4}, 3\times10^{-4}, \dots,9\times10^{-4}, 10^{-3}, 2\times10^{-3},\dots,2\times10^{-2}\right\rbrace,
\end{equation*}
while they are chosen equal to
\begin{equation*}
  \left\lbrace 10^{-5},2\times 10^{-5},3\times 10^{-5},\dots,9\times10^{-5},10^{-4},2\times10^{-4},\dots,4\times10^{-3}\right\rbrace,
\end{equation*}
for the RMGHMC algorithm. For each value of the time step~$\Delta t$, we run the RM(G)HMC algorithm for~$10^{7}$ iterations, starting from~$q^{0}$ defined in~\eqref{eq:q0}. We then count how many rejections were due to
\begin{itemize}
  \item not being able to solve for~$p^{n+1/2}$ in~\eqref{eq:GSV} in the \textit{forward} pass (labeled `Forward (momenta)');
  \item having succeeded but not being able to solve for~$q^{n+1}$ in~\eqref{eq:GSV} in the \textit{forward} pass (labeled `Forward (position)');
  \item having succeeded but not being able to solve for~$p^{n+1/2}$ in~\eqref{eq:GSV} in the \textit{backward} pass (labeled `Backward (momenta)');
  \item having succeeded but not being able to solve for~$q^{n+1}$ in~\eqref{eq:GSV} in the \textit{backward} pass (labeled `Backward (position)');
  \item having succeeded but not observing numerical reversibility~\eqref{eq:reversibiliy_check} (labeled `Reversibility');
  \item having succeeded but rejecting because of the Metropolis--Hastings accept/reject procedure (labeled `Metropolis--Hastings').
\end{itemize}
We then divide each number by the number of tries to obtain probabilities. 
Let us emphasize that the rejection event can be written as the union of these 6 disjoint rejection events.

\begin{figure}  
  \centering
  \begin{subfigure}{0.9\columnwidth}
    \centering
    \includegraphics[width=\columnwidth]{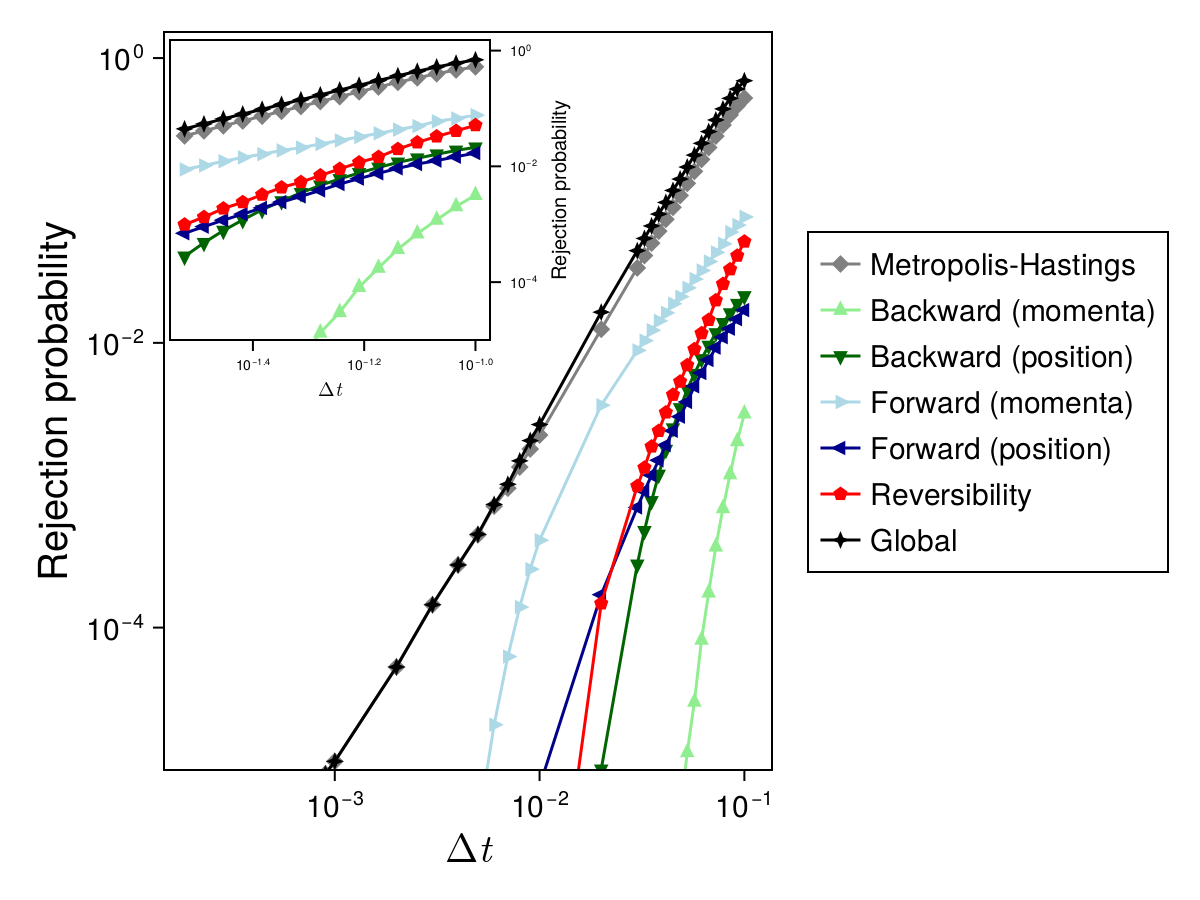}
    \caption{Rejection probabilities for~$\alpha=0.8$ for the RMHMC algorithm.}
    \label{fig:rejection_probabilities_RMHMC}
  \end{subfigure}\\
  \begin{subfigure}{0.9\columnwidth}
    \centering
    \includegraphics[width=\columnwidth]{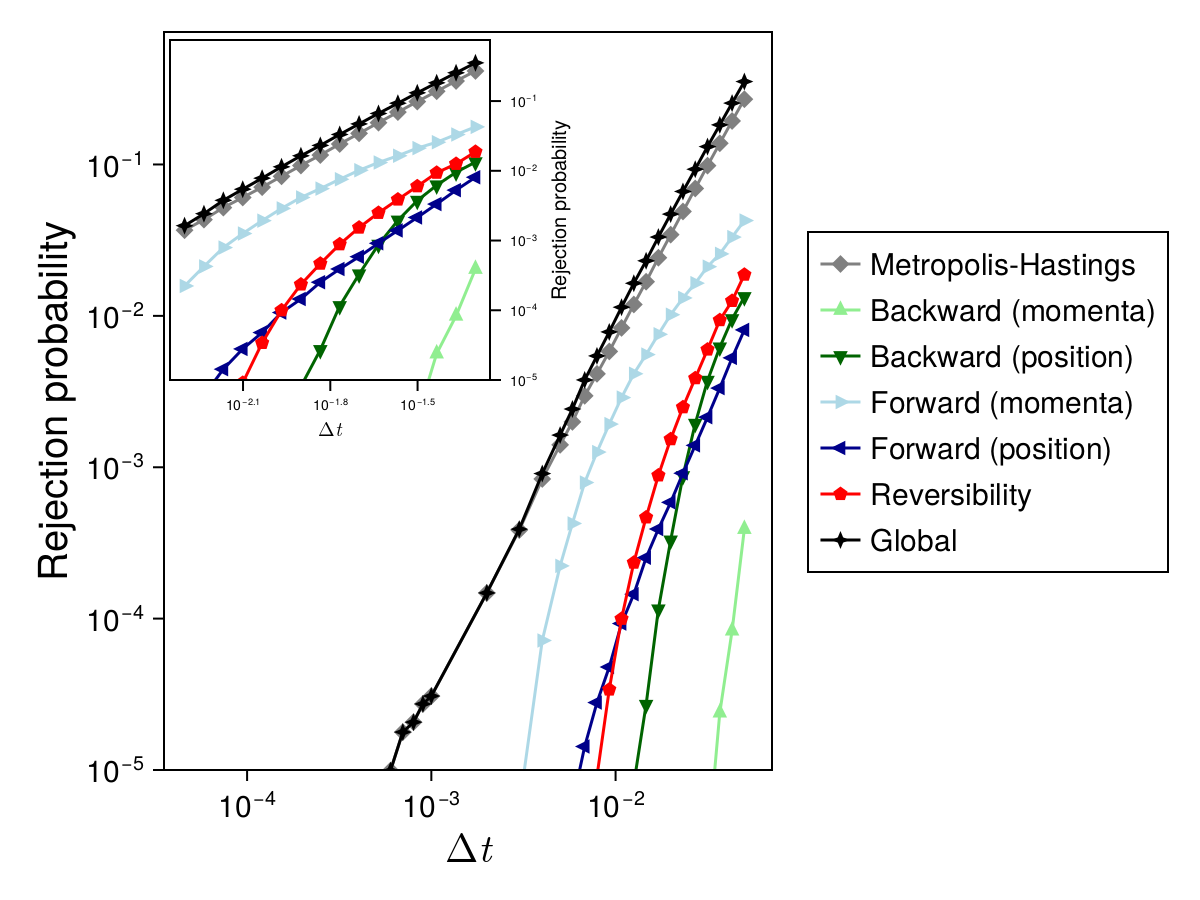}
    \caption{Rejection probabilities for~$\alpha=1.0$ for the RMGHMC algorithm.}
    \label{fig:rejection_probabilities_RMGHMC}
  \end{subfigure}
  \caption{Rejection probabilities associated to the RM(G)HMC algorithms for their respective optimal values of~$\alpha$.}
  \label{fig:rejection_probabilities}
\end{figure}

We show the results in Figure~\ref{fig:rejection_probabilities}, where only probabilities larger than~$10^{-5}$ are shown. The inset presents the results for the interval of time steps used in the numerical experiment, while the main axis shows the rejection probabilities for all values of the time step. Figure~\ref{fig:rejection_probabilities_RMHMC} corresponds to the RMHMC algorithm, while Figure~\ref{fig:rejection_probabilities_RMGHMC} corresponds to RMGHMC. In both cases, the optimal value for~$\alpha$ found above was used. We checked that similar behavior was observed for other values of~$\alpha$. 

In both cases, we observe that rejections due to nonreversibility issues happen for all values of the time step used for the numerical experiment above, in particular for the optimal time step. This shows the relevance of implementing reversibility checks in order to perform an unbiased sampling, especially when local exploration is optimized by adjusting the time step. Notice that rejection probabilities because of these nonreversibility issues eventually vanish in the small time step regime. This is consistent with other numerical experiments that were performed in~\cite{lelievre_2023_i}: the rejection probabilities seem to converge to~0 faster than the rejection probability of the sole Metropolis--Hastings accept/reject procedure, which scales as~$\rmO(\Delta t^{3})$. This also numerically confirms that RMHMC is a weakly consistent discretization of the overdamped Langevin dynamics~\eqref{eq:overdamped_langevin_diffusion}, see~\cite[Proposition~11]{lelievre_2023_i}.

We give in Table~\ref{tab:rejection_probabilities} the decomposition of the rejection probabilities for various values of~$\alpha$ and associated optimal time step for the RMHMC and RMGHMC algorithms. Note that, for both the RMHMC and RMGHMC algorithms, the optimal time step slightly decreases as~$\alpha$ increases (in particular when its value is larger than~1), even though we normalize the diffusion using the constant~$\kappa_\alpha$ in~\eqref{eq:optimal_diffusion_class}. This is related to the fact that~$\kappa_\alpha$ rapidly decreases to~0, since the integrand in the definition of~$\kappa_\alpha$ (see~\eqref{eq:kappa_alpha}) is equivalent to~$\rme^{(\alpha-1)\beta F}/\sigma$ when~$\alpha\to+\infty$, so that the normalization constant of the diffusion, which is the inverse of~$\kappa_\alpha$, rapidly converges to~0 when the value of~$\alpha$ is larger than 1. We checked that optimal time steps for values of~$\alpha$ below~1 (\emph{i.e.}~before the rapid convergence towards~0 happens) are similar.

We observe an increase in rejection probabilities because of nonreversibility issues when~$\alpha$ increases. This means that increasing the range of values of the diffusion (here, by increasing the value of~$\alpha$) has a great impact on nonreversibility, even when the time step is adjusted. It therefore advocates to keep track of the number of rejections because of nonreversibility issues during the simulation, for instance when targeting a specific acceptance rate, and to bound the diffusion values (using additional hyperparameters chosen by the practitioner depending on the application) when these rejections become too much of an issue.



Of course, the results presented here depend on our specific test problems, as well as on the implementation of the GSV solver (and to some extent on the way the Ornstein--Uhlenbeck is solved), and in particular on the choice of the hyperparameters used for Newton's algorithm, see Appendix~\ref{app:RMHMC_implementation}. An extensive numerical study on the effect of hyperparameters on the performance of RMHMC is given in~\cite{brofos_2021}.

\begin{table}
  \centering
  \hrule height 0.028cm
  \bigskip
  \begin{center}
    RMHMC algorithm
  \end{center}
  \resizebox{\textwidth}{!}{%
  \begin{tabular}{ccccccccc}
    \toprule
    \multirow{2}{*}[-0.1cm]{$\alpha$}& \multirow{2}{*}[-0.1cm]{$\Delta t$}& \multicolumn{2}{c}{Forward} & \multicolumn{2}{c}{Backward} & \multirow{2}{*}[-0.1cm]{Reversibility} & \multirow{2}{*}[-0.1cm]{Metropolis--Hastings} & \multirow{2}{*}[-0.1cm]{Global}\\
    \cmidrule(lr){3-4} \cmidrule(lr){5-6}
    && Momenta & Position & Momenta & Position\\\midrule
    $0.6$ & $7.86\times10^{-2}$ & $2.7\times10^{-2}$ & $4.8\times10^{-3}$ & $3.0\times10^{-6}$ & $4.9\times10^{-3}$ & $8.1\times10^{-3}$ & $0.35$ & $0.40$\\
    $0.8$ & $7.25\times10^{-2}$ & $4.3\times10^{-2}$ & $9.3\times10^{-3}$ & $3.7\times10^{-4}$ & $1.2\times10^{-2}$ & $2.0\times10^{-2}$ & $0.28$ & $0.37$\\
    $1.0$ & $6.18\times10^{-2}$ & $5.5\times10^{-2}$ & $1.2\times10^{-2}$ & $1.7\times10^{-3}$ & $1.6\times10^{-2}$ & $3.0\times10^{-2}$ & $0.18$ & $0.30$\\
    $1.5$ & $4.48\times10^{-2}$ & $8.2\times10^{-2}$ & $3.4\times10^{-2}$ & $1.0\times10^{-2}$ & $1.2\times10^{-2}$ & $3.9\times10^{-2}$ & $7.5\times10^{-2}$ & $0.25$\\
    \bottomrule
  \end{tabular}%
  }\\[0.2cm]
  \begin{center}
    RMGHMC algorithm
  \end{center}
  \resizebox{\textwidth}{!}{%
  \begin{tabular}{ccccccccc}
    \toprule
    \multirow{2}{*}[-0.1cm]{$\alpha$}& \multirow{2}{*}[-0.1cm]{$\Delta t$}& \multicolumn{2}{c}{Forward} & \multicolumn{2}{c}{Backward} & \multirow{2}{*}[-0.1cm]{Reversibility} & \multirow{2}{*}[-0.1cm]{Metropolis--Hastings} & \multirow{2}{*}[-0.1cm]{Global}\\
    \cmidrule(lr){3-4} \cmidrule(lr){5-6}
    && Momenta & Position & Momenta & Position\\\midrule
    $0.6$ & $3.68\times10^{-2}$ & $7.2\times10^{-3}$ & $4.4\times10^{-4}$ & $0.0$ & $4.1\times10^{-5}$ & $3.3\times10^{-4}$ & $0.12$ & $0.13$ \\
    $1.0$ & $3.15\times10^{-2}$ & $2.1\times10^{-2}$ & $2.1\times10^{-3}$ & $2.9\times10^{-6}$ & $3.7\times10^{-3}$ & $6.0\times10^{-3}$ & $9.8\times10^{-2}$ & $0.13$\\
    $1.4$ & $1.71\times10^{-2}$ & $1.8\times10^{-2}$ & $1.6\times10^{-3}$ & $2.9\times10^{-5}$ & $3.8\times10^{-3}$ & $7.9\times10^{-3}$ & $3.7\times10^{-2}$ & $6.8\times10^{-2}$\\
    $2.0$ & $7.92\times10^{-3}$ & $1.2\times10^{-2}$ & $1.2\times10^{-3}$ & $8.6\times10^{-5}$ & $2.1\times10^{-3}$ & $7.8\times10^{-3}$ & $2.2\times10^{-2}$ & $4.5\times10^{-2}$\\
    \bottomrule
  \end{tabular}%
  }
  \caption{Decomposition of the rejection probabilities for various values of~$\alpha$ and their corresponding optimal time steps~$\Delta t$ for the RM(G)HMC algorithms. ``Forward'' refers to non-convergent iterations to solve the implicit forward problem; ``Backward'' refers to a convergence of the implicit forward problem but a non-convergent implicit backward problem; ``Reversibility'' refers to convergent implicit forward and backward problems for which~\eqref{eq:SV_property} is not satisfied; upon acceptance in the three previous steps, ``Metropolis--Hastings'' refers to a rejection in the acceptance/rejection of Step~\ref{step:hmc_3} (respectively Step~\ref{step:ghmc_3}) in the RMHMC (respectively RMGHMC) algorithm. Finally, ``Global'' is the global rejection probability, namely the sum of all the previous columns.}
  \label{tab:rejection_probabilities}
\end{table}

\section{Conclusion and perspectives}

We introduced a class of diffusion matrices~\eqref{eq:optimal_diffusion_class} that modifies the effective diffusion of overdamped Langevin dynamics in order for the effective dynamics to be governed by an effective diffusion coefficient leading to a fast convergence. We provided a complete description of the methodology for one-dimensional collective variables. We list below a few perspectives opened by this work:
\begin{itemize}
  \item The methodology can be adapted to multidimensional collective variables upon finding a rationale for choosing the function~$a$, which, in principle, is matrix-valued.
  \item For one-dimensional collective variables, there may be an interest in introducing an additional multiplicative parameter in~$a_\alpha$ in order to balance the relative contributions of the projectors~$P$ and~$P^{\perp}$ in the definition of~\eqref{eq:optimal_diffusion_class}, see Remark~\ref{rem:free_energy_constant}.
  \item Other hyperparameters, such as the friction parameter~$\gamma$ in the RMGHMC algorithm, could be optimized to enhance convergence towards equilibrium.
  \item Another possible extension would be to directly optimize the diffusion in the latent space instead of relying on an explicit formula obtained in the homogenized limit, as in the parametric ansatz~\eqref{eq:optimal_diffusion_class}. More precisely, on could use the methodology developed in~\cite{cui_2024,lelievre_2024}, at least if the latent space is low dimensional (up to dimension three). In this case, Finite Element methods could be used to directly optimize the effective diffusion by maximizing the spectral gap of the operator (with adjoints taken in~$L^{2}(\xi\star\pi)$)
  \begin{equation*}
    \calL_{\sigma_a^{2}}=-\beta^{-1}\partial_{z}^{\star}\sigma_a^{2}\partial_z,
  \end{equation*}
  with respect to the map~$a$, where~$\sigma_a$ is defined in~\eqref{eq:drift_noise_effective_dynamics}.
  \item Since the mean force and free energy are learned in order to build the diffusion~\eqref{eq:optimal_diffusion_class}, one could use it in ABF methods in order to bias the dynamics, leading (once the estimated free energy is converged) to overdamped Langevin dynamics of the form
  \begin{equation*}
    \rmd q_t=\left(
        -D_{\alpha}(q_t)\nabla\left(
        V(q_t)-\upsilon F(\xi(q_t))
        \right)+\beta^{-1}\div D_{\alpha}(q_t)
        \right)\rmd t+\sqrt{2\beta^{-1}}D_{\alpha}(q_t)^{1/2}\rmd W_t,
  \end{equation*}
  with~$\upsilon$ a positive scaling parameter to be chosen. Of course, this idea can also be applied to the Langevin dynamics as presented in Section~\ref{sec:HMC}. It would be interesting to understand how much can be gained in practice by blending these two ideas: free energy biasing, and introducing a position dependent diffusion.
  \item Lastly, a major assumption for this work is the knowledge of a \textit{global information} through the collective variable~$\xi$, that identifies a direction in which the diffusion is modulated. Many works have been devoted to the definition and practical construction of optimal collective variables. In particular, learning collective variables using machine learning methods has been considered to this end, see for instance the review papers~\cite{Ferguson2018,Gkeka2020,sidkyreview,glielmoreview,Chen2021}. Using these tools may be beneficial to further improve sampling.
\end{itemize}

\paragraph{Acknowledgements.} The works of T.L., R.S. and G.S. benefit from fundings from the European Research Council (ERC) under the European Union’s Horizon 2020 research and innovation programme (project EMC2, grant agreement No 810367), and from the Agence Nationale de la Recherche through the grants ANR-19-CE40-0010-01 (QuAMProcs) and ANR-21-CE40-0006 (SINEQ).

\bibliographystyle{abbrv}
\bibliography{biblio}

\begin{thebibliography}{10}

\bibitem{bakry_2014}
D.~Bakry, I.~Gentil, and M.~Ledoux.
\newblock {\em {Analysis and Geometry of Markov Diffusion Operators}}.
\newblock Grundlehren der mathematischen Wissenschaften, Vol. 348. {Springer},
  2014.

\bibitem{beskos_2011}
A.~Beskos, F.~Pinski, J.~Sanz-Serna, and A.~Stuart.
\newblock {Hybrid Monte Carlo} on {Hilbert} spaces.
\newblock {\em Stochastic Processes and their Applications},
  121(10):2201--2230, 2011.

\bibitem{bou-rabee_2021}
N.~Bou-Rabee and A.~Eberle.
\newblock Two-scale coupling for preconditioned {H}amiltonian {M}onte {C}arlo
  in infinite dimensions.
\newblock {\em Stoch. Partial Differ. Equ. Anal. Comput.}, 9(1):207--242, 2021.

\bibitem{bou-rabee_2018}
N.~Bou-Rabee and J.~M. Sanz-Serna.
\newblock Geometric integrators and the {Hamiltonian Monte Carlo} method.
\newblock {\em Acta Numerica}, 27:113–206, 2018.

\bibitem{brofos_2021}
J.~A. Brofos and R.~R. Lederman.
\newblock On numerical considerations for {R}iemannian {M}anifold {H}amiltonian
  {M}onte {C}arlo.
\newblock {\em arXiv preprint}, 2111.09995, 2021.

\bibitem{bronasco_2025}
E.~Bronasco, B.~Leimkuhler, D.~Phillips, and G.~Vilmart.
\newblock Efficient {L}angevin sampling with position-dependent diffusion.
\newblock {\em arXiv preprint}, 2501.02943, 2025.

\bibitem{chak_2023}
M.~Chak, N.~Kantas, T.~Leli{\`e}vre, and G.~Pavliotis.
\newblock Optimal friction matrix for underdamped langevin sampling.
\newblock {\em ESAIM: Mathematical Modelling and Numerical Analysis},
  57(6):3335--3371, 2023.

\bibitem{Chen2021}
M.~Chen.
\newblock Collective variable-based enhanced sampling and machine learning.
\newblock {\em Eur. Phys. J. B}, 94:211, 2021.

\bibitem{comer_2015}
J.~Comer, J.~C. Gumbart, J.~H{\'e}nin, T.~Leli{\`e}vre, A.~Pohorille, and
  C.~Chipot.
\newblock The {Adaptive Biasing Force} method: {E}verything you always wanted
  to know but were afraid to ask.
\newblock {\em The Journal of Physical Chemistry B}, 119(3):1129--1151, 2015.

\bibitem{cui_2024}
T.~Cui, X.~Tong, and O.~Zahm.
\newblock Optimal {R}iemannian metric for {P}oincar{\'e} inequalities and how
  to ideally precondition {L}angevin dymanics.
\newblock {\em arXiv preprint}, 2404.02554, 2024.

\bibitem{darve_2001}
E.~Darve and A.~Pohorille.
\newblock {Calculating free energies using average force}.
\newblock {\em The Journal of Chemical Physics}, 115(20):9169--9183, 2001.

\bibitem{das_2022}
A.~Das, B.~Kuznets-Speck, and D.~T. Limmer.
\newblock Direct evaluation of rare events in active matter from variational
  path sampling.
\newblock {\em Phys. Rev. Lett.}, 128:028005, 2022.

\bibitem{dellago_1999}
C.~Dellago, P.~G. Bolhuis, and D.~Chandler.
\newblock On the calculation of reaction rate constants in the transition path
  ensemble.
\newblock {\em The Journal of Chemical Physics}, 110(14):6617--6625, 1999.

\bibitem{dolbeault_2012}
J.~Dolbeault, A.~Klar, C.~Mouhot, and C.~Schmeiser.
\newblock Exponential rate of convergence to equilibrium for a model describing
  fiber lay-down processes.
\newblock {\em Applied Mathematics Research eXpress}, 2013(2):165--175, 2012.

\bibitem{duane_1987}
S.~Duane, A.~D. Kennedy, B.~J. Pendleton, and D.~Roweth.
\newblock Hybrid {{Monte Carlo}}.
\newblock {\em Physics Letters B}, 195(2):216--222, 1987.

\bibitem{Ferguson2018}
A.~L. Ferguson.
\newblock Machine learning and data science in soft materials engineering.
\newblock {\em Journal of Physics: Condensed Matter}, 30:043002, 2018.

\bibitem{girolami_2011}
M.~Girolami and B.~Calderhead.
\newblock Riemann {M}anifold {{Langevin}} and {{Hamiltonian Monte Carlo}}
  methods.
\newblock {\em Journal of the Royal Statistical Society: Series B (Statistical
  Methodology)}, 73(2):123--214, 2011.

\bibitem{Gkeka2020}
P.~Gkeka, G.~Stoltz, A.~Barati~Farimani, Z.~Belkacemi, M.~Ceriotti, J.~D.
  Chodera, A.~R. Dinner, A.~L. Ferguson, J.~B. Maillet, H.~Minoux, C.~Peter,
  F.~Pietrucci, A.~Silveira, A.~Tkatchenko, Z.~Trstanova, R.~Wiewiora, and
  T.~Leli\`evre.
\newblock {{M}achine {L}earning {F}orce {F}ields and {C}oarse-{G}rained
  {V}ariables in {M}olecular {D}ynamics: {A}pplication to {M}aterials and
  {B}iological {S}ystems}.
\newblock {\em J. Chem. Theory Comput.}, 16(8):4757--4775, Aug 2020.

\bibitem{glielmoreview}
A.~Glielmo, B.~Husic, A.~Rodriguez, C.~Clementi, F.~No{\'e}, and A.~Laio.
\newblock Unsupervised learning methods for molecular simulation data.
\newblock {\em Chemical Reviews}, 121(16):9722–9758, 2021.

\bibitem{graham_2022}
M.~M. Graham, A.~H. Thiery, and A.~Beskos.
\newblock Manifold {M}arkov chain {Monte Carlo} methods for {B}ayesian
  inference in diffusion models.
\newblock {\em Journal of the Royal Statistical Society Series B: Statistical
  Methodology}, 84(4):1229--1256, 2022.

\bibitem{grothaus_2016}
M.~Grothaus and P.~Stilgenbauer.
\newblock Hilbert space hypocoercivity for the {L}angevin dynamics revisited.
\newblock {\em Methods Funct. Anal. Topology}, 22(2):152--168, 2016.

\bibitem{haario_1999}
H.~Haario, E.~Saksman, and J.~Tamminen.
\newblock Adaptive proposal distribution for {R}andom {W}alk {Metropolis}
  algorithm.
\newblock {\em Computational Statistics}, 14(3):375--395, 1999.

\bibitem{haario_2001}
H.~Haario, E.~Saksman, and J.~Tamminen.
\newblock An adaptive {M}etropolis algorithm.
\newblock {\em Bernoulli}, 7(2):223--242, 2001.

\bibitem{hairer_2006}
E.~Hairer, C.~Lubich, and G.~Wanner.
\newblock {\em {Geometric Numerical Integration}}, volume~31 of {\em Springer
  Series in Computational Mathematics}.
\newblock Springer-Verlag, Berlin, second edition, 2006.

\bibitem{hastings_1970}
W.~K. Hastings.
\newblock Monte {{Carlo}} sampling methods using {{Markov}} chains and their
  applications.
\newblock {\em Biometrika}, 57(1):97--109, 1970.

\bibitem{horowitz_1991}
A.~M. Horowitz.
\newblock A generalized guided {Monte} {Carlo} algorithm.
\newblock {\em Physics Letters B}, 268(2):247--252, 1991.

\bibitem{hottovy_2015}
S.~Hottovy, A.~McDaniel, G.~Volpe, and J.~Wehr.
\newblock The {Smoluchowski-Kramers} limit of stochastic differential equations
  with arbitrary state-dependent friction.
\newblock {\em Communications in Mathematical Physics}, 336(3):1259--1283,
  2015.

\bibitem{henin_2004}
J.~Hénin and C.~Chipot.
\newblock Overcoming free energy barriers using unconstrained molecular
  dynamics simulations.
\newblock {\em The Journal of Chemical Physics}, 121(7):2904--2914, 2004.

\bibitem{kirkwood_1935}
J.~G. Kirkwood.
\newblock Statistical mechanics of fluid mixtures.
\newblock {\em The Journal of Chemical Physics}, 3(5):300--313, 1935.

\bibitem{legoll_2010}
F.~Legoll and T.~Lelièvre.
\newblock Effective dynamics using conditional expectations.
\newblock {\em Nonlinearity}, 23(9):2131, 2010.

\bibitem{lelievre_2011}
T.~Leli\`evre and K.~Minoukadeh.
\newblock Long-time convergence of an adaptive biasing force method: {T}he
  bi-channel case.
\newblock {\em Arch. Ration. Mech. Anal.}, 202(1):1--34, 2011.

\bibitem{lelievre_2013}
T.~{Leli{\`e}vre}, F.~{Nier}, and G.~A. {Pavliotis}.
\newblock Optimal non-reversible linear drift for the convergence to
  equilibrium of a diffusion.
\newblock {\em Journal of Statistical Physics}, 152(2):237--274, 2013.

\bibitem{lelievre_2008}
T.~Leli\`evre, M.~Rousset, and G.~Stoltz.
\newblock Long-time convergence of an adaptive biasing force method.
\newblock {\em Nonlinearity}, 21(6):1155--1181, 2008.

\bibitem{lelievre_2010}
T.~Leli{\`e}vre, M.~Rousset, and G.~Stoltz.
\newblock {\em {Free Energy Computations: A Mathematical Perspective}}.
\newblock {Imperial College Press}, 2010.

\bibitem{lelievre_2023_i}
T.~Leli{\`e}vre, R.~Santet, and G.~Stoltz.
\newblock Unbiasing {Hamiltonian Monte Carlo} algorithms for a general
  {H}amiltonian function.
\newblock {\em Foundations of Computational Mathematics}, 2024.

\bibitem{lelievre_2016}
T.~Leli\`evre and G.~Stoltz.
\newblock Partial differential equations and stochastic methods in molecular
  dynamics.
\newblock {\em Acta Numer.}, 25:681--880, 2016.

\bibitem{lelievre_2022}
T.~Leli\`evre, G.~Stoltz, and W.~Zhang.
\newblock Multiple projection {M}arkov chain {M}onte {C}arlo algorithms on
  submanifolds.
\newblock {\em IMA J. Numer. Anal.}, 43(2):737--788, 2023.

\bibitem{lelievre_2024}
T.~Lelièvre, G.~A. Pavliotis, G.~Robin, R.~Santet, and G.~Stoltz.
\newblock Optimizing the diffusion coefficient of overdamped {L}angevin
  dynamics.
\newblock {\em Math. Comp.}, electronically published on May 22, 2025, DOI:
  https://doi.org/10.1090/mcom/4098 (to appear in print).

\bibitem{lelievre_2012}
T.~Lelièvre, M.~Rousset, and G.~Stoltz.
\newblock Langevin dynamics with constraints and computation of free energy
  differences.
\newblock {\em Math. Comp.}, 81:2071--2125, 2012.

\bibitem{lelievre_2019}
T.~Lelièvre, M.~Rousset, and G.~Stoltz.
\newblock Hybrid {{Monte Carlo}} methods for sampling probability measures on
  submanifolds.
\newblock {\em Numerische Mathematik}, 143(2):379--421, 2019.

\bibitem{metropolis_1953}
N.~Metropolis, A.~W. Rosenbluth, M.~N. Rosenbluth, A.~H. Teller, and E.~Teller.
\newblock Equation of state calculations by fast computing machines.
\newblock {\em The Journal of Chemical Physics}, 21(6):1087--1092, 1953.

\bibitem{nilmeier_2011}
J.~P. Nilmeier, G.~E. Crooks, D.~D.~L. Minh, and J.~D. Chodera.
\newblock Nonequilibrium candidate {Monte Carlo} is an efficient tool for
  equilibrium simulation.
\newblock {\em Proceedings of the National Academy of Sciences},
  108(45):E1009--E1018, 2011.

\bibitem{noble_2023}
M.~Noble, V.~D. Bortoli, and A.~Durmus.
\newblock Unbiased constrained sampling with self-concordant barrier
  {Hamiltonian Monte Carlo}.
\newblock In {\em Thirty-seventh Conference on Neural Information Processing
  Systems}, 2023.

\bibitem{roberts_2007}
G.~O. Roberts and J.~S. Rosenthal.
\newblock Coupling and ergodicity of adaptive {M}arkov chain {M}onte {C}arlo
  algorithms.
\newblock {\em J. Appl. Probab.}, 44(2):458--475, 2007.

\bibitem{roberts_2009}
G.~O. Roberts and J.~S. Rosenthal.
\newblock Examples of adaptive {MCMC}.
\newblock {\em J. Comput. Graph. Statist.}, 18(2):349--367, 2009.

\bibitem{roberts_2002}
G.~O. Roberts and O.~Stramer.
\newblock Langevin diffusions and {M}etropolis-{H}astings algorithms.
\newblock {\em Methodol. Comput. Appl. Probab.}, 4:337--357 (2003), 2002.

\bibitem{roberts_1996}
G.~O. Roberts and R.~L. Tweedie.
\newblock Exponential convergence of {L}angevin distributions and their
  discrete approximations.
\newblock {\em Bernoulli}, 2(4):341--363, 1996.

\bibitem{rossky_1978}
P.~J. Rossky, J.~D. Doll, and H.~L. Friedman.
\newblock Brownian dynamics as smart {M}onte {C}arlo simulation.
\newblock {\em J. Chem. Phys.}, 69(10):4628--4633, 1978.

\bibitem{roussel_2018}
J.~Roussel and G.~Stoltz.
\newblock Spectral methods for {L}angevin dynamics and associated error
  estimates.
\newblock {\em ESAIM Math. Model. Numer. Anal.}, 52(3):1051--1083, 2018.

\bibitem{sidkyreview}
H.~Sidky, W.~Chen, and A.~Ferguson.
\newblock Machine learning for collective variable discovery and enhanced
  sampling in biomolecular simulation.
\newblock {\em Molecular Physics}, 118(5):e1737742, 2020.

\bibitem{straub_1988}
J.~E. Straub, M.~Borkovec, and B.~J. Berne.
\newblock Molecular dynamics study of an isomerizing diatomic in a
  {Lennard‐Jones} fluid.
\newblock {\em The Journal of Chemical Physics}, 89(8):4833--4847, 1988.

\bibitem{lelievre_2007}
G.~S. T.~Lelièvre, M.~Rousset.
\newblock Computation of free energy differences through nonequilibrium
  stochastic dynamics: {T}he reaction coordinate case.
\newblock {\em Journal of Computational Physics}, 222(2):624--643, 2007.

\bibitem{wang_2024}
M.~Wang, D.~Su, and W.~Wang.
\newblock Averaging on macroscopic scales with application to
  {S}moluchowski--{K}ramers approximation.
\newblock {\em Journal of Statistical Physics}, 191(2):22, 2024.

\bibitem{zappa_2018}
E.~Zappa, M.~{Holmes-Cerfon}, and J.~Goodman.
\newblock Monte {{Carlo}} on {{manifolds}}: {{Sampling densities}} and
  {{integrating functions}}.
\newblock {\em Communications on Pure and Applied Mathematics},
  71(12):2609--2647, 2018.

\bibitem{zhang_2016}
W.~Zhang, C.~Hartmann, and C.~Schütte.
\newblock Effective dynamics along given reaction coordinates{,} and reaction
  rate theory.
\newblock {\em Faraday Discuss.}, 195:365--394, 2016.

\end{thebibliography}

\begin{appendices}
  
  \section{Proof of Proposition~\ref{prop:standard_effective_dynamics}}
  \label{app:proposition_1}
  We show that
  \begin{equation*}
    b(z)=-\sigma^{2}(z)F'(z)+\beta^{-1}(\sigma^{2})'(z).
  \end{equation*}
  To this end, we use the following lemma (see for instance~\cite[Lemma 3.10]{lelievre_2010} or~\cite[Lemma 2.2]{legoll_2010}).
  \begin{lemma}
    \label{lem:derivative}
    For any smooth function~$\chi:\calQ\to\bbR$, consider
    \begin{equation*}
      \chi^{\xi}(z)=\int_{\Sigma(z)}\chi\left\lVert\nabla\xi\right\rVert^{-1}\rmd\sigma_{\Sigma(z)}.
    \end{equation*}
    The derivative of~$\chi^{\xi}$ reads:
    \begin{equation*}
      \frac{\rmd\chi^{\xi}}{\rmd z}(z)=\int_{\Sigma(z)}\left(\frac{\nabla\xi\cdot\nabla\chi}{\left\lVert\nabla\xi\right\rVert^{2}}+\chi\div\left(\frac{\nabla\xi}{\left\lVert\nabla\xi\right\rVert^{2}}\right)\right)\left\lVert\nabla\xi\right\rVert^{-1}\rmd\sigma_{\Sigma(z)}.
    \end{equation*}
  \end{lemma}
  In view of~\eqref{eq:conditional_measure} and~\eqref{eq:free_energy}, it holds
  \begin{equation*}
    \rmd\pi^{\xi}=\frac{Z^{-1}\rme^{-\beta V}\left\lVert\nabla\xi\right\rVert^{-1}\rmd\sigma_{\Sigma(z)}}{\rme^{-\beta F}},
  \end{equation*}
  so that, for instance,
  \begin{equation*}
    \sigma^{2}(z)=\frac{\displaystyle Z^{-1}\int_{\Sigma(z)}\left\lVert\nabla\xi\right\rVert^{2}\rme^{-\beta V}\left\lVert\nabla\xi\right\rVert^{-1}\rmd\sigma_{\Sigma(z)}}{\rme^{-\beta F(z)}}.
  \end{equation*}
  It follows from Lemma~\ref{lem:derivative} with~$\chi=\left\lVert\nabla\xi\right\rVert^{2}\rme^{-\beta V}$, that
  \begin{align*}
    (\sigma^{2})'(z)
    &=\frac{\displaystyle Z^{-1}\int_{\Sigma(z)}\left(\frac{\nabla\xi\cdot\nabla\left(\left\lVert\nabla\xi\right\rVert^{2}\rme^{-\beta V}\right)}{\left\lVert\nabla\xi\right\rVert^{2}}+\left\lVert\nabla\xi\right\rVert^{2}\rme^{-\beta V}\div\left(\frac{\nabla\xi}{\left\lVert\nabla\xi\right\rVert^{2}}\right)\right)\left\lVert\nabla\xi\right\rVert^{-1}\rmd\sigma_{\Sigma(z)}}{\rme^{-\beta F(z)}}\\
    &\qquad+\beta F'(z)\frac{\displaystyle Z^{-1}\int_{\Sigma(z)}\left\lVert\nabla\xi\right\rVert^{2}\rme^{-\beta V}\left\lVert\nabla\xi\right\rVert^{-1}\rmd\sigma_{\Sigma(z)}}{\rme^{-\beta F(z)}}\\
    &=
    \frac{
      \displaystyle
      \int_{\Sigma(z)}\left(
        \frac{\nabla\xi\cdot\nabla\left(\left\lVert\nabla\xi\right\rVert^{2}\right)}{\left\lVert\nabla\xi\right\rVert^{2}}-\beta\nabla V\cdot\nabla\xi
      \right)Z^{-1}\rme^{-\beta V}\left\lVert\nabla\xi\right\rVert^{-1}\rmd\sigma_{\Sigma(z)}
    }{
      \rme^{-\beta F(z)}
    }\\
    &\qquad+\frac{
      \displaystyle
      \int_{\Sigma(z)}\left(
        \Delta\xi+\left\lVert\nabla\xi\right\rVert^{2}\frac{\left(-\nabla\xi\cdot\nabla\left(\left\lVert\nabla\xi\right\rVert^{2}\right)\right)}{\left\lVert\nabla\xi\right\rVert^{4}}
      \right)Z^{-1}\rme^{-\beta V}\left\lVert\nabla\xi\right\rVert^{-1}\rmd\sigma_{\Sigma(z)}
    }{
      \rme^{-\beta F(z)}
    }\\
    &\qquad+\beta F'(z)\sigma^{2}(z)\\
    &=\beta\frac{
      \displaystyle\int_{\Sigma(z)}\left(-\nabla V\cdot\nabla\xi+\beta^{-1}\Delta\xi\right)Z^{-1}\rme^{-\beta V}\left\lVert\nabla\xi\right\rVert^{-1}\rmd\sigma_{\Sigma(z)}
    }{
      \rme^{-\beta F(z)}
    }+\beta F'(z)\sigma^{2}(z)\\
    &=
    \beta b(z)+\beta F'(z)\sigma^{2}(z).
  \end{align*}
  This concludes the proof of the identity~\eqref{eq:standard_drift_noise_relation_effective_dynamics}.
  
  \section{Proof of Proposition~\ref{prop:effective_dynamics}}
  \label{app:proof_effective_dynamics}
  
  The derivation of the effective dynamics for general diffusion processes (and for multivalued collective variables) is done in~\cite[Section~3.2]{zhang_2016}. Let us make precise the computations in the specific case of the dynamics~\eqref{eq:overdamped_langevin_diffusion} with diffusion~$D(q)=P^{\perp}(q)+a(\xi(q))P(q)$, where~$a:\xi(\calQ)\to\bbR_{+}^{*}$ and~$P,P^{\perp}$ are defined in~\eqref{eq:P}-\eqref{eq:P_perp}. Using Itô's lemma, it holds
  \begin{align}
    \label{eq:ito_formula_effective_dynamics}
    \rmd\xi(q_t)
    &=\left[
      \nabla\xi(q_t)\cdot\left(
        -D(q_t)\nabla V(q_t)+\beta^{-1}\div D(q_t)
      \right)+\beta^{-1}D(q_t):\nabla^{2}\xi(q_t)
    \right]\rmd t\\
    &\nonumber \quad  + \sqrt{2\beta^{-1}}\nabla\xi(q_t)\cdot D(q_t)^{1/2}\rmd W_t.
  \end{align}
  In view of~\eqref{eq:div_P}, it holds
  \begin{equation}
    \label{eq:divP_scalar_nabla_xi}
    \nabla\xi(q)\cdot\div P(q)=-\nabla\xi(q)\cdot\div P^{\perp}(q)=\Delta\xi(q)-\frac{\nabla\xi(q)^{\sfT}\nabla^{2}\xi(q)\nabla\xi(q)}{\left\lVert\nabla \xi(q)\right\rVert^{2}}=P^{\perp}(q):\nabla^{2}\xi(q),
  \end{equation}
  so that, using~\eqref{eq:chain_rule_div},    \begin{align*}
    \nabla\xi(q)\cdot\div D(q)+D(q):\nabla^{2}\xi(q)
    &=
    \nabla\xi(q)\cdot\left(\div P^{\perp}(q)+a'(\xi(q))\nabla\xi(q)+a(\xi(q))\div P(q)\right)\\
    &\quad
    +P^{\perp}(q):\nabla^{2}\xi(q)+a(\xi(q))P(q):\nabla^{2}\xi(q)\\
    &=
    a'(\xi(q))\left\lVert\nabla\xi(q)\right\rVert^{2}+a(\xi(q))\left[\nabla\xi(q)\cdot\div P(q)+P(q):\nabla^{2}\xi(q)\right]\\
    &=
    a'(\xi(q))\left\lVert\nabla\xi(q)\right\rVert^{2}+a(\xi(q))\left[P^{\perp}(q):\nabla^{2}\xi(q)+P(q):\nabla^{2}\xi(q)\right]\\
    &=
    a'(\xi(q))\left\lVert\nabla\xi(q)\right\rVert^{2}+a(\xi(q))\Delta\xi(q).
  \end{align*}
  Using that~$D^{1/2}=P^{\perp}+\sqrt{a\circ\xi}\,P$ (similarly to~\eqref{eq:sqrt_inv_det_diffusion}) and~\eqref{eq:P_nabla_xi_relations}, the dynamics~\eqref{eq:ito_formula_effective_dynamics} therefore simplifies as
  \begin{align*}
    \rmd\xi(q_t)
    &=\left(
      -a(\xi(q_t))\left(
        \nabla V(q_t)\cdot\nabla\xi(q_t)-\beta^{-1}\Delta\xi(q_t)
      \right)+\beta^{-1}a'(\xi(q_t))\left\lVert\nabla\xi(q_t)\right\rVert^{2}
    \right)\rmd t\\
    &\nonumber\qquad+\sqrt{2\beta^{-1}a(\xi(q_t))}\nabla\xi(q_t)\cdot\rmd W_t,
  \end{align*}
  which is exactly~\eqref{eq:overdamped_langevin_dynamics_effective_dynamics}. Using that~$a$ and~$a'$ are functions of~$\xi(q)=z$, the conditional expectations of the drift and noise are thus given by
  \begin{equation*}
    \left\lbrace
    \begin{aligned}
      b_a(z)&=-a(z)\int_{\Sigma(z)}
        \left(\nabla V\cdot\nabla\xi- \beta^{-1}\Delta\xi\right)\rmd\pi^{\xi}
      +\beta^{-1}a'(z)\int_{\Sigma(z)}\left\lVert\nabla\xi\right\rVert^{2}\rmd\pi^{\xi},\\
      \sigma^{2}_a(z)&=a(z)\int_{\Sigma(z)}\left\lVert\nabla\xi\right\rVert^{2}\rmd\pi^{\xi},
    \end{aligned}
    \right.
  \end{equation*}
  which is exactly~\eqref{eq:drift_noise_effective_dynamics}. Lastly,~\eqref{eq:drift_noise_relation_effective_dynamics} follows directly from~\eqref{eq:standard_drift_noise_relation_effective_dynamics}. This concludes the proof.

  \section{Normalization constant of the diffusion}
  \label{app:kappa_alpha}

  Let us first recall that the computation of the normalization constant of the diffusion is essentially needed to obtain similar optimal time steps (when $\alpha$ is modified) for the presentation of our results. In practice, this normalization constant simply rescales the time step~$\Delta t$. In addition, when performing simulations where the free energy is learned on the fly, computing this normalization constant is useful in order to prevent unstable dynamics due to large values of the diffusion.

  We use the~$L^{p}$ normalization constraint~\eqref{eq:constraint_lp_norm} with~$p=1$ to compute the scalar~$\kappa_\alpha$ in~\eqref{eq:optimal_diffusion_class}. The main motivation to use the~$L^{p}$ constraint is that, for the class of diffusions~\eqref{eq:optimal_diffusion_class}, it can be rewritten as an integral over the latent space, so that its computation is numerically tractable. The choice~$p=1$ is motivated by the fact that the integral corresponds in this case to the average mean square displacement over a time step for the Brownian part of the dynamics. Moreover, it does not require to estimate an additional conditional expectation, see~\eqref{eq:conditioning_formula_normalizing_constant} below.

  Let us first recall the conditioning formula involving the two measures~\eqref{eq:conditional_measure}-\eqref{eq:image_measure}: for any test functions~$f:\xi(\calQ)\to\bbR$ and~$g:\calQ\to\bbR$,
  \begin{equation}
    \label{eq:conditioning_formula}
    \int_{\calQ}f(\xi(q))g(q)Z^{-1}\rme^{-\beta V(q)}\rmd q=\int_{\xi(\calQ)}f(z)\left(\int_{\Sigma(z)}g\,\rmd\pi^{\xi}(\cdot|z)\right)\rme^{-\beta F(z)}\rmd z.
  \end{equation}
  Since the squared Frobenius norm of the diffusion~\eqref{eq:optimal_diffusion_class}, given by
  \begin{equation*}
    \left\lvert D_{\alpha}(q)\right\rvert_{\rmF}^{2}=\kappa_{\alpha}^{2}\left((d-1)+a_{\alpha}(\xi(q))^{2}\right),
  \end{equation*}
  is a function of~$\xi(q)$, we can use the conditioning formula~\eqref{eq:conditioning_formula} to rewrite the~$L^{p}$ constraint~\eqref{eq:constraint_lp_norm} as
  \begin{align}
    \int_{\calQ}\left\lvert D(q)\right\rvert_{\rmF}^{p}\rme^{-\beta pV(q)}\rmd q
    &\nonumber
    =Z\int_{\calQ}\left\lvert D(q)\right\rvert_{\rmF}^{p}\rme^{-\beta(p-1)V(q)}Z^{-1}\rme^{-\beta V(q)}\rmd q\\
    &\label{eq:conditioning_formula_normalizing_constant}
    =Z\kappa_{\alpha}^{p}\int_{\xi(\calQ)}\left((d-1)+a_{\alpha}(z)^{2}\right)^{p/2}\left(
      \int_{\Sigma(z)}\rme^{-\beta (p-1)V}\rmd\pi^{\xi}(\cdot|z)
    \right)\rme^{-\beta F(z)}\rmd z.
  \end{align}
  Choosing~$p=1$, the integral with respect to the conditional measure on the right-hand side of~\eqref{eq:conditioning_formula_normalizing_constant} is simply equal to~1. The constant~$Z$, which is independent of the diffusion~$D_{\alpha}$, can be omitted (it simply rescales the time step). This does not alter the optimality of the homogenized diffusion~\eqref{eq:optimal_diffusion_homogenization}. The only effect is potentially having a suboptimal multiplicative constant in factor of~$\rme^{\alpha\beta F}$ when~$\alpha=1$ for the effective dynamics. In view of the constraint~\eqref{eq:constraint_lp_norm}, the scalar~$\kappa_{\alpha}$ is then given by~\eqref{eq:kappa_alpha}.

  \begin{remark}
    It is possible to use an~$L^{p}$ constraint with~$p>1$ by estimating the conditional measure on the right-hand side of~\eqref{eq:conditioning_formula_normalizing_constant}, typically using ABF-like methods as described in Section~\ref{sec:adaptive_learning}.
  \end{remark}
  
  \section{Computations related to the numerical experiment}
  \label{app:computations}

In this section, we provide the analytical formulas related to the collective variable and diffusion used in our numerical example defined in Section~\ref{sec:MALA_numerics}.
  
  \paragraph{Derivatives of the collective variable~$\xi$.}
  
  The second derivatives of $\xi$ write
  \begin{equation}
    \label{eq:xi_second_derivatives}
    \left\lbrace
    \begin{aligned}
      \partial_{x_1}^{2}\xi(q)     & = \frac{1}{2w\left\lVert q_2-q_1\right\rVert}\left[1-\frac{(x_1-x_2)^{2}}{\left\lVert q_2-q_1\right\rVert^{2}}\right],\\
      \partial_{y_1}^{2}\xi(q)     & = \frac{1}{2w\left\lVert q_2-q_1\right\rVert}\left[1-\frac{(y_1-y_2)^{2}}{\left\lVert q_2-q_1\right\rVert^{2}}\right],\\
      \partial_{x_1,y_1}^{2}\xi(q) & = -\frac{1}{2w\left\lVert q_2-q_1\right\rVert^{3}}(x_1-x_2)(y_1-y_2),\\
      \partial_{x_1,x_2}^{2}\xi(q) & =-\partial_{x_2}^{2}\xi(q)=-\partial_{x_1}^{2}\xi(q),\\
      \partial_{y_1,y_2}^{2}\xi(q) & =-\partial_{y_2}^{2}\xi(q)=-\partial_{y_1}^{2}\xi(q),\\
      \partial_{y_2,x_1}^{2}\xi(q) & = \partial_{y_1,x_2}^{2}\xi(q)=-\partial_{x_2,y_2}^{2}\xi(q)=-\partial_{x_1,y_1}^{2}\xi(q).
    \end{aligned}
    \right.
  \end{equation}
  
  \paragraph{Computing the divergence of the diffusion matrix.} Since~$\left\lVert\nabla\xi\right\rVert^2=1/(2w^{2})$ is constant, we can use the simplified expressions for the divergence of Remark~\ref{rem:MALA_gradient_xi_constant}. Using the various relations between the first and second derivatives (see~\eqref{eq:xi_first_derivatives} and~\eqref{eq:xi_second_derivatives}), the expression of the divergence of the diffusion matrix~\eqref{eq:divergence_diffusion} reduces to
  \begin{equation}
\label{eq:divergence_diffusion_numerical_experiment}
    \div D_{\alpha}(q) = 2w^{2}\kappa_{\alpha}\left[
      2\left(2w^{2}\rme^{\alpha\beta F\circ\xi(q)}-1\right)\left(\partial_{x_1}^2\xi(q)+\partial_{y_1}^2\xi(q)\right)+\alpha\beta F'\circ\xi(q)\rme^{\alpha\beta F\circ\xi(q)}
    \right]\nabla\xi(q).
  \end{equation}
  In particular, it holds~$\left[\div D_{\alpha}(q)\right]_3=-\left[\div D_{\alpha}(q)\right]_1$ and~$\left[\div D_{\alpha}(q)\right]_4=-\left[\div D_{\alpha}(q)\right]_2$.
  
  \section{Thermodynamic integration}
  \label{app:TI}
  
  In this section, we briefly review thermodynamic integration, and describe how this was implemented to run the numerical experiment in Section~\ref{sec:MALA_numerics}. This section is based on~\cite[Chapter~3]{lelievre_2010}.
  
  Thermodynamic integration is based on the following identity:
  \begin{equation}
    \label{eq:TI_identity}
    \forall z\in[z_{\rm min},z_{\rm max}],\qquad F(z)-F(z_{\rm min})=\int_{z_{\rm min}}^{z}F'(y)\rmd y,
  \end{equation}
  where the mean force~$F'$ is defined in~\eqref{eq:mean_force}, and the values of~$z_{\rm min},z_{\rm max}$ are set by the practitioner. Ideally, these values define an interval such that the typical values of the collective variable lie in it. The integral on the right-hand side of~\eqref{eq:TI_identity} is for instance approximated by a midpoint rule, using a grid of~$N_z$ evenly spaced point in the interval~$[z_{\rm min},z_{\rm max}]$:
  \begin{equation}
    \label{eq:TI_left_riemann_rule}
    \int_{z_{\rm min}}^{z}F'(y)\rmd y\approx \sum_{i=1}^{B(z)}F'(z^{i})\Delta z,\qquad z^{i}=z_{\rm min}+\left(i-\frac{1}{2}\right)\Delta z,\qquad \Delta z=\frac{z_{\rm max}-z_{\rm min}}{N_z},
  \end{equation}
  where~$B:[z_{\rm min},z_{\rm max})\to\left\lbrace 1,\dots,N_z\right\rbrace$ outputs the bin number as a function of the collective variable value:
  \begin{equation*}
    B(z)=1+\left\lfloor\frac{z-z_{\rm min}}{\Delta z}\right\rfloor.
  \end{equation*}
  Let~$1\leqslant i\leqslant N_{z}$. To estimate one value~$F'(z^{i})$, we need to compute an expectation with respect to the conditional measure~$\pi^{\xi}(\rmd q|z^{i})$ defined in~\eqref{eq:conditional_measure}. This measure rewrites
  \begin{equation*}
    \pi^{\xi}(\rmd q|z^{i})=\frac{\rme^{-\beta V^{\xi}(q)}\sigma_{\Sigma(z^{i})}(\rmd q)}{\displaystyle\int_{\Sigma(z^{i})}\rme^{-\beta V^{\xi}}\rmd\sigma_{\Sigma(z^{i})}},
  \end{equation*}
  where the modified potential~$V^{\xi}$ writes~$V^{\xi}=V+\beta^{-1}\ln\left\lVert\nabla\xi\right\rVert$. We therefore need to sample the probability measure whose density is proportional to~$\rme^{-\beta V^{\xi}}$ on~$\Sigma(z^{i})$. This can be done using trajectories of the projected dynamics
  \begin{equation*}
    \left\lbrace
    \begin{aligned}
      &\rmd q_t=-\nabla V^{\xi}(q_t)\rmd t+\sqrt{2\beta^{-1}}\rmd W_t+\nabla\xi(q_t)\rmd\lambda_t,\\
      &\text{with }(\lambda_t)_{t\geqslant0}\text{ an adapted process such that }\xi(q_t)=z^{i}.
    \end{aligned}
    \right.
  \end{equation*}
  Here, the process~$(\lambda_t)_{t\geqslant0}$ acts as a Lagrange multiplier for the constraints $(\xi(q_t)=z^{i})_{t \geqslant 0}$. This dynamics is numerically integrated with the following predictor-correct scheme: let~$q^{0}\in\Sigma(z^{i})$ and for~$n\geqslant 0$,
  \begin{equation}
    \label{eq:TI_update}
    \left\lbrace
    \begin{aligned}
      q^{n+1}&=q^{n}-\nabla V^{\xi}(q^{n})\Delta t+\sqrt{2\beta^{-1}\Delta t}\,\rmG^{n+1}+\nabla\xi(q^{n+1})\Delta\lambda^{n+1},\\
      \xi(q^{n+1})&=z^{i},
    \end{aligned}
    \right.
  \end{equation}
  where~$\rmG^{n+1}\sim\calN(0,\rmI_d)$. The value of the Lagrange multiplier~$\Delta\lambda^{n+1}$ is determined by the requirement that the constraint~$\xi(q^{n+1})=z^{i}$ is satisfied.
  
  \paragraph{Computation of the Lagrange multiplier for the numerical example of Section~\ref{sec:MALA_numerics}.} Since~$\left\lVert\nabla\xi\right\rVert$ is constant for our numerical experiment, we simply set~$V^{\xi}=V$. In practice, the update~\eqref{eq:TI_update} is performed by first computing the unconstrained move
  \begin{equation*}
    \widetilde{q}^{n+1}=q^{n}-\nabla V(q^{n})\Delta t+\sqrt{2\beta ^{-1}\Delta t}\,\rmG^{n+1},
  \end{equation*}
  and then solving for~$\Delta\lambda^{n+1}$. In view of~\eqref{eq:xi_first_derivatives}, only the components of~$q^{n+1}$ corresponding to the dimer differ from the ones of~$\widetilde{q}^{n+1}$. Denote by~$q_{1},q_{2}\in(\ell\bbT)^{2}$ the positions of the particles composing the dimer. Then,~\eqref{eq:TI_update} rewrites
  \begin{equation}
    \label{eq:TI_dimer_update}
    \left\lbrace
    \begin{aligned}
      q_1^{n+1}&=\widetilde{q}_1^{n+1}+\frac{q_1^{n+1}-q_2^{n+1}}{2w\left\lVert q_1^{n+1}-q_2^{n+1}\right\rVert}\Delta\lambda^{n+1},\\
      q_2^{n+1}&=\widetilde{q}_2^{n+1}-\frac{q_1^{n+1}-q_2^{n+1}}{2w\left\lVert q_1^{n+1}-q_2^{n+1}\right\rVert}\Delta\lambda^{n+1}.
    \end{aligned}
    \right.
  \end{equation}
  Subtracting both equalities leads to
  \begin{equation}
    \label{eq:lagrange_multiplier}
    \left(1-\frac{\Delta\lambda^{n+1}}{w\left\lVert q_1^{n+1}-q_2^{n+1}\right\rVert}\right)(q_2^{n+1}-q_1^{n+1})=\widetilde{q}_2^{n+1}-\widetilde{q}_1^{n+1}.
  \end{equation}
  From~$\xi(q^{n+1})=z^{i}$ it follows that~$\left\lVert q_1^{n+1}-q_2^{n+1}\right\rVert=2wz^{i}+r_0$. Taking the norm in~\eqref{eq:lagrange_multiplier} yields
  \begin{equation*}
    \left\lvert w(2wz^{i}+r_0)-\Delta\lambda^{n+1} \right\rvert=w\left\lVert\widetilde{q}_2^{n+1}-\widetilde{q}_1^{n+1}\right\rVert=w(2w\xi(\widetilde{q}^{n+1})+r_0).
  \end{equation*}
  It follows that
  \begin{equation*}
    \Delta\lambda^{n+1}=w(2wz^{i}+r_0)\pm w(2w\xi(\widetilde{q}^{n+1})+r_0).
  \end{equation*}
  We choose the solution with the minus sign as this choice leads to a consistent discretization in the limit~$\Delta t\to0$. Therefore, the Lagrange multiplier is equal to
  \begin{equation*}
    \Delta\lambda^{n+1}=2w^{2}\left(z^{i}-\xi(\widetilde{q}^{n+1})\right).
  \end{equation*}
  Notice that, by adding the two equalities in~\eqref{eq:TI_dimer_update}, it holds~$q_1^{n+1}+q_2^{n+1}=\widetilde{q}_1^{n+1}+\widetilde{q}_2^{n+1}$. Therefore, writing~$q_1^{n+1}=(q_1^{n+1}+q_2^{n+1}-(q_2^{n+1}-q_1^{n+1}))/2$ (and similarly for~$q_2^{n+1}$) and using~\eqref{eq:lagrange_multiplier} leads to
  \begin{equation}
    \label{eq:TI_update_new_positions}
    \left\lbrace
    \begin{aligned}
      q_1^{n+1}&=\frac{1}{2}\left(
        1+\frac{1}{1-\frac{\Delta\lambda^{n+1}}{w(2wz^{i}+r_0)}}
      \right)\widetilde{q}^{n+1}_1+\frac{1}{2}\left(
        1-\frac{1}{1-\frac{\Delta\lambda^{n+1}}{w(2wz^{i}+r_0)}}
      \right)\widetilde{q}_2^{n+1},\\
      q_2^{n+1}&=\frac{1}{2}\left(
        1-\frac{1}{1-\frac{\Delta\lambda^{n+1}}{w(2wz^{i}+r_0)}}
      \right)\widetilde{q}^{n+1}_1+\frac{1}{2}\left(
        1+\frac{1}{1-\frac{\Delta\lambda^{n+1}}{w(2wz^{i}+r_0)}}
      \right)\widetilde{q}_2^{n+1}.
    \end{aligned}
    \right.
  \end{equation}

  \begin{remark}
    Note that the update~\eqref{eq:TI_update_new_positions} is undefined when~$\Delta\lambda^{n+1}=w(2wz^{i}+r_0)$. This corresponds to~$\left\lVert\widetilde{q}_2^{n+1}-\widetilde{q}_1^{n+1}\right\rVert=0$, which is not observed in practice due to the strong repulsive force associated with the potential $V_{\rm DW}$ when the distance goes to~0 (see~\eqref{eq:dimer_potential}).
  \end{remark}
  
  \paragraph{Complying with the periodic boundary conditions.} Note that the norm~$\left\lVert q_2-q_1\right\rVert$ is computed taking the periodic boundary conditions into account. In practice, we first translate~$\widetilde{q}_1^{n+1}$ to the frame of reference of~$\widetilde{q}_2^{n+1}$, apply the update, and periodize.

  \paragraph{Hyperparameters for the numerical experiment.} To obtain the mean force and free energy profiles to run the simulations in Section~\ref{sec:MALA_numerics}, we used~$N_z=100$ points in the interval~$[z_{\rm min},z_{\rm max}]$ (see~\eqref{eq:TI_left_riemann_rule}) and used a time step~$\Delta t=2.5\times10^{-5}$. The initial configuration for the first level~$z^{0}=z_{\rm min}$ is set to~$q^{0}$ defined in~\eqref{eq:q0} except for~$y^{0}_{2}$ which is set to~$y^{0}_{1}+r_0+2wz_{\rm min}$ so that~$\xi(q^{0})=z_{\rm min}$. At the end of each simulation at the level~$z^{i}$, the final configuration is changed so that it becomes the initial configuration for the simulation for the level~$z^{i+1}$ (\emph{i.e.}~we slightly move the second component of the second particle so that the new constraint is satisfied). The simulation for each level~$z^{i}$ was run for a physical time of~$125$. The free energy is then reconstructed by integrating the mean force over the interval~$[z_{\rm min},z_{\rm max}]$.

  \section{Integrating the Hamiltonian dynamics with the GSV numerical scheme using Newton's method}
  \label{app:RMHMC_implementation}

  Let us make precise how the Hamiltonian dynamics are integrated using the GSV scheme~\eqref{eq:GSV} and Newton's method. In particular, this will highlight how we can limit the computational overhead compared to the standard HMC algorithm when the collective variable is a function of a small number of components.

  The gradients of the Hamiltonian function~\eqref{eq:H_diffusion} are given by
  \begin{align}
    \nabla_q H(q,p)
    &\nonumber=\nabla V(q)-\frac{1}{2}\frac{a_{\alpha}'(\xi(q))}{a_{\alpha}(\xi(q))}\nabla\xi(q)+\frac{\kappa_{\alpha}}{2}\frac{\left(\nabla\xi(q)\cdot p\right)^{2}a_{\alpha}'(\xi(q))}{\left\lVert\nabla\xi(q)\right\rVert^{2}}\nabla \xi(q)\\
    &\nonumber\quad+\kappa_{\alpha} (a_{\alpha}(\xi(q))-1)\frac{\nabla\xi(q)\cdot p}{\left\lVert\nabla\xi(q)\right\rVert^{2}}\nabla^{2}\xi(q)p\\
    &\label{eq:nabla_q_H}\quad-\kappa_{\alpha}(a_{\alpha}(\xi(q))-1)\frac{\left(\nabla\xi(q)\cdot p\right)^{2}}{\left\lVert\nabla\xi(q)\right\rVert^{4}}\nabla^{2}\xi(q)\nabla\xi(q).
  \end{align}
  and
  \begin{align*}
    \nabla_p H(q,p)=\kappa_{\alpha} p+\kappa_{\alpha} (a_{\alpha}(\xi(q))-1)\frac{\nabla\xi(q)\cdot p}{\left\lVert\nabla\xi(q)\right\rVert^{2}}\nabla\xi(q).
  \end{align*}

  \paragraph{Solving the implicit problem on the momenta.} Recall that the GSV scheme is defined in~\eqref{eq:GSV}. To solve for~$p^{n+1/2}$ we introduce the map~$g(\cdot;q^{n},p^{n})$ defined by
  \begin{align}
    g(p;q^{n},p^{n})
    &\nonumber=p-p^{n}+\frac{\Delta t}{2}\nabla V(q^{n})-\frac{\Delta t}{4}\frac{a_{\alpha}'(\xi(q^{n}))}{a_{\alpha}(\xi(q^{n}))}\nabla \xi(q^{n})+\frac{\kappa_{\alpha}\Delta t}{4}\frac{\left(\nabla\xi(q^{n})\cdot p\right)^{2}a_{\alpha}'(\xi(q^{n}))}{\left\lVert\nabla\xi(q^{n})\right\rVert^{2}}\nabla \xi(q^{n})\\
    &\label{eq:g_p}\quad+\frac{\kappa_{\alpha}\Delta t}{2}(a_{\alpha}(\xi(q^{n}))-1)\frac{\nabla\xi(q^{n})\cdot p}{\left\lVert\nabla\xi(q^{n})\right\rVert^{2}}\nabla^{2}\xi(q^{n})p\\
    &\nonumber\quad-\frac{\kappa_{\alpha}\Delta t}{2}(a_{\alpha}(\xi(q^{n}))-1)\frac{\left(\nabla\xi(q^{n})\cdot p\right)^{2}}{\left\lVert\nabla\xi(q^{n})\right\rVert^{4}}\nabla^{2}\xi(q^{n})\nabla\xi(q^{n}).
  \end{align}
  We want to find~$p^{n+1/2}$ that solves~$g(p^{n+1/2},q^{n},p^{n})=0$. The Jacobian matrix of~$g$ at a point~$p$ is given by
  \begin{align}
    \nabla g(p;q^{n},p^{n})
    &=\nonumber\rmI_{d}+\frac{\kappa_{\alpha}\Delta t}{2}\frac{\left(\nabla\xi(q^{n})\cdot p\right)a_{\alpha}'(\xi(q^{n}))}{\left\lVert\nabla\xi(q^{n})\right\rVert^{2}}\nabla \xi(q^{n})\otimes\nabla\xi(q^{n})\\
    &\nonumber\quad+\frac{\kappa_{\alpha}\Delta t}{2}\frac{a_{\alpha}(\xi(q^{n}))-1}{\left\lVert\nabla\xi(q^{n})\right\rVert^{2}}\left[\nabla^{2}\xi(q^{n})p\right]\otimes\nabla\xi(q^{n})\\
    &\nonumber\quad+\frac{\kappa_{\alpha}\Delta t}{2}(a_{\alpha}(\xi(q^{n}))-1)\frac{\nabla\xi(q^{n})\cdot p}{\left\lVert\nabla\xi(q^{n})\right\rVert^{2}}\nabla^{2}\xi(q^{n})\\
    &\label{eq:nabla_g_p}\quad-\kappa_{\alpha}\Delta t (a_{\alpha}(\xi(q^{n}))-1)\frac{\nabla\xi(q^{n})\cdot p}{\left\lVert\nabla\xi(q^{n})\right\rVert^{4}}\left[\nabla^{2}\xi(q^{n})\nabla\xi(q^{n})\right]\otimes\nabla\xi(q^{n}).
  \end{align}
  We then define the Newton sequence as follows: let~$p^{n+1/2,0}=p^{n+1/2}$ and~$p^{n+1/2,1}=p^{n}-\frac{\Delta t}{2}\nabla_q H(q^{n},p^{n})$, and for any~$i\geqslant 1$, 
  \begin{equation*}
    p^{n+1/2,i+1}=p^{n+1/2,i}-\left(\nabla g(p^{n+1/2,i};q^{n},p^{n})\right)^{-1}g(p^{n+1/2,i};q^{n},p^{n}).
  \end{equation*}
  In practice, the linear system
  \begin{equation}
    \label{eq:newton_sequence}
    \nabla g(p^{n+1/2,i};q^{n},p^{n})u=-g(p^{n+1/2,i};q^{n},p^{n}),
  \end{equation}
  is solved for the unknown~$u=p^{n+1/2,i+1}-p^{n+1/2,i}$. The Newton sequence is constructed iteratively, until a maximum number of iterations~$N_{\rm Newton}$ has been attained or when stopping criteria of the following form are met: 
  \begin{equation}
    \label{eq:newton_thresholds}
    \left\lVert p^{n+1/2,i+1}-p^{n+1/2,i}\right\rVert<\eta_{\rm Newton,Cauchy},\qquad \left\lVert g(p^{n+1/2,i+1})\right\rVert<\eta_{\rm Newton,root},
  \end{equation}
  where~$\eta_{\rm Newton,Cauchy},\eta_{\rm Newton,root}>0$ are user-defined thresholds. When the maximum number of iterations is attained without obtaining convergence, or when the matrix~$\nabla g$ is not (numerically) invertible, the computation fails and~$S(q^{n},p^{n})$ is returned (see~\eqref{eq:GSV_rev}).

  \paragraph{Solving the implicit problem on the positions.} If the first implicit problem succeeded, one can proceed and try to compute the new position~$q^{n+1}$. Let us define to this end the map
  \begin{align}
    h(q;q^{n},p^{n+1/2})
    &\nonumber=q-q^{n}-\kappa_{\alpha}\Delta tp^{n+1/2}-\frac{\kappa_{\alpha}\Delta t}{2}(a_{\alpha}(\xi(q^{n}))-1)\frac{\nabla\xi(q^{n})\cdot p^{n+1/2}}{\left\lVert\nabla\xi(q^{n})\right\rVert^{2}}\nabla\xi(q^{n})\\
    &\label{eq:h_q}\quad-\frac{\kappa_{\alpha}\Delta t}{2}(a_{\alpha}(\xi(q))-1)\frac{\nabla\xi(q)\cdot p^{n+1/2}}{\left\lVert\nabla\xi(q)\right\rVert^{2}}\nabla\xi(q).
  \end{align}
  We want to find~$q^{n+1}$ that solves~$h(q^{n+1};q^{n},p^{n+1/2})=0$. The Jacobian of~$h$ at a point~$q$ is given by
  \begin{align}
    \nabla h(q;q^{n},p^{n+1/2})
    &\nonumber=\rmI_d-\frac{\kappa_{\alpha}\Delta t}{2}\frac{\left(\nabla\xi(q)\cdot p^{n+1/2}\right)a_{\alpha}'(\xi(q))}{\left\lVert\nabla\xi (q)\right\rVert^{2}}\nabla\xi(q)\otimes\nabla \xi(q)\\
    &\nonumber\quad-\frac{\kappa_{\alpha}\Delta t}{2}\frac{a_{\alpha}(\xi(q))-1}{\left\lVert\nabla\xi(q)\right\rVert^{2}}\nabla\xi(q)\otimes\left[\nabla^{2}\xi(q) p^{n+1/2}\right]\\
    &\label{eq:nabla_h_q}\quad+\kappa_{\alpha}\Delta t(a_{\alpha}(\xi(q))-1)\frac{\nabla\xi(q)\cdot p^{n+1/2}}{\left\lVert\nabla\xi(q)\right\rVert^{4}}\left[\nabla^{2}\xi(q)\nabla\xi(q)\right]\otimes\nabla\xi(q)\\
    &\nonumber\quad-\frac{\kappa_{\alpha}\Delta t}{2}(a_{\alpha}(\xi(q))-1)\frac{\nabla\xi(q)\cdot p^{n+1/2}}{\left\lVert\nabla\xi(q)\right\rVert^{2}}\nabla^{2}\xi(q).
  \end{align}
  We then follow a similar strategy as for the first implicit problem. If the computations fail, the flow again returns~$S(q^{n},p^{n})$. We emphasize that the forces~$-\nabla V$ need to be computed only once per integration of the Hamiltonian dynamics.

  \begin{remark}
    When~$\left\lVert\nabla\xi\right\rVert$ is constant, any term of the form~$\nabla^{2}\xi(q)\nabla\xi(q)$ appearing in~\eqref{eq:nabla_q_H},~\eqref{eq:g_p},~\eqref{eq:nabla_g_p},~\eqref{eq:nabla_h_q} vanishes.
  \end{remark}

  \paragraph{Numerical reversibility check.}  If all the Newton methods needed to define~$\varphi_{\Delta t}^{\rm GSV}\circ S\circ\varphi_{\Delta t}^{\rm GSV}(q,p)$ have converged, then the reversibility check~\eqref{eq:SV_property} is performed. In practice, it is implemented with a user defined threshold~$\eta_{\rm rev}$: the property~\eqref{eq:SV_property} \textit{numerically} holds when 
  \begin{equation}
    \label{eq:reversibiliy_check}
    \left\lVert S\circ\varphi_{\Delta t}^{\rm GSV}\circ S\circ\varphi_{\Delta t}^{\rm GSV}(q,p)-(q,p)\right\rVert<\eta_{\rm rev}.
  \end{equation}

  \paragraph{Optimization of the implementation of RMHMC.}
  In view of the various formulas needed to run Newton's method, the additional computations (compared to the standard HMC scheme) only require gradients or Hessians of the collective variable. When the collective variable is a function of only a few components of the position~$q$, one can therefore limit computational overheads and the associated storage capacity, as already discussed for MALA in Section~\ref{sec:MALA_numerics}. As a guiding example, consider the physical system used for our numerical experiment in Section~\ref{sec:MALA_numerics}. The gradient of~$\xi$ only contains 4 nonzero components (see~\eqref{eq:xi_first_derivatives}) and the Hessian of~$\xi$ only contains a~$4\times4$ nonzero submatrix, so that only the first four components in~\eqref{eq:g_p}-\eqref{eq:h_q} are non trivially modified (compared to the standard HMC scheme). Furthermore, the matrices~$\nabla g$ and~$\nabla h$ appearing in~\eqref{eq:nabla_g_p} and~\eqref{eq:nabla_h_q} are block diagonal with a structure similar to the diffusion, see~\eqref{eq:diffusion_block_diagonal}. This means that the update of the Newton sequence described in~\eqref{eq:newton_sequence} can be done by (i) solving only a~$4\times4$ linear system for the first 4 components and (ii) updating the other components as the update is explicit for these components. These observations can be generalized to any collective variable which is a function of only~$k\leqslant d$ components of the positions.

  \paragraph{Optimization of the implementation of RMGHMC.} Adding on what has already been stated for the RMHMC algorithm, the steps associated to an Ornstein--Uhlenbeck integration can also be optimized: for the diffusion~\eqref{eq:optimal_diffusion_class}, the update~\eqref{eq:OU_update_rmhmc} can be performed without solving a linear system, as the matrix
  \begin{equation*}
    \rmI_d+\frac{\Delta t}{4}\gamma D_{\alpha}(q^{n})=\left(1+\frac{\Delta t}{4}\gamma\kappa_\alpha\right)P^{\perp}(q^{n})+\left(1+\frac{\Delta t}{4}\gamma\kappa_{\alpha}a\circ\xi(q^{n})\right)P(q^{n}),
  \end{equation*}
  has the explicit inverse
  \begin{equation*}
    \left(1+\frac{\Delta t}{4}\gamma\kappa_\alpha\right)^{-1}P^{\perp}(q^{n})+\left(1+\frac{\Delta t}{4}\gamma\kappa_{\alpha}a\circ\xi(q^{n})\right)^{-1}P(q^{n}).
  \end{equation*}

  \paragraph{Hyperparameters used for the numerical experiment in Section~\ref{sec:HMC_numerics}.}
  The maximum number of iterations for Newton's method is set to~$N_{\rm Newton}=100$, and we use~$\eta_{\rm Newton, Cauchy}=\eta_{\rm Newton,root}=10^{-12}$ for the tolerance checks in Newton's method, see~\eqref{eq:newton_thresholds}. The tolerance for the reversibility check is set to~$\eta_{\rm rev}=10^{-6}$, see~\eqref{eq:reversibiliy_check}. Note that in our case, we set the tolerance for the reversibility check in~\eqref{eq:reversibiliy_check} smaller than the tolerances for Newton's method in~\eqref{eq:newton_thresholds} in order to take into account the accumulating round-off errors when assessing numerical reversibility.  We checked that, in our numerical experiment, when~\eqref{eq:reversibiliy_check} does not hold, the value of the left-hand side is orders of magnitude larger than the tolerance used.
    
\end{appendices}

\end{document}